\documentclass[11pt]{amsart}
\usepackage{cite,mathrsfs,amsmath,amssymb,amsthm,bbm,esint,enumerate}
\usepackage[english]{babel}
\usepackage[T1]{fontenc}
\usepackage[utf8]{inputenc}

\usepackage[theoremfont,largesc]{newtxtext}
\usepackage[nopatch=eqnum]{microtype}

\usepackage[hang,flushmargin,bottom]{footmisc}
\usepackage[letterpaper,textwidth=5.7in,top=1.1in,bottom=1.1in,centering]{geometry}

\definecolor{oiBlue}{RGB}{0, 114, 178}

\usepackage[colorlinks=true,urlcolor=magenta,citecolor=oiBlue,linkcolor=magenta,breaklinks=true]{hyperref}
\usepackage{orcidlink}

\newtheorem{lem}{Lemma}[section]
\newtheorem{thrm}[lem]{Theorem}
\newtheorem{prop}[lem]{Proposition}
\newtheorem{cor}[lem]{Corollary}
\theoremstyle{definition}

\theoremstyle{remark}

\newcommand{\eq}[2]{\begin{equation}\label{#1}#2\end{equation}}
\newcommand{\pde}[2]{\begin{cases}#1\smallskip\\#2\end{cases}}

\newcommand{\R}{\mathbb{R}}
\newcommand{\C}{\mathbb{C}}
\newcommand{\Z}{\mathbb{Z}}

\newcommand{\Schwartz}{\mathcal{S}}
\newcommand{\Test}{C_c^\infty}

\renewcommand{\Re}{\operatorname{Re}}

\renewcommand{\epsilon}{\varepsilon}
\newcommand{\mc}{\mathcal}
\newcommand{\mb}{\mathbf}
\newcommand{\mr}{\mathrm}
\newcommand{\mf}{\mathfrak}

\DeclareMathOperator{\supp}{supp}
\DeclareMathOperator{\sgn}{sgn}
\newcommand{\<}{\langle}
\renewcommand{\>}{\rangle}
\newcommand{\p}{\partial}

\newcommand{\bpf}{\begin{proof}}
\newcommand{\epf}{\end{proof}}
\newcommand{\qtq}[1]{\quad\text{#1}\quad}

\newcommand{\qt}[1]{\quad\text{#1}}

\newcommand{\LHS}[1]{\text{LHS\eqref{#1}}}
\newcommand{\RHS}[1]{\text{RHS\eqref{#1}}}

\newcommand{\I}{\mf I}
\DeclareMathOperator{\tr}{tr}
\newcommand{\cP}{\mc P}
\newcommand{\cH}{\mc H}

\newcommand{\sprack}[1]{^{[#1]}}
\newcommand{\vk}{\varkappa}
\newcommand{\op}{\mr{op}}

\newcommand{\err}{\mb{err}}
\newcommand{\Err}{\mb{Err}}

\newcommand{\X}{\mb X}
\newcommand{\Y}{\mb Y}
\newcommand{\cG}{\mc G}
\newcommand{\bB}{\mb B}

\newcommand{\fB}{\mf B}
\newcommand{\bv}{\mb v}

\newcommand{\cR}{\mc R}
\newcommand{\cQ}{\mc Q}

\title{On the long wave limit of the good Boussinesq equation}

\author[B. Harrop-Griffiths]{Benjamin Harrop-Griffiths \orcidlink{0000-0002-0215-7613}}
\address{Department of Mathematics and Statistics, Georgetown University, DC 20057, USA}
\email{benjamin.harropgriffiths@georgetown.edu}
\thanks{}

\keywords{}
\subjclass{}

\numberwithin{equation}{section}
\allowdisplaybreaks

\begin{document}

\begin{abstract}
We consider the good Boussinesq equation on the line.  Sufficiently regular solutions are known to converge in the long wave limit to a pair of counterpropagating solutions to the Korteweg--de Vries equation.  We prove that this convergence holds for solutions in $H^{-\frac12}$, the largest Sobolev space in which the good Boussinesq equation is well-posed.  As part of our analysis, we revisit the well-posedness theory, proving new a priori estimates and a corresponding small data global well-posedness result.
\end{abstract}

\maketitle

\section{Introduction}

We consider solutions \(u\colon \R_t\times \R_x\to \R\) of the ``good'' Boussinesq equation,
\eq{GB}{\tag{gB}
u_{tt} - \tfrac43 u_{xx} + \tfrac13 u_{xxxx} + \tfrac43(u^2)_{xx} = 0.
}
The equation \eqref{GB} appeared as a model for ion waves in plasmas in \cite{Kako} and for solitary waves in shape-memory alloys in \cite{PhysRevB.36.3031}.

The good Boussinesq equation \eqref{GB} differs by a sign from the ``bad'' Boussinesq equation,
\eq{bB}{\tag{bB}
u_{tt} - \tfrac43 u_{xx} - \tfrac13 u_{xxxx} + \tfrac43(u^2)_{xx} = 0,
}
which first appeared in \cite{MR3363411} as an effective model for shallow water waves.  Unfortunately, \eqref{bB} is extraordinarily badly behaved: even the linearization of \eqref{bB} at zero is ill-posed! Indeed, it was shown in \cite{charlier2023boussinesqsequationwaterwaves,MR5130833} that there exist Schwartz solutions of \eqref{bB} that blow up on arbitrarily short timescales.

Our primary interest in this article is the long wave limit of \eqref{GB}.  Introducing a small parameter \(0<\epsilon\ll1\), we consider solutions of the form
\eq{u-rescaled}{
u(t,x) = \epsilon^2 u^\epsilon(\epsilon^3 t,\epsilon x),
}
which yields the rescaled equation
\eq{GB-rescaled}{\tag{gB\textsuperscript{$\epsilon$}}
\epsilon^2u_{tt}^\epsilon - \tfrac 4{3\epsilon^2}u_{xx}^\epsilon + \tfrac13u_{xxxx}^\epsilon + \tfrac43\bigl((u^\epsilon)^2\bigr)_{xx} = 0.
}
A formal multiple scales expansion shows that to leading order
\eq{formal LW}{
u^\epsilon \approx w^{\epsilon,+} + w^{\epsilon,-},
}
where \(w^{\epsilon,\pm}\) are solutions of the Korteweg--de Vries (KdV) equations,
\eq{KdV}{\tag{KdV$^{\epsilon,\pm}$}
w_t^{\epsilon,\pm} \pm \tfrac2{\sqrt 3\epsilon^2}w_x^{\epsilon,\pm}  = \pm \tfrac1{4\sqrt 3}\bigl(w_{xxx}^{\epsilon,\pm} + 8w^{\epsilon,\pm}w_x^{\epsilon,\pm}\bigr).
}
Note that if \(w^{\epsilon,\pm}\) satisfies \eqref{KdV} then
\eq{static rescaling}{
w(t,x) = w^{\epsilon,\pm}\bigl(\pm t,x + \tfrac2{\sqrt 3\epsilon^2}t\bigr)
}
is a solution of the ``usual'' KdV equation,
\eq{KdV stat}{\tag{KdV}
w_t = \tfrac1{4\sqrt 3}\bigl(w_{xxx} + 8ww_x\bigr).
}

In recent work, Hong and Yang \cite{MR4767456} and then Hong and Jang \cite{hong2026longtimevaliditykortewegdevries} proved the validity of the approximation \eqref{formal LW} for initial data in \(H^{0+}\). This improved previous high regularity results for unidirectional waves given in \cite{MR3922812,MR3702025}. In previous work, Schneider \cite{MR1628680} proved a similar approximation holds for a well-posed modification of \eqref{bB} with initial data in weighted Sobolev spaces. This was subsequently refined in \cite{MR1968371}.

The first contribution of this article is that the long wave approximation \eqref{formal LW} holds at optimal regularity:
\begin{thrm}[The long wave limit]\label{t:long wave}
Let \(s\geq -\frac12\).
\begin{enumerate}[(i)]
\item If \((u_0,u_1)\in H^s\times \p_x H^{s-1}\) and for all \(0<\epsilon\ll1\) sufficiently small \((u^\epsilon,u_t^\epsilon)\in C(\R;H^s\times \p_x H^{s-1})\) is the solution of \eqref{GB-rescaled} with initial data \((u^\epsilon,u_t^\epsilon)(0) = (u_0,u_1)\) then there exist solutions \(w^{\epsilon,\pm} \in C(\R;H^s)\) of \eqref{KdV} so that for all \(T>0\) we have
\eq{long wave approx}{
\lim_{\epsilon \to 0}\sup_{|t|\leq T}\|u^\epsilon(t) - w^{\epsilon,+}(t) - w^{\epsilon,-}(t)\|_{H^s} = 0.
}
\item Conversely, if \(w_0^\pm\in H^s\) and \(w^{\epsilon,\pm} \in C(\R;H^s)\) are the solutions of \eqref{KdV} with initial data \(w^{\epsilon,\pm}(0) = w_0^\pm\) then for all \(0<\epsilon \ll1\) sufficiently small there exists a solution \((u^\epsilon,u_t^\epsilon)\in C(\R;H^s\times\p_x H^{s-1})\) of \eqref{GB-rescaled} satisfying \eqref{long wave approx} for all \(T>0\).
\end{enumerate}
\end{thrm}

Theorem~\ref{t:long wave} extends the results of \cite{MR4767456,hong2026longtimevaliditykortewegdevries} in two key ways.  First, we are able to lower the regularity of the initial data to \(H^{-\frac12}\).  As we will discuss next, this is optimal in light of the well-posedness theory.  Second, we prove convergence in the same topology as the initial data.  All of the aforementioned long wave results, \cite{MR4767456,hong2026longtimevaliditykortewegdevries,MR3922812,MR3702025,MR1628680,MR1968371}, involve a loss of regularity between the data and the space in which convergence is obtained.

As we will discuss in more detail below, the method we use to prove Theorem~\ref{t:long wave} is quite robust.  The argument combines uniform (in \(\epsilon\)) low regularity well-posedness for \eqref{GB-rescaled} and \eqref{KdV} with an energy estimate for smooth solutions.  In particular, we expect this approach could be adapted to obtain similarly sharp long wave approximation results for a variety of other models.

Theorem~\ref{t:long wave} proves that the approximation \eqref{formal LW} is valid in Sobolev spaces on the scale of the rescaled equation \eqref{GB-rescaled}. Rescaling to the original variables according to \eqref{u-rescaled}, the identity \eqref{long wave approx} implies that
\[
\sup_{|t|\leq \frac T{\epsilon^3}}\|u(t,x) - \epsilon^2 w^{\epsilon,+}(\epsilon^3t,\epsilon x) - \epsilon^2 w^{\epsilon,-}(\epsilon^3t,\epsilon x)\|_{H^s} = o\bigl(\epsilon^{\frac32+\min\{s,0\}}\bigr),
\]
as \(\epsilon\to0\).

Let us now justify our assertion that Theorem~\ref{t:long wave} constitutes the ``optimal regularity'' for \eqref{GB}.  The well-posedness of \eqref{GB} has been the subject of a large body of work, including \cite{MR954673,MR1251854,MR2512853,MR2654246,MR3032295,MR2566562,MR1399198,MR1095753,MR3953031,MR2756068}.  The optimal local well-posedness result on the line (and the circle) was obtained in \cite{MR3016208}, where Kishimoto proved that \eqref{GB} is locally well-posed for initial data \((u,u_t)(0)\in H^s\times H^{s-2}\) for all \(s\geq -\frac12\).

Kishimoto's result is optimal in two regards.   First, he also proved in \cite{MR3016208} that for \(s<-\frac12\) the solution map for \eqref{GB} fails to be continuous at zero.  This demonstrates that there is no hope of local well-posedness in larger Sobolev spaces.  Second, Sachs \cite{MR1048956} implemented the argument of \cite{MR604036} to show that there exists initial data for which \eqref{GB} blows up in finite time.  This proves that we should not expect to have global well-posedness for arbitrary initial data.

By writing \((u,u_t) = (q,r_x)\), we may formally recast \eqref{GB} as a system,
\eq{GB-system}{
\pde{
q_t = r_x,
}{
r_t = \bigl(\tfrac 43 q - \tfrac 1 3 q_{xx} - \tfrac43 q^2\bigr)_x.
}
}
Every solution \((q,r)\in C(\R;H^s\times H^{s-1})\) of \eqref{GB-system} yields a solution \((u,u_t)\in C(\R;H^s\times H^{s-2})\) of \eqref{GB}; the converse is false.  The system \eqref{GB-system} was proved to be locally well-posed in \(H^s\times H^{s-1}\) for all \(s>\frac52\) in \cite{MR954673} and then all \(s\geq0\) in \cite{MR1251854}.

The system \eqref{GB-system} is the Hamiltonian flow of
\eq{energy}{
\cH[q,r] := \int_\R \tfrac23 q^2 + \tfrac16 q_x^2 +  \tfrac12r^2 - \tfrac49q^3\,dx,
}
with respect to the Poisson bracket
\eq{Poisson bracket}{
\{\cG,\cH\} := \int_\R\frac{\delta \cG}{\delta q}\left(\frac{\delta \cH}{\delta r}\right)_x + \frac{\delta \cG}{\delta r}\left(\frac{\delta \cH}{\delta q}\right)_x\,dx.
}
The Hamiltonian \(\cH\) controls \(\|(q,r)\|_{H^1\times L^2}^2\) for sufficiently small solutions and is conserved under the flow.  This was used by Linares in \cite{MR1251854} to prove that \eqref{GB-system} is globally well-posed for sufficiently small initial data in \(H^1\times L^2\).

As the reader will have no doubt surmised from this discussion, there is a gap between the existing well-posedness theory for \eqref{GB} and our implicit assertion in Theorem~\ref{t:long wave} that \eqref{GB-rescaled} is globally well-posed in \(H^s\times \p_x H^{s-1}\) for sufficiently small \(0<\epsilon\ll1\).  Our next result will close this gap:
\begin{thrm}[Global well-posedness]\label{t:GWP}
There exists \(\delta>0\) so that for all \(s\geq -\frac12\) and \((q,r)(0)\in H^s\times H^{s-1}\) satisfying
\eq{init small}{
\|(q,r)(0)\|_{H^{-\frac12}\times H^{-\frac32}}<\delta,
}
the equation \eqref{GB-system} is globally well-posed in \(H^s\times H^{s-1}\).
\end{thrm}

While preparing this article for submission, we learned that an alternative proof of Theorem~\ref{t:GWP} using the method of commuting flows has been independently obtained by M.~Shan~\cite{MSPC}. We also learned that Shan had previously announced similar a priori estimates to those given in Theorem~\ref{t:AP} in \cite{MStalk}.

At this juncture, let us briefly clarify our notion of well-posedness.  It follows from \cite{MR954673} that for all \(s>\frac52\) the system \eqref{GB-system} is unconditionally locally well-posed in \(H^s\times H^{s-1}\). Precisely, for all initial data \((q_0,r_0)\in H^s\times H^{s-1}\) we may find a maximal open interval \(I\) containing zero and a unique distributional solution \((q,r)\in C(I;H^s\times H^{s-1})\) of \eqref{GB-system} with initial data \((q,r)(0) = (q_0,r_0)\).  Moreover, for any compact time interval \(J\subseteq I\) containing zero there exists a neighborhood \(N\subseteq H^s\times H^{s-1}\) of \((q_0,r_0)\) so that the solution map \((t,q_0,r_0)\mapsto \Phi_t(q_0,r_0) := (q,r)(t)\) is continuous from \(J\times N\to H^s\times H^{s-1}\).

Denoting the ball
\eq{balls}{
B_\delta := \bigl\{(q,r)\in H^{-\frac12}\times H^{-\frac32}:\|(q,r)\|_{H^{-\frac12}\times H^{-\frac32}}<\delta\bigr\},
}
for \(s>\frac52\) we will prove that if \((q,r)(0)\in (H^s\times H^{s-1})\cap B_\delta\) then the maximal time of existence \(I = \R\).
When \(-\frac12\leq s\leq \frac52\), we will instead prove that the solution map \(\Phi\) has a unique continuous extension to a map \(\R\times \bigl((H^s\times H^{s-1})\cap B_\delta\bigr)\to H^s\times H^{s-1}\).  Our proof will also show that for all \(s\geq -\frac12\) and \(T>0\) the map \((H^s\times H^{s-1})\cap B_\delta\ni(q_0,r_0)\mapsto \Phi_t(q_0,r_0)\in H^s\times H^{s-1}\) is locally Lipschitz continuous, uniformly for all \(|t|\leq T\).

As solutions of \eqref{GB-system} yield solutions of \eqref{GB}, Theorem~\ref{t:GWP} immediately yields global well-posedness (in the same sense) for \eqref{GB}:
\begin{cor}\label{c:GWP}
There exists \(\delta>0\) so that for all \(s\geq -\frac12\) and \((u,u_t)(0)\in H^s\times \p_x H^{s-1}\) satisfying
\eq{init small-2}{
\|(u,u_t)(0)\|_{H^{-\frac12}\times \p_x H^{-\frac32}}<\delta,
}
the equation \eqref{GB} is globally well-posed in \(H^s\times \p_x H^{s-1}\).
\end{cor}

It is straightforward to verify that if \(s\geq -\frac12\) and \(0<\epsilon\leq 1\) then for any \((u_0,u_1)\in H^s\times \p_xH^{s-1}\) we have
\[
\|(\epsilon^2u_0,\epsilon^5u_1)(\epsilon x)\|_{H^{-\frac12}\times \p_xH^{-\frac32}}\lesssim \epsilon \|(u_0,u_1)\|_{H^{-\frac12}\times \p_xH^{-\frac32}}.
\]
This immediately yields a second corollary that resolves the question of existence of global solutions to \eqref{GB-rescaled} that arises in Theorem~\ref{t:long wave}:
\begin{cor}\label{cor:rescaled well-posedness}
For all \(s\geq -\frac12\) and \((u_0,u_1)\in H^s\times \p_x H^{s-1}\) there exists \(\epsilon_0>0\), depending only on \(\|(u_0,u_1)\|_{H^{-\frac12}\times \p_xH^{-\frac32}}\), so that \eqref{GB-rescaled} with initial data \((u^\epsilon,u_t^\epsilon)(0) = (u_0,u_1)\) is globally well-posed in \(H^s\times \p_xH^{s-1}\) for all \(0<\epsilon\leq \epsilon_0\).
\end{cor}

Let us now turn our discussion to the proof of Theorems~\ref{t:long wave} and~\ref{t:GWP}.  Our starting point will be to consider the local well-posedness of \eqref{GB-system}.  As we have already mentioned, the result of Kishimoto~\cite{MR3016208} does not directly imply local well-posedness for the system \eqref{GB-system} in \(H^s\times H^{s-1}\).  Consequently, we must first prove a new local well-posedness result:
\begin{thrm}[Local well-posedness]\label{t:LWP} For all \(R>0\) we may find a time \(T\simeq \<R\>^{-6}\) so that for all \(s\geq-\frac12\) and initial data \((q,r)(0)\in (H^s\times H^{s-1})\cap B_R\) the equation \eqref{GB-system} is locally well-posed in \(H^s\times H^{s-1}\) on the time interval \([-T,T]\).
\end{thrm}

We prove Theorem~\ref{t:LWP} in Section~\ref{s:LWP} using a contraction mapping argument in similar function spaces to those adopted by Kishimoto in \cite{MR3016208}, which were in turn inspired by those appearing in \cite{MR2204680}.  In fact, we will prove a more refined local well-posedness result in Theorem~\ref{t:LWP-details}; we will return to this point shortly.

With local well-posedness in hand, Theorem~\ref{t:GWP} now follows from appropriate a priori estimates:
\begin{thrm}[A priori estimates]	\label{t:AP}
There exists \(\delta>0\) so that if \(I\subseteq \R\) is an open interval containing zero and \((q,r)\in C(I;H^3\times H^2)\) is a solution of \eqref{GB-system} satisfying \eqref{init small} then for all \(-\frac12\leq s<2\) and \(t\in I\) we have the estimate
\eq{AP bound}{
\|(q,r)(t)\|_{H^s\times H^{s-1}}\simeq_s \|(q,r)(0)\|_{H^s\times H^{s-1}}.
}
\end{thrm}

The restriction \(-\frac12\leq s<2\) in Theorem~\ref{t:AP} (and Theorem~\ref{t:equicontinuity} below) is almost certainly illusory.  We believe it should be possible to use similar arguments to \cite{MR3874652,MR4747862} to obtain a corresponding result for all \(s\geq -\frac12\).

The proof of Theorem~\ref{t:AP} will rely on a remarkable fact about the system \eqref{GB-system}: it is completely integrable! In particular, it possesses a Lax Pair: unbounded operators
\eq{good Lax pair}{
L := \p_x^3 - \p_x + q\p_x + \p_x q + r\qtq{and}P := \p_x^2 + \tfrac 43 q
}
on \(L^2\) so that \((q,r)\) solve \eqref{GB-system} if and only if
\[
L_t = [P,L].
\]
The Lax pair \eqref{good Lax pair} can be readily derived from the corresponding Lax pair for \eqref{bB} appearing in \cite{zakharov};  it appears explicitly in \cite{MR730182}.

The existence of a Lax pair frequently yields a corresponding inverse scattering transform that can be used to solve the underlying PDE.  Unfortunately, as pointed out in \cite{MR2411922}, developing a comprehensive inverse scattering transform for Lax operators associated to Boussinesq equations is an extremely challenging problem. 

If we take \(q = p + \frac12\) in \eqref{GB-system} then we obtain the system
\eq{GB-system for scat}{
\pde{
p_t = r_x,
}{
r_t = - \bigl(\tfrac 1 3 p_{xx} + \tfrac43 p^2\bigr)_x.
}
}
An inverse scattering approach to \eqref{GB-system for scat} was developed in \cite{MR4481092}. This was used to investigate the asymptotic behavior of generic, Schwartz, solitonless initial data in \cite{MR4632376}, with further results obtained in \cite{wang2025longtimeasymptoticsgoodboussinesq}.  Let us also pause to mention the article \cite{MR4692724} that demonstrates a remarkable nonlinear ``Miura map'' taking solutions of the quadratic derivative nonlinear Schr\"odinger equation,
\eq{quad NLS}{
i\psi_t = \psi_{xx} - \bar\psi\bar\psi_x,
}
to solutions of \eqref{GB-system for scat}. The equation \eqref{quad NLS} also possesses a Lax pair, which was used to obtain a priori estimates in \cite{klaus}.

Our a priori estimates for \eqref{GB-system} are based on the conservation of the (negative) logarithm of the renormalized perturbation determinant, (formally) defined by
\eq{perturbation determinant formal}{
\alpha(\kappa;q,r) := -\log\det\!{}_2\left[1 + i\tfrac1{\sqrt{iL_0 + i\kappa^3-i\kappa}}(L - L_0)\tfrac1{\sqrt{iL_0 + i\kappa^3-i\kappa}}\right],
}
for real \(|\kappa|\geq 2\), where the free Lax operator
\eq{free Lax}{
L_0 := \p_x^3 - \p_x.
}

The use of the logarithm of the (renormalized) perturbation determinant (or equivalently transmission coefficient) to prove optimal a priori bounds in Sobolev spaces was pioneered in the works \cite{MR3820439,MR3874652}.  Subsequently, multiple articles have adopted a related approach to obtain a priori estimates for a plethora of completely integrable PDEs; see, e.g., \cite{MR4515528,MR4628747,MR4923538,MR4270277,klaus,MR4743514,MR5083420}.

In Section~\ref{s:AP}, we will prove that if we integrate \(\alpha(\kappa) + \alpha(-\kappa)\) in \(\kappa\) against appropriate weights, we recover conserved quantities that control \(\|(q,r)\|_{H^s\times H^{s-1}}^2\) for \(-\frac12\leq s<1\).  This yields the corresponding a priori estimates of Theorem~\ref{t:AP} for this range of \(s\) and consequently proves Theorem~\ref{t:GWP}. The \(s=1\) case of \eqref{AP bound} is obtained by combining the conservation of \(\cH\) with the a priori estimates for \(-\frac12\leq s<1\). The remaining a priori estimates are obtained in a similar manner to the \(-\frac12\leq s<1\) case, after correcting \(\alpha(\kappa) + \alpha(-\kappa)\) by a \(\kappa\)-dependent multiple of \(\cH\).

When proving a priori estimates using the Lax pair structure, it is frequently possible to prove that the set of orbits emanating from a bounded and equicontinuous set of initial data is bounded and equicontinuous.  This crucial observation underpins Killip and Vi\c san's method of commuting flows \cite{MR3990604}, and has led to optimal well-posedness results for a host of completely integrable PDEs; see, e.g., \cite{MR4304314,MR4628747,MR4726498,MR4380591,MR4743514,MR4922705,MR5031430,klaus2026wellposednesskdvhierarchy,killip2026hamiltonianformulationcontinuumcalogeromoser}.  Our final result proves an analogous result for \eqref{GB-system}:

\begin{thrm}\label{t:equicontinuity}
If \(-\frac12\leq s<2\) and \(\cQ\subseteq (H^s\times H^{s-1})\cap B_\delta\) is a bounded and equicontinuous set of initial data then the set
\[
\cQ^* = \bigl\{\Phi_t(q,r)\in H^s\times H^{s-1}:t\in \R\text{ and }(q,r)\in \cQ\bigr\}
\]
is bounded and equicontinuous.
\end{thrm}

An inspection of the proof of Theorem~\ref{t:GWP} shows that for all \(T>0\) and \(|t|\leq T\) the map \((q,r)\mapsto \Phi_t(q,r)\) is uniformly Lipschitz continuous for \((q,r)\) in a bounded subset of \((H^s\times H^{s-1})\cap B_\delta\). Combining this with the fact that the solution map commutes with spatial translations (or equivalently that smooth solutions conserve the momentum \(\cP[q,r] := \int_\R qr\,dx\)) suffices to ensure that the set
\[
\cQ_T^* = \bigl\{\Phi_t(q,r)\in H^s\times H^{s-1}:|t|\leq T\text{ and }(q,r)\in \cQ\bigr\}
\]
is bounded and equicontinuous.  The key novelty of Theorem~\ref{t:equicontinuity} is that equicontinuity holds for all \(t\in \R\).

With the question of well-posedness settled, we finally turn to the proof of Theorem~\ref{t:long wave} in Section~\ref{s:long wave}.  If we rescale \eqref{GB-system} in a manner consistent with \eqref{u-rescaled} by taking
\[
(q,r)(t,x) = (\epsilon^2q^\epsilon,\epsilon^4r^\epsilon)(\epsilon^3 t,\epsilon x),
\]
we obtain the system
\eq{GB-system rescaled}{
\pde{
q_t^\epsilon = r_x^\epsilon,
}{
\epsilon^2r_t^\epsilon = \bigl(\tfrac 4{3\epsilon^2} q^\epsilon - \tfrac13 q_{xx}^\epsilon - \tfrac43 (q^\epsilon)^2\bigr)_x.
}
}
Rather than working directly with \eqref{GB-system rescaled}, in a similar fashion to \cite{MR4767456,hong2026longtimevaliditykortewegdevries} we symmetrize by introducing new dependent variables,
\eq{d:veps}{
v^{\epsilon,\pm} := \tfrac12\bigl(q^\epsilon \mp \epsilon^2\Lambda_\epsilon^{-1}r^\epsilon\bigr),
}
where
\eq{d:lambda}{
\Lambda_\epsilon := \sqrt{\tfrac 43 - \tfrac13\epsilon^2\p_x^2},
}
so that \eqref{GB-system rescaled} becomes
\eq{GB-system-2}{
v_t^{\epsilon,\pm} = \mp\tfrac1{\epsilon^2}\Lambda_\epsilon v_x^{\epsilon,\pm} \pm \tfrac23\Lambda_\epsilon^{-1}\p_x\bigl[(v^{\epsilon,+}+v^{\epsilon,-})^2\bigr].
}

In Theorem~\ref{t:LWP-details} we prove that the equation \eqref{GB-system-2} is locally well-posed in \(H^s\) for all \(s\geq -\frac12\), where the timescale can be chosen to be uniform for all sufficiently small (relative to the size of the initial data) \(\epsilon\).  This is our aforementioned refinement of Theorem~\ref{t:LWP}.
The fact that we can prove uniform-in-\(\epsilon\) estimates for solutions should not be a surprise: if we believe the long wave limit is valid then uniform bounds must hold!  Nevertheless, proving this is quite challenging.  We are required to revisit the local theory, keeping careful track of the powers of \(\epsilon\) throughout.  This is by far the most technical part of this article, culminating in the delicate nonlinear estimates of Proposition~\ref{p:bilinear estimates}.

Our proof of local well-posedness in Theorem~\ref{t:LWP-details} shows that, for sufficiently small \(0<\epsilon\ll1\), the map \((v^{\epsilon,+},v^{\epsilon,-})(0)\mapsto (v^{\epsilon,+},v^{\epsilon,-})(t)\) is locally Lipschitz in \(H^s\times H^s\) with a Lipschitz constant that is uniform in \(\epsilon\) on the (similarly uniform) timescale on which local well-posedness holds.  This uniform continuity of the solution map lies at the heart of the proof of Theorem~\ref{t:long wave} --- it reduces the problem to showing that the long wave limit holds for some approximating sequence of initial data.  An analogous uniform-in-\(\epsilon\) continuity result for solutions of \eqref{KdV} is also true.  This follows from the continuity of the solution map for \eqref{KdV stat} and the change of variables \eqref{static rescaling}.

To complete the proof, it suffices to show that \eqref{long wave approx} is true in \(H^s\) for initial data in \(H^m\) with \(m\gg s\). (We take \(m = \lceil s+5\rceil\).) Here, we approximate \(\Lambda_\epsilon\) by the leading terms in its Taylor expansion in \eqref{GB-system-2}, utilizing the additional regularity to control the remainder.  We then apply Gronwall's inequality to obtain an energy estimate that bounds the difference between solutions of \eqref{GB-system-2} and solutions of \eqref{KdV}.

A key challenge in both the energy estimate and the final approximation argument is controlling the size of solutions to \eqref{GB-rescaled} (or equivalently \eqref{GB-system rescaled} or \eqref{GB-system-2}) uniformly in \(\epsilon\).  The a priori estimates of Theorem~\ref{t:AP} are not scale invariant: they are adapted to the scale \(\epsilon=1\).  In particular, they do not yield uniform-in-\(\epsilon\) bounds.  Instead, we use that the a priori estimates for \eqref{KdV stat} proved in \cite{MR3820439,MR3874652} yield uniform control of \eqref{KdV}.  These uniform estimates can then be transferred to the solution of \eqref{GB-rescaled} by bootstrapping the identity \eqref{long wave approx}.

The remainder of this paper comprises three, relatively self-contained sections. Section~\ref{s:LWP} addresses the question of local well-posedness, proving both Theorem~\ref{t:LWP} and its generalization, Theorem~\ref{t:LWP-details}. Section~\ref{s:AP} provides the a priori estimates of Theorems~\ref{t:AP} and~\ref{t:equicontinuity}. Combining Theorem~\ref{t:AP} with local well-posedness then yields Theorem~\ref{t:GWP}. Finally, Section~\ref{s:long wave} addresses the central question in this article, Theorem~\ref{t:long wave}, essentially using the results of Sections~\ref{s:LWP} and~\ref{s:AP} as a black box.

\subsection*{Acknowledgments \& statement of AI usage} This material is based upon work supported by the National Science Foundation under Award No.~DMS-2406816.  Part of this work was completed while the author participated in the program \emph{Dispersive Integrable Equations: Pathfinders in Hamiltonian PDE} at the Institut Henri Poincar\'e.  The author is grateful to the organizers for their invitation and the IHP for their hospitality.

The author acknowledges the use of AI and computer algebra tools in preparing this manuscript.  Several Google Gemini models were used to assist with preliminary investigations, bookwork, and literature review.  The polynomials \(b\), \(c\), \(d\) in Proposition~\ref{p:alpha} were computed with help from Wolfram Mathematica.  The first draft of this article was written entirely by the author.  Proofreading was then carried out by the author in parallel with OpenAI's GPT-6 Astra. This led to several minor revisions, a simplification in the proof of Lemma~\ref{l:modulation}, and our discovery of the announcement of \cite{MStalk}.  All other mathematical arguments and all text in this article were written solely by the author, who takes full responsibility for its content.

\section{Local well-posedness}\label{s:LWP}

We start this section by introducing some notation that will be used throughout this article. Given a Schwartz function \(f\in \Schwartz(\R)\), we denote its Fourier transform by
\[
\hat f(\xi) := \tfrac1{\sqrt{2\pi}}\int_\R f(x)e^{-ix\xi}\,dx.
\]
If \(f\in \Schwartz(\R^2)\) is a Schwartz function of space and time, we denote its spacetime Fourier transform by
\[
\widetilde f(\tau,\xi) := \tfrac1{2\pi}\iint_{\R^2}f(t,x)e^{-i(t\tau+x\xi)}\,dx\,dt.
\]

For real-valued \(|\kappa|\geq 1\) and \(s\in \R\), we take
\[
\|f\|_{H_\kappa^s}^2 := \int_\R (\kappa^2 + \xi^2)^s|\hat f(\xi)|^2\,d\xi,
\]
so that \(H^s = H_1^s\).  We let \(\p_xH^s:=\{f_x:f\in H^s\}\) denote the space of derivatives of \(H^s\) functions, which we endow with the norm \(\|f_x\|_{\p_xH^s} = \|f\|_{H^s}\). Given Banach spaces \(X\) and \(Y\) with corresponding norms \(\|\cdot\|_X\) and \(\|\cdot \|_Y\), we set
\[
\|(f,g)\|_{X\times Y}^2:=\|f\|_X^2 + \|g\|_Y^2.
\]

We take \(\chi\in \Test(\R)\) to be a fixed even function satisfying \(0\leq \chi(x)\leq 1\) that is identically \(1\) on \([-1,1]\) with compact support in \((-2,2)\). For any constant \(R>0\) we denote the rescaled function \(\chi_{\leq R}(x) := \chi\bigl(\frac xR\bigr)\) and let \(\chi_{>R} := 1 - \chi_{\leq R}\).

For dyadic \(N\in 2^\Z\) we define the spatial Littlewood--Paley projection
\[
\widehat{P_{\leq N} f}(\xi) := \chi_{\leq N}(\xi)\hat f(\xi),
\]
and set \(P_{>N} = 1 - P_{\leq N}\), \(P_N := P_{\leq N}P_{>\frac N2}\).  We will require both inhomogeneous and homogeneous Littlewood--Paley decompositions in this article.  However, as we will primarily work with inhomogeneous decompositions, unless otherwise stated we denote \(f_1 = P_{\leq 1}f\) and \(f_N = P_N f\) for dyadic \(N\geq 2\), so that
\[
f = \sum_{N\geq 1}f_N.
\]
For all dyadic \(N\in 2^\Z\) we denote \(f_{\leq N} = P_{\leq N}f\) and \(f_{>N} = P_{>N}f\).

For \(f\in C([-T,T];L^2)\) we take
\[
\|f\|_{C_TL_x^2} := \sup_{|t|\leq T}\|f(t)\|_{L^2},
\]
and for bounded \(f\in C(\R;L^2)\) we denote
\[
\|f\|_{C_tL_x^2} := \sup_{t\in \R}\|f(t)\|_{L^2}.
\]

Finally, when considering solutions of \eqref{GB-system-2}, we will frequently take \(v^\epsilon := (v^{\epsilon,+},v^{\epsilon,-})\), with a similar convention for other pairs of equations with right/left-moving components, such as \eqref{KdV}.

With our notation clarified, let us turn to the main goal of this section: proving local well-posedness of the system \eqref{GB-system-2}.  Theorem~\ref{t:LWP} will then follow directly from the case \(\epsilon = 1\) after we note that if \(v^{\epsilon,\pm}\) are defined as in \eqref{d:veps} then the map
\eq{translate between systems}{
H^s\times H^{s-1}\ni (q,r)\mapsto (v^{1,+},v^{1,-})\in H^s\times H^s
}
is linear, invertible, and satisfies
\[
\|(q,r)\|_{H^s\times H^{s-1}} \simeq \|(v^{1,+},v^{1,-})\|_{H^s\times H^s}.
\]

To state our local well-posedness result, let us introduce the ball
\eq{fBR def}{
\fB_R := \Bigl\{(v^+,v^-)\in H^{-\frac12}\times H^{-\frac12}:\|(v^+,v^-)\|_{H^{-\frac12}\times H^{-\frac12}}\leq R\Bigr\}.
}
We then have the following:

\begin{thrm}\label{t:LWP-details}
If \(R>0\) then for all \(\epsilon>0\) there exists a time \(T_\epsilon\simeq \<R\<R\epsilon\>\>^{-3}\) so that for all \(s\geq -\frac12\) and all \(v_0 = (v_0^+,v_0^-)\in (H^s\times H^s)\cap \fB_R\) the system \eqref{GB-system-2} with initial data \(v^\epsilon(0) = (v^{\epsilon,+},v^{\epsilon,-})(0) = v_0\) is locally well-posed in \(H^s\times H^s\) on the time interval \([-T_\epsilon,T_\epsilon]\).

Further, if we denote the solution map by \(\Upsilon_t^\epsilon(v_0):=v^\epsilon(t)\) then for all \(v_0,w_0\in (H^s\times H^s)\cap \fB_R\) there exists a time \(T\simeq \<R\>^{-3}\) so that
\begin{align}
&\sup_{0<\epsilon\leq \frac1R}\sup_{|t|\leq T}\|\Upsilon_t^\epsilon(v_0) - \Upsilon_t^\epsilon(w_0)\|_{H^s\times H^s}\label{uniform lipschitz}\\
&\qquad\lesssim_{s,R}\|v_0-w_0\|_{H^s\times H^s} + \Bigl[\|v_0\|_{H^s\times H^s} + \|w_0\|_{H^s\times H^s}\Bigr]\|v_0 - w_0\|_{H^{-\frac12}\times H^{-\frac12}}.\notag
\end{align}
\end{thrm}

A key facet of this result is that the estimate \eqref{uniform lipschitz} holds uniformly in \(\epsilon\).  As discussed in the introduction, this will be a vital ingredient for the proof of Theorem~\ref{t:long wave}. Unfortunately, proving this uniformity will cause several technical headaches, especially in the proof of our main multilinear estimate in Proposition~\ref{p:bilinear estimates}.

To prove Theorem~\ref{t:LWP-details}, we start by using the Duhamel formula to write \eqref{GB-system-2} with initial data \((v^{\epsilon,+},v^{\epsilon,-})(0) = (v_0^+,v_0^-)\) in the form
\eq{Duhamel}{
v^{\epsilon,\sigma}(t) = S^{\epsilon,\sigma}(t)v_0^\sigma + \tfrac{2\sigma}3\int_0^t S^{\epsilon,\sigma}(t - \tau)\Lambda_\epsilon^{-1}\p_x\bigl[\bigl(v^{\epsilon,+}(\tau)+v^{\epsilon,-}(\tau)\bigr)^2\bigr]\,d\tau,
}
where \(\sigma\in\{\pm\}\) and the linear propagator \(S^{\epsilon,\sigma}(t) := e^{-\frac\sigma{\epsilon^2}t\Lambda_\epsilon\p_x}\).  We then apply the contraction mapping theorem to this formulation of \eqref{GB-system-2} in suitable function spaces.

A natural first attempt at constructing appropriate function spaces would be to consider the corresponding \(X^{s,b}\)-spaces, first introduced in the context of well-posedness by Bourgain in \cite{MR1209299,MR1215780}.  Indeed, as the equation \eqref{GB-system-2} behaves like the KdV equation at low frequency, we might hope to control low frequency interactions as in the proof of local well-posedness for the KdV equation in \cite{MR1329387}.  Unfortunately, as observed in \cite{MR2512853}, it is impossible to use \(X^{s,b}\) spaces to handle the high frequency interactions when \(s<-\frac14\).  Instead, we follow Kishimoto \cite{MR3016208} in adopting similar spaces to those used in the study of the quadratic nonlinear Schr\"odinger equation by Bejenaru--Tao in \cite{MR2204680}.  These spaces trade modulation for regularity away from the characteristic curve and suffice to control the troublesome high \(\times\) high \(\to\) low interactions.

In contrast with the results of \cite{MR3016208}, we wish to keep careful track of the various powers of \(\epsilon\) to ensure uniformity.  This makes it prudent to work directly with the formulation \eqref{GB-system-2} rather than transforming the equation into a quadratic nonlinear Schr\"odinger equation.  We also modify the function spaces slightly, eschewing \(\ell^1\)-type summation in modulation in favor of lower powers of the modulation weight \(b\).  Finally, to enable us to work with spaces adapted to the unit timescale, we allow any \(\epsilon>0\) and do not restrict our attention to \(0<\epsilon\leq1\).

For \(b\in \R\), \(\epsilon>0\), and \(\sigma\in \{\pm\}\), we define
\begin{align*}
\|u\|_{X_{\epsilon,\sigma}^b} &:= \|\<\tau+\sigma\phi_\epsilon(\xi)\>^b\widetilde u\|_{L_{\xi,\tau}^2},\\
\|u\|_{Y_{\epsilon,\sigma}^b} &:= \|\<\tau+\sigma\phi_\epsilon(\xi)\>^b\widetilde u\|_{L_\xi^2L_\tau^1},\\
\|u\|_{Z_{\epsilon,\sigma}^b} &:= \|\<\tau+\sigma\phi_\epsilon(\xi)\>^b\widetilde u\|_{L_\xi^2L_\tau^\infty},
\end{align*}
where we denote the symbol of the operator \(\frac1{i\epsilon^2}\Lambda_\epsilon\p_x\) by
\eq{phi epsilon}{
\phi_\epsilon(\xi) := \tfrac1{\epsilon^3}\phi(\epsilon \xi),\qtq{where}\phi(\xi) :=  \tfrac1{\sqrt 3}\xi\sqrt{4 + \xi^2}.
}
We then take \(\X_{\epsilon,\sigma}^{s,b}\) to be the completion of \(\Schwartz(\R^2)\) under the norm
\begin{align*}
\|u\|_{\X_{\epsilon,\sigma}^{s,b}}^2:= \sum_{N\geq 1}N^{2s}\inf\Bigl\{\|v\|_{X_{\epsilon,\sigma}^b}^2 + \Bigl[\|w\|_{Y_{\epsilon,\sigma}^0}^2 + \<\phi_\epsilon(N)\>^{2b}\|w\|_{X_{\epsilon,\sigma}^0}^2\Bigr]:u_N = v + w\Bigr\}.
\end{align*}
We also introduce the space \(\Y_{\epsilon,\sigma}^{s,b} = \bigl(\<\p_t + \frac\sigma{\epsilon^2}\Lambda_\epsilon\p_x\>\X_{\epsilon,\sigma}^{s,b}\bigr)^*\), with norm
\begin{align*}
\|u\|_{\Y_{\epsilon,\sigma}^{s,b}}^2:= \sum_{N\geq 1}N^{-2s}\Bigl[\|u_N\|_{X_{\epsilon,\sigma}^{1-b}}^2 + \inf\Bigl\{\|v\|_{Z_{\epsilon,\sigma}^1}^2 + \<\phi_\epsilon(N)\>^{-2b}\|w\|_{X_{\epsilon,\sigma}^1}^2:u_N = v + w\Bigr\}\Bigr].
\end{align*}
We note that this construction ensures that
\begin{align}
\|u\|_{\X_{\epsilon,\sigma}^{s,b}}^2\simeq \sum_{N\geq 1}N^{2s}\|u_N\|_{\X_{\epsilon,\sigma}^{0,b}}^2\qtq{and}\|u\|_{\Y_{\epsilon,\sigma}^{s,b}}^2\simeq \sum_{N\geq 1}N^{-2s}\|u_N\|_{\Y_{\epsilon,\sigma}^{0,b}}^2.\label{XY sq sum}
\end{align}

We start the proof of Theorem~\ref{t:LWP-details} with a trio of linear estimates:
\begin{lem}
For \(s\in \R\), \(b>\frac12\), and \(\sigma\in \{\pm\}\) we have the estimates
\begin{align}
\|\chi(t) S^{\epsilon,\sigma}(t) f\|_{\X_{\epsilon,\sigma}^{s,b}}&\lesssim_{s,b} \|f\|_{H^s},\label{NRG1}\\
\biggl\|\int_\R \chi(t)\chi\bigl(\tfrac{t'}4\bigr)\sgn(t')&S^{\epsilon,\sigma}(t-t')F(t')\,dt'\biggr\|_{\X_{\epsilon,\sigma}^{s,b}}\label{NRG2a}\\
&\lesssim_{s,b} \sup\left\{\left|\iint_{\R^2} F(t) w(t)\,dx\,dt\right|:\|w\|_{\Y_{\epsilon,\sigma}^{s,b}}\leq 1\right\},\notag\\
\biggl\|\int_\R \chi\bigl(\tfrac{t-t'}4\bigr)\sgn(t-t')&S^{\epsilon,\sigma}(t-t')F(t')\,dt'\biggr\|_{\X_{\epsilon,\sigma}^{s,b}}\label{NRG2b}\\
&\lesssim_{s,b} \sup\left\{\left|\iint_{\R^2} F(t) w(t)\,dx\,dt\right|:\|w\|_{\Y_{\epsilon,\sigma}^{s,b}}\leq 1\right\},\notag
\end{align}
uniformly for all \(\epsilon>0\).
\end{lem}
\bpf
The proofs of \eqref{NRG1}, \eqref{NRG2a}, and \eqref{NRG2b} follow essentially standard arguments; see, e.g., \cite[Proposition~2.12]{MR2233925}. However, for completeness we spell out the details here.

To prove \eqref{NRG1}, we first compute that the spacetime Fourier transform
\[
\tfrac1{2\pi}\iint_{\R^2}\chi(t) S^{\epsilon,\sigma}(t) f_N e^{-i(t\tau + x\xi)}\,dx\,dt = \hat \chi\bigl(\tau + \sigma\phi_\epsilon(\xi)\bigr)\hat f_N(\xi).
\]
We then have
\begin{align*}
\|\chi(t) S^{\epsilon,\sigma}(t) f_N\|_{\X_{\epsilon,\sigma}^{0,b}}&\leq\|\chi(t) S^{\epsilon,\sigma}(t) f_N\|_{X_{\epsilon,\sigma}^b}\\
&\lesssim \|\<\tau\>^b\hat \chi(\tau)\|_{L_\tau^2}\|f_N\|_{L^2}\lesssim_b \|f_N\|_{L^2},
\end{align*}
and the estimate \eqref{NRG1} now follows from \eqref{XY sq sum}.

For \eqref{NRG2a}, we use that \(S^{\epsilon,\sigma}(t)\) is a group with \eqref{NRG1} and Plancherel's Theorem to bound
\begin{align*}
\left\|\int_\R \chi(t)\chi\bigl(\tfrac{t'}4\bigr)\sgn(t')S^{\epsilon,\sigma}(t-t')F_N(t')\,dt'\right\|_{\X_{\epsilon,\sigma}^{0,b}}\lesssim_b \left\|\int_\R\overline{h\bigl(\tau + \sigma\phi_\epsilon(\xi)\bigr)}\widetilde F_N(\tau)\,d\tau\right\|_{L_\xi^2},
\end{align*}
where
\[
h(\tau) := \int_\R\chi\bigl(\tfrac t4\bigr)\sgn(t)e^{-it\tau}\,dt.
\]
Integrating by parts, we compute that
\eq{h bound}{
\bigl|h(\tau)\bigr|\lesssim \<\tau\>^{-1}.
}
For any decomposition \(F_N = G + H\) we can then bound
\begin{align*}
\left\|\int_\R\overline{h\bigl(\tau + \sigma\phi_\epsilon(\xi)\bigr)}\widetilde G(\tau)\,d\tau\right\|_{L_\xi^2} &\lesssim_b \|\<\tau + \sigma\phi_\epsilon(\xi)\>^{b-1}\widetilde G\|_{L_{\tau,\xi}^2},\\
\left\|\int_\R\overline{h\bigl(\tau + \sigma\phi_\epsilon(\xi)\bigr)}\widetilde H(\tau)\,d\tau\right\|_{L_\xi^2} &\lesssim \|\<\tau + \sigma\phi_\epsilon(\xi)\>^{-1}\widetilde H\|_{L_\xi^2L_\tau^1},
\end{align*}
where the first of these estimates uses the Cauchy--Schwarz inequality and that \(b>\frac12\). This gives us the estimate
\begin{align*}
\left\|\int_\R \chi(t)\chi\bigl(\tfrac{t'}4\bigr)\sgn(t')S^{\epsilon,\sigma}(t-t')F_N(t')\,dt'\right\|_{\X_{\epsilon,\sigma}^{0,b}}\lesssim_b \bigl\|\bigl\<\p_t + \tfrac\sigma{\epsilon^2}\Lambda_\epsilon\p_x\bigr\>^{-1}F_N\bigr\|_{\X_{\epsilon,\sigma}^{0,b}},
\end{align*}
and \eqref{NRG2a} then follows from \eqref{XY sq sum} and duality.

Turning to \eqref{NRG2b}, for any decomposition \(F_N = G + H\) we use Plancherel's Theorem and \eqref{h bound} to obtain
\begin{align*}
\left\|\int_\R \chi\bigl(\tfrac{t-t'}4\bigr)\sgn(t-t')S^{\epsilon,\sigma}(t-t')G(t')\,dt'\right\|_{X_{\epsilon,\sigma}^b}&\lesssim \bigl\|\bigl\<\tau + \sigma\phi_\epsilon(\xi)\bigr\>^b h\bigl(\tau + \sigma\phi_\epsilon(\xi)\bigr)\widetilde G\bigr\|_{L_{\tau,\xi}^2}\\
&\lesssim \|\<\tau + \sigma\phi_\epsilon(\xi)\>^{b-1}\widetilde G\|_{L_{\tau,\xi}^2},\\
\left\|\int_\R \chi\bigl(\tfrac{t-t'}4\bigr)\sgn(t-t')S^{\epsilon,\sigma}(t-t')H(t')\,dt'\right\|_{Y_{\epsilon,\sigma}^0}&\lesssim \bigl\|h\bigl(\tau + \sigma\phi_\epsilon(\xi)\bigr)\widetilde H\bigr\|_{L_\xi^2L_\tau^1}\\
&\lesssim \|\<\tau + \sigma\phi_\epsilon(\xi)\>^{-1}\widetilde H\|_{L_\xi^2L_\tau^1}.
\end{align*}
Combining these gives us the estimate
\begin{align*}
\left\|\int_\R \chi\bigl(\tfrac{t-t'}4\bigr)\sgn(t-t')S^{\epsilon,\sigma}(t-t')F_N(t')\,dt'\right\|_{\X_{\epsilon,\sigma}^{0,b}}\lesssim \|\<\p_t + \tfrac\sigma{\epsilon^2}\Lambda_\epsilon\p_x\>^{-1}F_N\|_{\X_{\epsilon,\sigma}^{0,b}},
\end{align*}
from which we obtain \eqref{NRG2b} using \eqref{XY sq sum} and duality once more.
\epf

With suitable linear estimates in hand, let us now turn our attention to nonlinear estimates.  Our starting point is an estimate for the resonance function associated to \eqref{GB-system-2}:
\begin{lem}\label{l:modulation}
Let \(\xi_1,\xi_2,\xi_3\in \R\) satisfy \(\xi_1 + \xi_2 + \xi_3 = 0\). Then we have the estimates
\begin{align}
|\phi_\epsilon(\xi_1) + \phi_\epsilon(\xi_2) + \phi_\epsilon(\xi_3)| &\gtrsim \tfrac{\xi_{\min}\xi_{\max}^2}{1+\epsilon\xi_{\max}},\label{modulation aligned}\\
|\phi_\epsilon(\xi_1) - \phi_\epsilon(\xi_2) - \phi_\epsilon(\xi_3)| &\gtrsim \tfrac1{\epsilon^2}\xi_{\min}(1+\epsilon\xi_{\max}),\label{modulation misaligned}
\end{align}
uniformly for all \(\epsilon>0\), where \(\xi_{\max} = \max\{|\xi_1|,|\xi_2|,|\xi_3|\}\) and \(\xi_{\min} = \min\{|\xi_1|,|\xi_2|,|\xi_3|\}\).
\end{lem}
\bpf
It suffices to consider the case that \(\epsilon = 1\); the general case follows from rescaling according to \eqref{phi epsilon}.

Before proceeding, let us compute that
\begin{align}
\phi'(\xi) &= \tfrac1{\sqrt 3}\Bigl[\sqrt{4+\xi^2} + \tfrac{\xi^2}{\sqrt{4+\xi^2}}\Bigr],\label{phi d1}\\
\phi''(\xi) &= \tfrac{2\xi(6 + \xi^2)}{\sqrt 3(4 + \xi^2)^{3/2}},\label{phi d2}
\end{align}
and, for future use, that
\eq{phi d1 corrected}{
\phi'(\xi) - \tfrac 2{\sqrt 3} = \tfrac1{\sqrt 3}\Bigl[\tfrac{\xi^2}{2+\sqrt {4+\xi^2}} + \tfrac{\xi^2}{\sqrt {4+\xi^2}}\Bigr].
}

Turning to the proof of \eqref{modulation aligned}, as \(\phi\) is odd and \(\xi_1+\xi_2+\xi_3 = 0\), we may assume without loss of generality that \(0\leq \xi_1\leq \xi_2\) so that \(\xi_{\min} = \xi_1\) and \(\xi_{\max} = - \xi_3\geq \xi_2\geq -\tfrac12\xi_3\).  From \eqref{phi d2}, for all \(\xi\geq 0\) we have
\[
\phi''(\xi)\simeq \tfrac{\xi}{1+\xi}.
\]
We may then bound
\begin{align*}
\phi(\xi_1+\xi_2) - \phi(\xi_1) - \phi(\xi_2) &= \int_0^{\xi_1}\int_0^{\xi_2}\phi''(x+y)\,dy\,dx\\
&\gtrsim \int_0^{\xi_1}\int_{\frac12\xi_2}^{\xi_2}\tfrac{x+y}{1+x+y}\,dy\,dx\gtrsim \int_0^{\xi_1}\tfrac{\xi_2^2}{1+\xi_2}\,dx\simeq \tfrac{\xi_{\min}\xi_{\max}^2}{1+\xi_{\max}},
\end{align*}
which yields \eqref{modulation aligned}.

For \eqref{modulation misaligned}, we may assume without loss of generality that \(0\leq \xi_2\leq |\xi_3|\). If \(\xi_{\min} = \xi_2\) then we may use that \(\phi\) is odd and increasing to bound
\[
|\phi(\xi_1) - \phi(\xi_2) - \phi(\xi_3)|\geq \phi(\xi_{\max}) \gtrsim  \xi_{\min}(1+\xi_{\max}).
\]
If we instead have that \(\xi_{\min} = |\xi_1|\) then we must have \(0\leq \xi_1\leq \xi_2\) and \(\xi_{\max} = -\xi_3\geq \xi_2\geq -\frac12\xi_3\).  If \(\xi_1 = 0\), the estimate is vacuous.  Otherwise, we may find some \(\eta\in(\xi_2,\xi_1+\xi_2)\) so that
\[
|\phi(\xi_1) - \phi(\xi_2) - \phi(\xi_3)|\geq \phi(\xi_1+\xi_2) - \phi(\xi_2) = \phi'(\eta)\xi_1\gtrsim \xi_{\min}(1+\xi_{\max}),
\]
where we have used \eqref{phi d1} to bound
\[
\phi'(\eta)\gtrsim 1+\eta\gtrsim 1 + \xi_{\max}.
\]
This completes the proof of \eqref{modulation misaligned}.
\epf

Next, we prove several frequency-localized linear and bilinear estimates in the spaces \(\X_{\epsilon,\sigma}^{0,b}\) and \(\Y_{\epsilon,\sigma}^{0,b}\).  The construction of our spaces will yield different estimates at low and high modulations.  In light of this, for \(M>0\) we define the spacetime Fourier multiplier
\[
\widetilde{(Q_{\leq M}^{\epsilon,\sigma}u)}(\tau,\xi) = \chi_{\leq M}(\tau + \sigma\phi_\epsilon(\xi))\widetilde u(\tau,\xi),
\]
which projects to low modulations, and take the corresponding projection to high modulations to be \(Q_{>M}^{\epsilon,\sigma} = 1 - Q_{\leq M}^{\epsilon,\sigma}\). We then have the following:
\begin{lem}\label{l:the beast}
Let \(\frac12<b\leq 1\). Then, for dyadic \(N\geq 1\), \(\sigma\in \{\pm\}\), and \(2\leq p< \frac2{2b-1}\) we have the estimates
\begin{align}
\|u_N\|_{L_t^\infty L_x^2} &\lesssim_b \|u_N\|_{\X_{\epsilon,\sigma}^{0,b}},\label{X-L2}\\
\|u_N\|_{L_t^pL_x^2} & \lesssim_{b,p} \|u_N\|_{\Y_{\epsilon,\sigma}^{0,b}},\label{Y-L2}
\end{align}
and for all \(M>0\) we have the high modulation estimates
\begin{align}
\|Q_{>M}^{\epsilon,\sigma}u_N\|_{L_{t,x}^2}&\lesssim_b \bigl\<\min\{M,\phi_\epsilon(N)\}\bigr\>^{-b}\|u_N\|_{\X_{\epsilon,\sigma}^{0,b}},\label{X-hi}\\
\|Q_{>M}^{\epsilon,\sigma}u_N\|_{L_{t,x}^2}&\lesssim_b  \<M\>^{b-1}\|u_N\|_{\Y_{\epsilon,\sigma}^{0,b}}.\label{Y-hi}
\end{align}

For dyadic \(N>1\) and \(M>0\) we have the Strichartz estimates,
\begin{align}
\|u_N\|_{L_t^4L_x^\infty} & \lesssim_b \bigl(\tfrac{1+\epsilon N}N\bigr)^{\frac14}\|u_N\|_{\X_{\epsilon,\sigma}^{0,b}},\label{X-L4}\\
\|Q_{\leq M}^{\epsilon,\sigma}u_N\|_{L_t^4L_x^\infty} & \lesssim_b \<M\>^{b-\frac12}\bigl(\tfrac{1+\epsilon N}N\bigr)^{\frac14}\|u_N\|_{\Y_{\epsilon,\sigma}^{0,b}},\label{Y-L4}
\end{align}
and, for any dyadic \(N_2\geq 32\), \(N_1\leq \tfrac1{32}N_2\), and any \(M_1,M_2>0\), we have the bilinear Strichartz estimates
\begin{align}
&\|Q_{\leq M_1}^{\epsilon,\sigma_1} u_1\times Q_{\leq M_2}^{\epsilon,\sigma_2} v_{N_2}\|_{L_{t,x}^2}\label{XX-bil}\\
&\qquad\lesssim_b \Bigl[1 + \tfrac{\<M_1\>^b}{\<\phi_\epsilon(1)\>^b}\Bigr]\Bigl[1+\tfrac{\<M_2\>^b}{\<\phi_\epsilon(N_2)\>^b}\Bigr]\bigl(\tfrac{1+\epsilon N_2}{N_2^2}\bigr)^{\frac12}\|u_1\|_{\X_{\epsilon,\sigma_1}^{0,b}}\|v_{N_2}\|_{\X_{\epsilon,\sigma_2}^{0,b}},\notag\\
&\|Q_{\leq M_1}^{\epsilon,\sigma_1} u_{\leq N_1}\times Q_{\leq M_2}^{\epsilon,\sigma_2} v_{N_2}\|_{L_{t,x}^2}\label{YX-bil}\\
&\qquad\lesssim_b \<M_1\>^{b-\frac12}\Bigl[1+\tfrac{\<M_2\>^b}{\<\phi_\epsilon(N_2)\>^b}\Bigr]\bigl(\tfrac{1+\epsilon N_2}{ N_2^2}\bigr)^{\frac12}\|u_{\leq N_1}\|_{\Y_{\epsilon,\sigma_1}^{0,b}}\|v_{N_2}\|_{\X_{\epsilon,\sigma_2}^{0,b}},\notag
\end{align}
for any \(\sigma_1,\sigma_2\in\{\pm\}\).

Finally, we may find \(0<\nu\leq 1\) so that if \(M = \nu N_1N_3\frac{N_2}{1+\epsilon N_2}>0\) for any dyadic \(N_3\) and dyadic \(1<N_1\leq N_2\) then we have the non-resonance estimate
\begin{align}
&\left|\iint_{\R^2}\chi(t)\times Q_{\leq M}^{\epsilon,\sigma_1}u_{N_1}\times Q_{\leq M}^{\epsilon,\sigma_2}v_{N_2}\times Q_{\leq M}^{\epsilon,\sigma_3}P_{N_3}w\,dx\,dt\right|\label{nonresonance}\\
&\qquad \lesssim_{b,k}\<M\>^{-k}\min\{N_1^{\frac12},N_3^{\frac12}\}\|u_{N_1}\|_{\X_{\epsilon,\sigma_1}^{0,b}}\|v_{N_2}\|_{\X_{\epsilon,\sigma_2}^{0,b}}\|P_{N_3}w\|_{\Y_{\epsilon,\sigma_3}^{0,b}},\notag
\end{align}
for any \(k\geq 1\) and \(\sigma_1,\sigma_2,\sigma_3\in\{\pm\}\).

All of these estimates are uniform in \(\epsilon>0\).
\end{lem}
\bpf
We start by applying Plancherel's Theorem, Minkowski's integral inequality, and the Hausdorff--Young inequality to prove that for any \(2\leq p\leq \infty\) we have
\[
\|u\|_{L_t^pL_x^2} = \|\hat u\|_{L_t^pL_\xi^2}\leq \|\hat u\|_{L_\xi^2L_t^p}\leq \|\widetilde u\|_{L_\xi^2L_\tau^{p'}}.
\]
Taking \(p=\infty\) and using that \(b>\frac12\) gives us
\eq{embedding A}{
\|u\|_{L_t^\infty L_x^2} \lesssim\|u\|_{Y_{\epsilon,\sigma}^0}\lesssim_b\|u\|_{X_{\epsilon,\sigma}^b},
}
which yields \eqref{X-L2}. If we instead have that \(2\leq p< \frac2{2b-1}\) then 
\[
\|u\|_{L_t^p L_x^2}\leq \|\widetilde u\|_{L_\xi^2L_\tau^{p'}}\lesssim_{b,p} \|u\|_{X_{\epsilon,\sigma}^{1-b}},
\]
giving us \eqref{Y-L2}.

The estimates \eqref{X-hi} and \eqref{Y-hi} follow immediately from the fact that for all \(b'\geq 0\) we have
\eq{embedding B}{
\|Q_{>M}^{\epsilon,\sigma} u\|_{L_{t,x}^2} = \|\chi_{>M}(\tau + \sigma\phi_\epsilon(\xi))\widetilde u\|_{L_{\tau,\xi}^2}\lesssim_{b'} \<M\>^{-b'}\|u\|_{X_{\epsilon,\sigma}^{b'}}.
}

For the remaining estimates, it will be useful to first clarify the roles of the two components of the \(\X_{\epsilon,\sigma}^{s,b}\) norm. To this end, for \(M>0\) we apply Bernstein's inequality, respectively \eqref{embedding A} and \eqref{embedding B}, to bound
\begin{gather*}
\|Q_{\leq M}^{\epsilon,\sigma}u_N\|_{X_{\epsilon,\sigma}^b}^2 \lesssim_b \tfrac{\<M\>^{2b}}{\<\phi_\epsilon(N)\>^{2b}}\<\phi_\epsilon(N)\>^{2b}\|u_N\|_{X_{\epsilon,\sigma}^0}^2,\\
\|Q_{>M}^{\epsilon,\sigma}u_N\|_{Y_{\epsilon,\sigma}^0}^2 + \<\phi_\epsilon(N)\>^{2b}\|Q_{>M}^{\epsilon,\sigma}u_N\|_{X_{\epsilon,\sigma}^0}^2\lesssim_b\bigl[1 + \tfrac{\<\phi_\epsilon(N)\>^{2b}}{\<M\>^{2b}}\bigr]\|u_N\|_{X_{\epsilon,\sigma}^b}^2.
\end{gather*}
It follows that
\begin{gather}
\|Q_{\leq \phi_\epsilon(N)}^{\epsilon,\sigma}u_N\|_{X_{\epsilon,\sigma}^b}^2\lesssim_b \|u_N\|_{\X_{\epsilon,\sigma}^{0,b}}^2,\label{X-low-mod}\\
\|Q_{>\phi_\epsilon(N)}^{\epsilon,\sigma}u_N\|_{Y_{\epsilon,\sigma}^0}^2 + \<\phi_\epsilon(N)\>^{2b}\|Q_{>\phi_\epsilon(N)}^{\epsilon,\sigma}u_N\|_{X_{\epsilon,\sigma}^0}^2\lesssim_b \|u_N\|_{\X_{\epsilon,\sigma}^{0,b}}^2.\label{X-hi-mod}
\end{gather}

To prove the Strichartz estimates \eqref{X-L4} and \eqref{Y-L4}, for \(N>1\) we first apply \cite[Theorem 2.1]{MR1101221} with \eqref{phi d2} to obtain the Strichartz estimate
\eq{free Strichartz}{
\|S^{\epsilon,\sigma}(t) f_N\|_{L_t^4L_x^\infty} \lesssim \bigl(\tfrac{1+\epsilon N}N\bigr)^{\frac14}\|f_N\|_{L^2}.
}
We now employ the transference principle (see, e.g., \cite[Lemma 2.9]{MR2233925}) to obtain \eqref{Y-L4}. Precisely, taking
\[
F(\tau,x) = \tfrac1{\sqrt{2\pi}}\int_\R S^{\epsilon,\sigma}(-t)Q_{\leq M}^{\epsilon,\sigma}u_N(t,x)e^{-it\tau}\,dt,
\]
we may write
\[
Q_{\leq M}^{\epsilon,\sigma}u_N(t) = \tfrac1{\sqrt{2\pi}}\int_\R S^{\epsilon,\sigma}(t)F(\tau)e^{it\tau}\,d\tau.
\]
Applying \eqref{free Strichartz} and Plancherel's Theorem then yields the estimate
\begin{align*}
\|Q_{\leq M}^{\epsilon,\sigma}u_N\|_{L_t^4L_x^\infty} &\lesssim \bigl(\tfrac{1+\epsilon N}N\bigr)^{\frac14}\|F\|_{L_\tau^1L_x^2}\\
&\lesssim \bigl(\tfrac{1+\epsilon N}N\bigr)^{\frac14}\|\chi_{\leq M}(\tau)\widetilde u_N(\tau - \sigma\phi_\epsilon(\xi),\xi)\|_{L_\tau^1L_\xi^2}\\
&\lesssim_b \bigl(\tfrac{1+\epsilon N}N\bigr)^{\frac14}\<M\>^{b-\frac12}\|u_N\|_{X_{\epsilon,\sigma}^{1-b}},
\end{align*}
giving us \eqref{Y-L4}.

An identical argument proves that
\[
\|Q_{\leq M}^{\epsilon,\sigma}u_N\|_{L_t^4L_x^\infty}\lesssim_b \bigl(\tfrac{1+\epsilon N}N\bigr)^{\frac14}\|Q_{\leq M}^{\epsilon,\sigma}u_N\|_{X_{\epsilon,\sigma}^b},
\]
so by taking \(M = \phi_\epsilon(N)\) and applying \eqref{X-low-mod}, we obtain
\[
\|Q_{\leq \phi_\epsilon(N)}^{\epsilon,\sigma}u_N\|_{L_t^4L_x^\infty}\lesssim_b \bigl(\tfrac{1+\epsilon N}N\bigr)^{\frac14}\|u_N\|_{\X_{\epsilon,\sigma}^{0,b}}.
\]
For the high modulation contribution to the \(\X_{\epsilon,\sigma}^{0,b}\)-norm, we instead apply Bernstein's inequality, H\"older's inequality, \eqref{embedding A}, and \eqref{X-hi-mod} to bound
\begin{align*}
\|Q_{>\phi_\epsilon(N)}^{\epsilon,\sigma}u_N\|_{L_t^4L_x^\infty}&\lesssim N^{\frac12}\|Q_{>\phi_\epsilon(N)}^{\epsilon,\sigma}u_N\|_{L_t^4L_x^2}\\
&\lesssim  N^{\frac12}\|Q_{>\phi_\epsilon(N)}^{\epsilon,\sigma}u_N\|_{L_t^\infty L_x^2}^{\frac12}\|Q_{>\phi_\epsilon(N)}^{\epsilon,\sigma}u_N\|_{L_{t,x}^2}^{\frac12}\\
&\lesssim N^{\frac12}\<\phi_\epsilon(N)\>^{-\frac b2}\|Q_{>\phi_\epsilon(N)}^{\epsilon,\sigma}u_N\|_{Y_{\epsilon,\sigma}^0}^{\frac12}\Bigl[\<\phi_\epsilon(N)\>^b\|Q_{>\phi_\epsilon(N)}^{\epsilon,\sigma}u_N\|_{X_{\epsilon,\sigma}^0}\Bigr]^{\frac12}\\
&\lesssim_b \bigl(\tfrac{1+\epsilon N}N\bigr)^{\frac14}\|u_N\|_{\X_{\epsilon,\sigma}^{0,b}},
\end{align*}
which completes the proof of \eqref{X-L4}.

Next we consider the bilinear Strichartz estimates \eqref{XX-bil} and \eqref{YX-bil}.  Employing the transference principle with \eqref{X-low-mod} and \eqref{X-hi-mod} as in our proof of the Strichartz estimates \eqref{X-L4} and \eqref{Y-L4}, it suffices to prove a corresponding bilinear Strichartz estimate for free waves,
\begin{align}\label{free bilinear Strichartz}
\|S^{\epsilon,\sigma_1}(t) f_{\leq N_1}\times S^{\epsilon,\sigma_2}(t) g_{N_2}\|_{L_{t,x}^2}\lesssim \bigl(\tfrac{1+\epsilon N_2}{N_2^2}\bigr)^{\frac12}\|f_{\leq N_1}\|_{L^2}\|g_{N_2}\|_{L^2},
\end{align}
where \(\sigma_j\in \{\pm\}\).

To prove \eqref{free bilinear Strichartz}, we take a test function \(h\in L_{t,x}^2\) and consider the operator
\[
K := \int_\R S^{\epsilon,\sigma_1}(-t)P_{\leq 2 N_1}h(t)P_{\simeq N_2}S^{\epsilon,\sigma_2}(t)\,dt,
\]
on \(L^2\), where \(P_{\simeq N} = P_{\leq 2N}P_{> \frac14 N}\).  We then compute that
\begin{align*}
\left|\iint_{\R^2} S^{\epsilon,\sigma_1}(t)f_{\leq N_1}\times S^{\epsilon,\sigma_2} (t)g_{N_2}\times h\,dx\,dt\right|&= \left|\int_\R f_{\leq N_1}\times Kg_{N_2}\,dx\right|\\
&\leq \|f_{\leq N_1}\|_{L^2}\|g_{N_2}\|_{L^2}\|K\|_{\I_2},
\end{align*}
where the Hilbert--Schmidt norm
\begin{align*}
\|K\|_{\I_2}^2 &= \iint_{\R^2} |\chi_{\leq 2N_1}(\xi)|^2|\chi_{\simeq N_2}(\eta)|^2 \bigl|\widetilde h\bigl(\sigma_2 \phi_\epsilon(\eta) - \sigma_1\phi_\epsilon(\xi),\xi-\eta\bigr)\bigr|^2\,d\xi\,d\eta,
\end{align*}
with \(\chi_{\simeq N}(\xi) = \chi_{\leq 2N}(\xi)\chi_{>\frac14N}(\xi)\).  Setting \((\tau,\zeta) = (\sigma_2\phi_\epsilon(\eta) - \sigma_1\phi_\epsilon(\xi),\xi-\eta)\), we compute that the Jacobian determinant
\[
\bigl|\tfrac{\p(\tau,\zeta)}{\p(\xi,\eta)}\bigr| = |\sigma_1\phi_\epsilon'(\xi) - \sigma_2\phi_\epsilon'(\eta)|.
\]

For \(\xi\in \supp\chi_{\leq 2N_1}\) and \(\eta\in \supp \chi_{\simeq N_2}\) with \(N_1\leq \frac1{32}N_2\) and \(N_2>1\), we have \(2|\xi|\leq |\eta|\simeq N_2\).   From \eqref{phi d1 corrected}, we compute that for \(\zeta\geq 0\) we have
\begin{align*}
0\leq \phi_\epsilon'(\zeta) - \tfrac 2{\sqrt 3 \epsilon^2} = \tfrac1{2\sqrt 3}\Bigl[\tfrac{(2\zeta)^2}{4+\sqrt{16 + (2\epsilon\zeta)^2}} + \tfrac{(2\zeta)^2}{\sqrt{16 + (2\epsilon\zeta)^2}}\Bigr]\leq \tfrac 12\bigl(\phi_\epsilon'(2\zeta) - \tfrac 2{\sqrt 3 \epsilon^2}\bigr).
\end{align*}
By also using \eqref{phi d2}, we see that \( \phi_\epsilon'(\zeta) - \frac 2{\sqrt 3 \epsilon^2}\) is non-negative, even, and non-decreasing when \(\zeta\geq 0\). It follows that when \(\sigma_1 = \sigma_2\) we have
\[
\bigl|\tfrac{\p(\tau,\zeta)}{\p(\xi,\eta)}\bigr| \simeq \phi_\epsilon'(\eta) - \tfrac 2{\sqrt 3 \epsilon^2} \simeq \tfrac{N_2^2}{1+\epsilon N_2},
\]
where we have once again used \eqref{phi d1 corrected}. When \(\sigma_1 = -\sigma_2\), we instead have
\[
\bigl|\tfrac{\p(\tau,\zeta)}{\p(\xi,\eta)}\bigr|\simeq \bigl(\phi_\epsilon'(\eta) - \tfrac 2{\sqrt 3 \epsilon^2}\bigr) + \tfrac1{\epsilon^2} \gtrsim\tfrac{N_2^2}{1+\epsilon N_2}.
\]
In either case we may bound
\[
\|K\|_{\I_2}^2\lesssim \tfrac{1+\epsilon N_2}{N_2^2}\|h\|_{L_{t,x}^2}^2,
\]
which yields \eqref{free bilinear Strichartz}.

Finally, we turn to the non-resonance estimate \eqref{nonresonance}.  We set \(\vec\xi = (\xi_1,\xi_2,\xi_3)\), \(\vec\sigma = (\sigma_1,\sigma_2,\sigma_3)\), and
\[
\Xi_\epsilon(\vec\xi,\vec\sigma) = \sigma_1\phi_\epsilon(\xi_1) + \sigma_2\phi_\epsilon(\xi_2) + \sigma_3\phi_\epsilon(\xi_3).
\]
We then compute that
\begin{align*}
\LHS{nonresonance}&\simeq \Biggl|\iint\limits_{\substack{\tau_0 + \tau_1+\tau_2+\tau_3 = 0\\\xi_1+\xi_2+\xi_3 = 0}}\hat \chi\bigl(\tau_0 + \Xi_\epsilon(\vec\xi,\vec\sigma)\bigr) \chi_{\leq M}(\tau_1)\chi_{\leq M}(\tau_2)\chi_{\leq M}(\tau_3)\\
&\qquad\qquad\qquad\qquad\times\widetilde u_{N_1}\bigl(\tau_1 - \sigma_1\phi_\epsilon(\xi_1),\xi_1\bigr) \widetilde v_{N_2}\bigl(\tau_2 - \sigma_2\phi_\epsilon(\xi_2),\xi_2\bigr)\\
&\qquad\qquad\qquad\qquad\times\widetilde{P_{N_3}w}\bigl(\tau_3 - \sigma_3\phi_\epsilon(\xi_3),\xi_3\bigr)\,d\vec\xi\,d\vec\tau\Biggr|,
\end{align*}
where \(\vec\tau = (\tau_0,\tau_1,\tau_2,\tau_3)\). For \(\vec\tau,\vec\xi\) in the support of the integrand we may apply Lemma~\ref{l:modulation} to bound
\[
\bigl|\Xi_\epsilon(\vec\xi,\vec\sigma)\bigr|\gtrsim N_1N_3\tfrac{N_2}{1+\epsilon N_2}.
\]
Consequently, provided \(0<\nu\ll1\) is sufficiently small we have
\[
\bigl|\hat \chi\bigl(\tau_0 + \Xi_\epsilon(\vec\xi,\vec\sigma)\bigr) \bigr|\lesssim_k \<M\>^{-k}
\]
for any \(k\geq 1\). Applying H\"older's inequality, if \(N_1\leq N_3\) we have
\begin{align*}
\LHS{nonresonance} &\lesssim_k\<M\>^{-k}\|\widetilde u_{N_1}\|_{L_\xi^1L_\tau^1}\|\widetilde v_{N_2}\|_{L_\xi^2L_\tau^1}\bigl\|\chi_{\leq M}\bigl(\tau + \sigma_3\phi_\epsilon(\xi)\bigr)\widetilde{P_{N_3}w}\bigr\|_{L_\xi^2L_\tau^1}\\
&\lesssim_{b,k}\<M\>^{b-\frac12-k}N_1^{\frac12}\|u_{N_1}\|_{\X_{\epsilon,\sigma_1}^{0,b}}\|v_{N_2}\|_{\X_{\epsilon,\sigma_2}^{0,b}}\|P_{N_3}w\|_{\Y_{\epsilon,\sigma_3}^{0,b}},
\end{align*}
whereas if \(N_3\leq N_1\) we bound
\begin{align*}
\LHS{nonresonance} &\lesssim_k\<M\>^{-k}\|\widetilde u_{N_1}\|_{L_\xi^2L_\tau^1}\|\widetilde v_{N_2}\|_{L_\xi^2L_\tau^1}\bigl\|\chi_{\leq M}\bigl(\tau + \sigma_3\phi_\epsilon(\xi)\bigr)\widetilde{P_{N_3}w}\bigr\|_{L_\xi^1L_\tau^1}\\
&\lesssim_{b,k}\<M\>^{b-\frac12-k}N_3^{\frac12}\|u_{N_1}\|_{\X_{\epsilon,\sigma_1}^{0,b}}\|v_{N_2}\|_{\X_{\epsilon,\sigma_2}^{0,b}}\|P_{N_3}w\|_{\Y_{\epsilon,\sigma_3}^{0,b}}.
\end{align*}
Combining these yields \eqref{nonresonance}.
\epf

We now combine Lemmas~\ref{l:modulation} and~\ref{l:the beast} to prove our main nonlinear estimate:
\begin{prop}\label{p:bilinear estimates}
Let
\[
I[u,v,w] := \left|\iint_{\R^2}\chi(t)\,\Lambda_\epsilon^{-1}\p_x(uv)\,w\,dx\,dt\right|.
\]
Then, for \(s\geq -\frac12\) and \(\frac12<b<\frac 7{12}\) we have the estimate
\begin{align}\label{bilinear estimates}
I[u,v,w]&\lesssim_{s,b} \<\epsilon\>^{\frac12}\Bigl[\|u\|_{\X_{\epsilon,\sigma_1}^{-\frac12,b}}\|v\|_{\X_{\epsilon,\sigma_2}^{s,b}} + \|u\|_{\X_{\epsilon,\sigma_1}^{s,b}}\|v\|_{\X_{\epsilon,\sigma_2}^{-\frac12,b}}\Bigr]\|w\|_{\Y_{\epsilon,\sigma_3}^{s,b}},
\end{align}
uniformly for all \(\epsilon>0\) and \(\sigma_j\in \{\pm\}\).
\end{prop}
\bpf
We first decompose by frequency,
\begin{align*}
I[u,v,w]&\leq \sum_{N_1,N_2,N_3\geq 1}I[u_{N_1},v_{N_2},w_{N_3}].
\end{align*}
By symmetry, it suffices to restrict our attention to the case that \(N_1\leq N_2\). We also note that the highest two frequencies must be comparable, i.e., \(\frac14N_2\leq \max\{N_1,N_3\}\leq 4N_2\).

We now divide into several cases depending on the size of \(N_1,N_2,N_3\).  In several of the cases we will further decompose according to frequency and modulation.  In each case we will endeavor to bound the resulting integral in terms of a constant multiple of
\[
J(N_1,N_2,N_3) := \|u_{N_1}\|_{\X_{\epsilon,\sigma_1}^{-\frac12,b}}\|v_{N_2}\|_{\X_{\epsilon,\sigma_2}^{s,b}}\|w_{N_3}\|_{\Y_{\epsilon,\sigma_3}^{s,b}}.
\]
We will call such an estimate \emph{acceptable} if, after summation over admissible (according to the case) \(N_1,N_2,N_3\), we obtain a quantity that is bounded above by \RHS{bilinear estimates}. We suppress the dependence of implicit constants on \(s\) and \(b\) throughout.

\medskip

\noindent\emph{Case 1: Low frequencies, \(N_2 \leq 64\).} We start by addressing the case that all frequencies are low. Applying Bernstein's inequality followed by \eqref{X-L2} and \eqref{Y-L2}, we have
\begin{align*}
I[u_{N_1},v_{N_2},w_{N_3}] &\lesssim \|u_{N_1}\|_{L_{t,x}^\infty}\|v_{N_2}\|_{L_t^\infty L_x^2}\|w_{N_3}\|_{L_{t,x}^2}\\
&\lesssim \|u_{N_1}\|_{L_t^\infty L_x^2}\|v_{N_2}\|_{L_t^\infty L_x^2}\|w_{N_3}\|_{L_{t,x}^2}\\
&\lesssim J(N_1,N_2,N_3),
\end{align*}
where we have used that \(N_3\leq 256\).  This estimate is readily seen to be acceptable as \(N_1,N_2,N_3\lesssim1\).

\medskip

\noindent\emph{Case 2: Low $\times$ high $\to$ high interactions, \(N_2\geq 128\) and \(N_1\leq 16N_3\).} Next, we turn to the more straightforward low $\times$ high $\to$ high interactions. In this case we have \(\frac1{64} N_2\leq N_3\leq 4N_2\), so \(N_2\) and \(N_3\) are comparable. Here, we will subdivide into two subcases depending on whether \(N_1=1\) or \(N_1>1\):

\medskip

\noindent\emph{Case 2(a): \(N_2\geq 128\) and \(1 = N_1 \leq 16N_3\).} In light of the estimates \eqref{X-low-mod} and \eqref{X-hi-mod}, we start by decomposing the inputs \(u_{N_1}\) and \(v_{N_2}\) according to their modulation as
\begin{alignat*}{4}
u_{N_1} &= u_{N_1}^\ell + u_{N_1}^h,&&\qtq{where}u_{N_1}^\ell = Q_{\leq \phi_\epsilon(N_1)}^{\epsilon,\sigma_1}u_{N_1}&&\qtq{and}u_{N_1}^h = Q_{>\phi_\epsilon(N_1)}^{\epsilon,\sigma_1}u_{N_1},\\
v_{N_2} &= v_{N_2}^\ell + v_{N_2}^h,&&\qtq{where}v_{N_2}^\ell = Q_{\leq \phi_\epsilon(N_2)}^{\epsilon,\sigma_2}v_{N_2}&&\qtq{and}v_{N_2}^h = Q_{>\phi_\epsilon(N_2)}^{\epsilon,\sigma_2}v_{N_2}.
\end{alignat*}
We then use the triangle inequality to estimate
\begin{align*}
I[u_1,v_{N_2},w_{N_3}] &\leq I[u_1,v_{N_2}^h,w_{N_3}] + I[u_1^h,v_{N_2}^\ell,w_{N_3}] + I[u_1^\ell,v_{N_2}^\ell,w_{N_3}].
\end{align*}

For the first summand we apply Bernstein's inequality, \eqref{X-L2}, \eqref{X-hi}, and \eqref{Y-L2} to bound
\begin{align*}
I[u_1,v_{N_2}^h,w_{N_3}]&\lesssim \tfrac{N_3}{1+\epsilon N_3}\|u_1\|_{L_{t,x}^\infty}\|v_{N_2}^h\|_{L_{t,x}^2}\|w_{N_3}\|_{L_{t,x}^2}\\
&\lesssim \tfrac{N_3}{1+\epsilon N_3}\|u_1\|_{L_t^\infty L_x^2}\|v_{N_2}^h\|_{L_{t,x}^2}\|w_{N_3}\|_{L_{t,x}^2}\\
&\lesssim \tfrac{N_3}{1+\epsilon N_3} \<\phi_\epsilon(N_2)\>^{-b}J(1,N_2,N_3)\\
&\lesssim J(1,N_2,N_3),
\end{align*}
which is acceptable.

For the second summand, we again use Bernstein's inequality at low frequency, followed by \eqref{X-hi}, \eqref{X-L2}, and \eqref{Y-L2} to obtain
\begin{align*}
I[u_1^h,v_{N_2}^\ell,w_{N_3}]&\lesssim \tfrac{N_3}{1+\epsilon N_3}\|u_1^h\|_{L_t^2L_x^\infty}\|v_{N_2}^\ell\|_{L_t^\infty L_x^2}\|w_{N_3}\|_{L_{t,x}^2}\\
&\lesssim \tfrac{N_3}{1+\epsilon N_3}\|u_1^h\|_{L_{t,x}^2}\|v_{N_2}^\ell\|_{L_t^\infty L_x^2}\|w_{N_3}\|_{L_{t,x}^2}\\
&\lesssim \tfrac{N_3}{1+\epsilon N_3}\<\phi_\epsilon(1)\>^{-b}J(1,N_2,N_3)\\
&\lesssim J(1,N_2,N_3),
\end{align*}
which is again acceptable.

For the final summand, we use \eqref{XX-bil} and \eqref{Y-L2} to estimate
\begin{align*}
I[u_1^\ell,v_{N_2}^\ell,w_{N_3}] &\lesssim \tfrac{N_3}{1+\epsilon N_3}\|u_1^\ell v_{N_2}^\ell\|_{L_{t,x}^2}\|w_{N_3}\|_{L_{t,x}^2}\lesssim J(1,N_2,N_3),
\end{align*}
which is once again acceptable.

\medskip

\noindent\emph{Case 2(b): \(N_2\geq 128\) and \(1 < N_1 \leq 16N_3\).}  Noting that we must also have \(N_3>1\) in this case, we choose \(0<\nu\leq 1\) and \(M = \nu N_1N_3\frac{N_2}{1+\epsilon N_2}\) as in \eqref{nonresonance}. We then decompose \(u_{N_1}\), \(v_{N_2}\), and \(w_{N_3}\) by modulation as
\begin{alignat*}{4}
u_{N_1} &= u_{N_1}^L + u_{N_1}^H,&&\qtq{where}u_{N_1}^L = Q_{\leq M}^{\epsilon,\sigma_1}u_{N_1}&&\qtq{and}u_{N_1}^H = Q_{>M}^{\epsilon,\sigma_1}u_{N_1},\\
v_{N_2} &= v_{N_2}^L + v_{N_2}^H,&&\qtq{where}v_{N_2}^L = Q_{\leq M}^{\epsilon,\sigma_2}v_{N_2}&&\qtq{and}v_{N_2}^H = Q_{>M}^{\epsilon,\sigma_2}v_{N_2},\\
w_{N_3} &= w_{N_3}^L + w_{N_3}^H,&&\qtq{where}w_{N_3}^L = Q_{\leq M}^{\epsilon,\sigma_3}w_{N_3}&&\qtq{and}w_{N_3}^H = Q_{> M}^{\epsilon,\sigma_3}w_{N_3},
\end{alignat*}
and bound
\begin{align*}
I[u_{N_1},v_{N_2},w_{N_3}] &\leq I[u_{N_1},v_{N_2},w_{N_3}^H] + I[u_{N_1},v_{N_2}^H,w_{N_3}^L]\\
&\quad + I[u_{N_1}^H,v_{N_2}^L,w_{N_3}^L] + I[u_{N_1}^L,v_{N_2}^L,w_{N_3}^L].
\end{align*}

For the first summand, we use \eqref{X-L2}, \eqref{X-L4}, and \eqref{Y-hi} to estimate
\begin{align*}
I[u_{N_1},v_{N_2},w_{N_3}^H] &\lesssim \tfrac{N_3}{1+\epsilon N_3}\|u_{N_1}\|_{L_t^\infty L_x^2}\|v_{N_2}\|_{L_t^4L_x^\infty}\|w_{N_3}^H\|_{L_{t,x}^2}\\
&\lesssim N_1^{\frac12}\bigl(\tfrac{N_3}{1+\epsilon N_3}\bigr)^{\frac34}\<M\>^{b-1}J(N_1,N_2,N_3)\\
&\lesssim \bigl(\tfrac{N_1}{N_2}\bigr)^{\frac1{12}}J(N_1,N_2,N_3),
\end{align*}
where the final inequality uses the fact that \(b\leq \frac7{12}\).  This estimate is readily seen to be acceptable after summing over \(N_1\leq N_2\) and then \(N_2\simeq N_3\).

For the second summand, we observe that \(M\lesssim \phi_\epsilon(N_2)\) and hence we may apply \eqref{X-L2}, \eqref{X-hi}, and \eqref{Y-L4} to bound
\begin{align*}
I[u_{N_1},v_{N_2}^H,w_{N_3}^L] & \lesssim \tfrac{N_3}{1+\epsilon N_3}\|u_{N_1}\|_{L_t^\infty L_x^2}\|v_{N_2}^H\|_{L_{t,x}^2}\|w_{N_3}^L\|_{L_t^4L_x^\infty}\\
&\lesssim N_1^{\frac12}\bigl(\tfrac{N_3}{1+\epsilon N_3}\bigr)^{\frac34}\<M\>^{-\frac12}J(N_1,N_2,N_3)\\
&\lesssim \bigl(\tfrac{N_1}{N_2}\bigr)^{\frac1{12}}J(N_1,N_2,N_3),
\end{align*}
which is identical to the first summand.

For the third summand, we note that \(M\lesssim \phi_\epsilon(N_1)\) when \(N_2<\epsilon^{-1}\) and \(M\gtrsim \phi_\epsilon(N_1)\) when \(N_2\geq \epsilon^{-1}\).  Using this, we may apply \eqref{X-hi}, \eqref{X-L4}, and \eqref{Y-L2} to estimate
\begin{align*}
I[u_{N_1}^H,v_{N_2}^L,w_{N_3}^L]&\lesssim \tfrac{N_3}{1+\epsilon N_3}\|u_{N_1}^H\|_{L_{t,x}^2}\|v_{N_2}^L\|_{L_t^4 L_x^\infty}\|w_{N_3}^L\|_{L_t^4L_x^2}\\
&\lesssim \begin{cases}N_1^{\frac12}N_3^{\frac34}\<M\>^{-b}J(N_1,N_2,N_3)&\qt{if \(N_2< \epsilon^{-1}\)}\\
\epsilon^{-\frac34}N_1^{\frac12}\<\phi_\epsilon(N_1)\>^{-b}J(N_1,N_2,N_3)&\qt{if \(N_2\geq \epsilon^{-1}\)}
\end{cases}\\
&\lesssim \begin{cases}\bigl(\tfrac{N_1}{N_2}\bigr)^{\frac1{12}}J(N_1,N_2,N_3)&\qt{if \(N_2<\epsilon^{-1}\)},\\
\tfrac{(\epsilon N_1)^{\frac1{12}}}{(1+\epsilon N_1)^{\frac5{12}}}J(N_1,N_2,N_3)&\qt{if \(N_2\geq \epsilon^{-1}\)},
\end{cases}
\end{align*}
which is acceptable.

For the final summand we apply \eqref{nonresonance} with \(k=1\) to obtain
\begin{align*}
I[u_{N_1}^L,v_{N_2}^L,w_{N_3}^L]&\lesssim \tfrac1{N_2}J(N_1,N_2,N_3),
\end{align*}
which is evidently also acceptable.

\medskip

\noindent\emph{Case 3: High $\times$ high $\to$ low interactions, \(N_2\geq 128\) and \(N_3\leq \tfrac1{32}N_1\).} We now turn to the more challenging high $\times$ high $\to$ low interactions, where \(N_1\) and \(N_2\) are now comparable.  This time, we not only distinguish between the case \(N_3>1\) and \(N_3=1\), but also between the case \(N_2<\epsilon^{-1}\) and \(N_2 \geq \epsilon^{-1}\).

\medskip

\noindent\emph{Case 3(a): \(128\leq N_2<\epsilon^{-1}\) and \(1<N_3\leq \tfrac1{32}N_1\).} Adopting the modulation decomposition from Case 2(a), we start by decomposing
\begin{align*}
I[u_{N_1},v_{N_2},w_{N_3}] &\leq I[u_{N_1}^h,v_{N_2},w_{N_3}] + I[u_{N_1}^\ell,v_{N_2}^h,w_{N_3}] + I[u_{N_1}^\ell,v_{N_2}^\ell,w_{N_3}].
\end{align*}

For the first summand, we apply \eqref{X-hi}, \eqref{X-L4}, and \eqref{Y-L2} to bound
\begin{align*}
I[u_{N_1}^h,v_{N_2},w_{N_3}]& \lesssim N_3\|u_{N_1}^h\|_{L_{t,x}^2}\|v_{N_2}\|_{L_t^4L_x^\infty}\|w_{N_3}\|_{L_t^4L_x^2}\\
&\lesssim N_1^{\frac14-s}N_3^{1+s}\<\phi_\epsilon(N_1)\>^{-b}J(N_1,N_2,N_3)\\
&\lesssim \bigl(\tfrac{N_3}{N_1}\bigr)^{1+s}J(N_1,N_2,N_3),
\end{align*}
which is readily seen to be acceptable by first summing in \(N_3\lesssim N_1\) and then \(N_1\simeq N_2\).

The second summand is treated in a parallel manner, reversing the roles of \(u_{N_1}\) and \(v_{N_2}\) to obtain an identical contribution,
\begin{align*}
I[u_{N_1}^\ell,v_{N_2}^h,w_{N_3}]\lesssim N_3\|u_{N_1}^\ell\|_{L_t^4L_x^\infty}\|v_{N_2}^h\|_{L_{t,x}^2}\|w_{N_3}\|_{L_t^4L_x^2}\lesssim \bigl(\tfrac{N_3}{N_2}\bigr)^{1+s}J(N_1,N_2,N_3).
\end{align*}

For the third summand, we proceed as in Case 2(b) and further decompose in modulation,
\begin{align*}
I[u_{N_1}^\ell,v_{N_2}^\ell,w_{N_3}] &\leq I[u_{N_1}^\ell,v_{N_2}^\ell,w_{N_3}^H] + I[u_{N_1}^\ell,v_{N_2}^{\ell,H},w_{N_3}^L]\\
&\quad + I[u_{N_1}^{\ell,H},v_{N_2}^{\ell,L},w_{N_3}^L] + I[u_{N_1}^{\ell,L},v_{N_2}^{\ell,L},w_{N_3}^L],
\end{align*}
where we again take \(M = \nu N_1N_3\frac{N_2}{1+\epsilon N_2}\) as in \eqref{nonresonance} and observe that \(M\lesssim \phi_\epsilon(N_1)\simeq \phi_\epsilon(N_2)\).

For the first term, we use \eqref{X-L2}, \eqref{X-L4}, and \eqref{Y-hi} to bound
\begin{align*}
I[u_{N_1}^\ell,v_{N_2}^\ell,w_{N_3}^H] &\lesssim N_3\|u_{N_1}^\ell\|_{L_t^\infty L_x^2}\|v_{N_2}^\ell\|_{L_t^4L_x^\infty}\|w_{N_3}^H\|_{L_{t,x}^2}\\
&\lesssim N_3^{1+s}N_2^{\frac14-s}\<M\>^{b-1}J(N_1,N_2,N_3)\\
&\lesssim \bigl(\tfrac{N_3}{N_2}\bigr)^{\frac7{12}+s}J(N_1,N_2,N_3),
\end{align*}
which is acceptable. Similarly, the second term is readily estimated using \eqref{YX-bil} and \eqref{X-hi} to obtain
\begin{align*}
I[u_{N_1}^\ell,v_{N_2}^{\ell,H},w_{N_3}^L]&\lesssim \|u_{N_1}^\ell\times \Lambda_\epsilon^{-1}\p_x w_{N_3}^L\|_{L_{t,x}^2}\|v_{N_2}^{\ell,H}\|_{L_{t,x}^2}\\
&\lesssim N_1^{-\frac12-s}N_3^{1+s}\<M\>^{-\frac12}J(N_1,N_2,N_3)\\
&\lesssim \bigl(\tfrac{N_3}{N_1}\bigr)^{\frac56+s}J(N_1,N_2,N_3),
\end{align*}
and by reversing the roles of \(u_{N_1}\) and \(v_{N_2}\) we get an identical estimate for the third term,
\begin{align*}
I[u_{N_1}^{\ell,H},v_{N_2}^{\ell,L},w_{N_3}^L]\lesssim \|v_{N_2}^{\ell,L}\times \Lambda_\epsilon^{-1}\p_x w_{N_3}^L\|_{L_{t,x}^2}\|u_{N_1}^{\ell,H}\|_{L_{t,x}^2}\lesssim
\bigl(\tfrac{N_3}{N_2}\bigr)^{\frac56+s}J(N_1,N_2,N_3),
\end{align*}
both of which are acceptable. For the final term we invoke \eqref{nonresonance} with \(k=1\) to get
\begin{align*}
I[u_{N_1}^{\ell,L},v_{N_2}^{\ell,L},w_{N_3}^L]&\lesssim N_2^{-1}J(N_1,N_2,N_3),
\end{align*}
which is acceptable once again.

\medskip

\noindent\emph{Case 3(b): \(N_2\geq \max\{128,\epsilon^{-1}\}\) and \(1<N_3\leq \tfrac1{32}N_1\).}  We proceed as in Case 3(a) and decompose
\begin{align*}
I[u_{N_1},v_{N_2},w_{N_3}] &\leq  I[u_{N_1}^h,v_{N_2},w_{N_3}] + I[u_{N_1}^\ell,v_{N_2}^h,w_{N_3}] + I[u_{N_1}^\ell,v_{N_2}^\ell,w_{N_3}].
\end{align*}

For the first and second summands, we use the same estimates as in Case 3(a) to bound
\begin{align*}
&I[u_{N_1}^h,v_{N_2},w_{N_3}] + I[u_{N_1}^\ell,v_{N_2}^h,w_{N_3}]\\
&\qquad\lesssim
\begin{cases}\tfrac{(\epsilon N_3)^{1+s}}{1+\epsilon N_3}(\epsilon N_2)^{-\frac13-s}J(N_1,N_2,N_3)&\qt{if \(s\geq -\frac13\)},\\
\tfrac{(\epsilon N_3)^{1+s}}{1+\epsilon N_3}\epsilon^{-\frac{3s+1}2}J(N_1,N_2,N_3)&\qt{if \(s\in[-\frac12,-\frac13)\)}.
\end{cases}
\end{align*}
To verify that this is acceptable, it is best to treat the cases \(0<\epsilon\leq 1\) and \(\epsilon>1\) separately. In each of these two cases, we may then consider the sum in \(N_3\) followed by the sum in \(N_1\simeq N_2\) for each of the cases \(-\frac12\leq s<-\frac13\), \(-\frac13\leq s<0\), \(s=0\), and \(s>0\).

For the third summand, we again decompose
\begin{align*}
I[u_{N_1}^\ell,v_{N_2}^\ell,w_{N_3}] &\leq I[u_{N_1}^\ell,v_{N_2}^\ell,w_{N_3}^H] + I[u_{N_1}^\ell,v_{N_2}^{\ell,H},w_{N_3}^L]\\
&\quad + I[u_{N_1}^{\ell,H},v_{N_2}^{\ell,L},w_{N_3}^L] + I[u_{N_1}^{\ell,L},v_{N_2}^{\ell,L},w_{N_3}^L].
\end{align*}
The second, third, and fourth terms are estimated just as in Case 3(a), which yields the bounds
\begin{align*}
I[u_{N_1}^\ell,v_{N_2}^{\ell,H},w_{N_3}^L]& + I[u_{N_1}^{\ell,H},v_{N_2}^{\ell,L},w_{N_3}^L]\\
&\lesssim \begin{cases}
\tfrac{(\epsilon N_3)^{\frac56+s}}{1+\epsilon N_3}(\epsilon N_2)^{-\frac16-s}J(N_1,N_2,N_3)&\qt{if \(s\geq -\frac16\)},\\
\tfrac{(\epsilon N_3)^{1+2s}}{1+\epsilon N_3}\epsilon^{-\frac12-3s}J(N_1,N_2,N_3)&\qt{if \(s\in[-\frac12, -\frac16)\)},
\end{cases}\\
I[u_{N_1}^{\ell,L},v_{N_2}^{\ell,L},w_{N_3}^L]&\lesssim \tfrac{(\epsilon N_3)^{\frac12+s}}{1+\epsilon N_3} (\epsilon N_2)^{-\frac12-s} \epsilon J(N_1,N_2,N_3),
\end{align*}
where we have used that \(M\lesssim \phi_\epsilon(N_1)\simeq \phi_\epsilon(N_2)\) in the first of these two estimates. By once again considering several cases depending on the size of \(\epsilon\) and \(s\), one may verify that both of these bounds are acceptable.

The remaining term requires the full power of the \(\Y_{\epsilon,\sigma_3}^{s,b}\)-norm. We start by applying Young's convolution inequality and H\"older's inequality on the Fourier side, followed by the definitions of our norms and \eqref{embedding A} to obtain
\begin{align*}
&I[u_{N_1}^\ell,v_{N_2}^\ell,w_{N_3}^H]\\
&\lesssim \tfrac{N_3}{1+\epsilon N_3}\|\hat \chi\|_{L_\tau^1}\|\widetilde u_{N_1}^\ell\|_{L_\xi^2L_\tau^1}\|\widetilde v_{N_2}^\ell\|_{L_\xi^2L_\tau^1}\|\widetilde w_{N_3}^H\|_{L_\xi^1L_\tau^\infty}\\
&\lesssim \tfrac{N_3^{\frac32}}{1+\epsilon N_3}\<M\>^{-1}\|u_{N_1}^\ell\|_{Y_{\epsilon,\sigma_1}^0} \|v_{N_2}^\ell\|_{Y_{\epsilon,\sigma_2}^0} \|w_{N_3}^H\|_{Z_{\epsilon,\sigma_3}^1}\\
&\lesssim \begin{cases}\tfrac{(\epsilon N_3)^{\frac56+s}}{(1+\epsilon N_3)}(\epsilon N_2)^{-\frac16-s}\|u_{N_1}\|_{\X_{\epsilon,\sigma_1}^{-\frac12,b}}\|v_{N_2}\|_{\X_{\epsilon,\sigma_2}^{s,b}}\Bigl[N_3^{-s}\|w_{N_3}\|_{Z_{\epsilon,\sigma_3}^1}\Bigr]&\qt{if \(s\geq -\frac16\)}\\
\tfrac{(\epsilon N_3)^{1+2s}}{1+\epsilon N_3}\epsilon^{-\frac12-3s}\|u_{N_1}\|_{\X_{\epsilon,\sigma_1}^{-\frac12,b}}\|v_{N_2}\|_{\X_{\epsilon,\sigma_2}^{s,b}}\Bigl[N_3^{-s}\|w_{N_3}\|_{Z_{\epsilon,\sigma_3}^1}\Bigr]&\qt{if \(s\in[-\frac12,-\frac16)\)},\end{cases}
\end{align*}
where we have used that \(M\gtrsim \phi_\epsilon(N_3)\).

For the remaining part of the norm, if \(\phi_\epsilon(N_3)\leq 1\) then we apply \eqref{X-L2} and Bernstein's inequality to estimate
\begin{align*}
&I[u_{N_1}^\ell,v_{N_2}^\ell,w_{N_3}^H]\\
&\lesssim \tfrac{N_3}{1+\epsilon N_3}\|u_{N_1}^\ell\|_{L_t^\infty L_x^2} \|v_{N_2}^\ell\|_{L_t^\infty L_x^2} \|w_{N_3}^H\|_{L_t^2L_x^\infty}\\
&\lesssim \tfrac{(\epsilon N_3)^{\frac32+s}}{(1+\epsilon N_3)}(\epsilon N_2)^{\frac12-s}\epsilon^{-2}\<M\>^{-1}\|u_{N_1}\|_{\X_{\epsilon,\sigma_1}^{-\frac12,b}}\|v_{N_2}\|_{\X_{\epsilon,\sigma_2}^{s,b}}\Bigl[\tfrac{N_3^{-s}}{\<\phi_\epsilon(N_3)\>^b}\|w_{N_3}\|_{X_{\epsilon,\sigma_3}^1}\Bigr]\\
&\lesssim \begin{cases}\tfrac{(\epsilon N_3)^{\frac56+s}}{1+\epsilon N_3}(\epsilon N_2)^{-\frac16-s}\|u_{N_1}\|_{\X_{\epsilon,\sigma_1}^{-\frac12,b}}\|v_{N_2}\|_{\X_{\epsilon,\sigma_2}^{s,b}}\Bigl[\tfrac{N_3^{-s}}{\<\phi_\epsilon(N_3)\>^b}\|w_{N_3}\|_{X_{\epsilon,\sigma_3}^1}\Bigr]&\text{if \(s\geq -\tfrac16\)}\\
\tfrac{(\epsilon N_3)^{1+2s}}{1+\epsilon N_3}\epsilon^{-\frac12-3s}\|u_{N_1}\|_{\X_{\epsilon,\sigma_1}^{-\frac12,b}}\|v_{N_2}\|_{\X_{\epsilon,\sigma_2}^{s,b}}\Bigl[\tfrac{N_3^{-s}}{\<\phi_\epsilon(N_3)\>^b}\|w_{N_3}\|_{X_{\epsilon,\sigma_3}^1}\Bigr]&\text{if \(s\in[-\tfrac12,-\tfrac16)\)},\end{cases}
\end{align*}
whereas if \(\phi_\epsilon(N_3)>1\) we apply \eqref{X-L2} and \eqref{X-L4} to bound
\begin{align*}
&I[u_{N_1}^\ell,v_{N_2}^\ell,w_{N_3}^H]\\
&\lesssim \tfrac{N_3}{1+\epsilon N_3}\|u_{N_1}^\ell\|_{L_t^\infty L_x^2} \|v_{N_2}^\ell\|_{L_t^4 L_x^\infty} \|w_{N_3}^H\|_{L_{t,x}^2}\\
&\lesssim \tfrac{(\epsilon N_3)^{1+s}(\epsilon N_2)^{\frac12-s}}{(1+\epsilon N_3)\epsilon^{\frac54}}\<M\>^{-1}\<\phi_\epsilon(N_3)\>^b\|u_{N_1}\|_{\X_{\epsilon,\sigma_1}^{-\frac12,b}}\|v_{N_2}\|_{\X_{\epsilon,\sigma_2}^{s,b}}\Bigl[\tfrac{N_3^{-s}}{\<\phi_\epsilon(N_3)\>^b}\|w_{N_3}\|_{X_{\epsilon,\sigma_3}^1}\Bigr]\\
&\lesssim\!\begin{cases}\tfrac{(\epsilon N_3)^{\frac7{12}+s}(\epsilon N_2)^{\frac1{12}-s-b}}{(1+\epsilon N_3)^{1-b}}\|u_{N_1}\|_{\X_{\epsilon,\sigma_1}^{-\frac12,b}}\|v_{N_2}\|_{\X_{\epsilon,\sigma_2}^{s,b}}\Bigl[\tfrac{N_3^{-s}}{\<\phi_\epsilon(N_3)\>^b}\|w_{N_3}\|_{X_{\epsilon,\sigma_3}^1}\Bigr]&\!\!\!\!\text{if \(s\geq \tfrac1{12}-b\)},\\
\tfrac{(\epsilon N_3)^{\frac12+2s+b}\epsilon^{\frac14-3s-3b}}{(1+\epsilon N_3)^{1-b}}\|u_{N_1}\|_{\X_{\epsilon,\sigma_1}^{-\frac12,b}}\|v_{N_2}\|_{\X_{\epsilon,\sigma_2}^{s,b}}\Bigl[\tfrac{N_3^{-s}}{\<\phi_\epsilon(N_3)\>^b}\|w_{N_3}\|_{X_{\epsilon,\sigma_3}^1}\Bigr]&\!\!\!\!\text{if \(s\in[-\tfrac12,\tfrac1{12}\!-\!b)\)}.\end{cases}
\end{align*}

Decomposing into several cases once again shows that by combining these estimates we obtain an acceptable bound.

\medskip

\noindent\emph{Case 3(c):  \(128\leq N_2<\epsilon^{-1}\) and \(1=N_3\leq \tfrac1{32}N_1\).} In this case we further decompose
\[
I[u_{N_1},v_{N_2},w_1] \leq \sum_{N_3\leq 2}I[u_{N_1},v_{N_2},w_{N_3}],
\]
where the sum is over \emph{all} dyadic \(N_3\leq 2\) and, by an abuse of notation, we take \(w_{N_3} = P_{N_3}w_1\) for the remainder of this case (and for Case 3(d)). We similarly take
\[
J(N_1,N_2,N_3) = \|u_{N_1}\|_{\X_{\epsilon,\sigma_1}^{-\frac12,b}}\|v_{N_2}\|_{\X_{\epsilon,\sigma_2}^{s,b}}\|P_{N_3}w_1\|_{\Y_{\epsilon,\sigma_3}^{s,b}},
\]
so that
\eq{summation helper for lf case}{
\sum_{N_1\simeq N_2\gtrsim 1}\sup_{N_3\leq 2}J(N_1,N_2,N_3)\lesssim \|u\|_{\X_{\epsilon,\sigma_1}^{-\frac12,b}}\|v\|_{\X_{\epsilon,\sigma_2}^{s,b}}\|w\|_{\Y_{\epsilon,\sigma_3}^{s,b}}.
}

Proceeding precisely as in Case 3(a), and retaining the choice of \(M = \nu N_1N_3\frac{N_2}{1+\epsilon N_2}\) but now with \(N_3\leq 2\), we obtain the estimates
\begin{align*}
I[u_{N_1}^h,v_{N_2},w_{N_3}] + I[u_{N_1}^\ell,v_{N_2}^h,w_{N_3}]&\lesssim \tfrac{N_3}{N_2^{1+s}}J(N_1,N_2,N_3),
\\
I[u_{N_1}^\ell,v_{N_2}^\ell,w_{N_3}^H] &\lesssim \tfrac{N_3^{\frac7{12}}}{N_2^{\frac7{12}+s}}J(N_1,N_2,N_3),
\\
I[u_{N_1}^\ell,v_{N_2}^{\ell,H},w_{N_3}^L]+
I[u_{N_1}^{\ell,H},v_{N_2}^{\ell,L},w_{N_3}^L]&\lesssim \tfrac{N_3^{\frac56}}{N_2^{\frac56+s}}J(N_1,N_2,N_3),\\
I[u_{N_1}^{\ell,L},v_{N_2}^{\ell,L},w_{N_3}^L]&\lesssim \tfrac{N_3^{\frac12}}{N_2^{\frac12+s}}J(N_1,N_2,N_3),
\end{align*}
all of which are seen to be acceptable using \eqref{summation helper for lf case}.

\medskip

\noindent\emph{Case 3(d): \(N_2\geq \max\{128,\epsilon^{-1}\}\) and \(1=N_3\leq \tfrac1{32}N_1\).} Decomposing \(w_1\) as in Case 3(c) (with an identical notational convention) and then estimating each term as in Case 3(b), we similarly obtain the estimates
\begin{align*}
I[u_{N_1}^h,v_{N_2},w_{N_3}] &+ I[u_{N_1}^\ell,v_{N_2}^h,w_{N_3}] \\&\lesssim  \begin{cases}\tfrac{\epsilon N_3}{1+\epsilon N_3}(\epsilon N_2)^{-\frac13-s} \epsilon^s J(N_1,N_2,N_3)&\qt{if \(s\geq -\frac13\)},\\
\tfrac{\epsilon N_3}{1+\epsilon N_3}\epsilon^{-\frac{s+1}2}J(N_1,N_2,N_3)&\qt{if \(s\in[-\frac12,-\frac13)\)},
\end{cases}\\
I[u_{N_1}^\ell,v_{N_2}^{\ell,H},w_{N_3}^L] &+ I[u_{N_1}^{\ell,H},v_{N_2}^{\ell,L},w_{N_3}^L]\\&\lesssim \begin{cases}
\tfrac{(\epsilon N_3)^{\frac56}}{1+\epsilon N_3}(\epsilon N_2)^{-\frac16-s}\epsilon^sJ(N_1,N_2,N_3)&\qt{if \(s\geq -\frac16\)},\\
\tfrac{(\epsilon N_3)^{1+s}}{1+\epsilon N_3}\epsilon^{-\frac12-2s}J(N_1,N_2,N_3)&\qt{if \(s\in[-\frac12, -\frac16)\)},
\end{cases}\\
I[u_{N_1}^{\ell,L},v_{N_2}^{\ell,L},w_{N_3}^L]&\lesssim \tfrac{(\epsilon N_3)^{\frac12}}{1+\epsilon N_3} (\epsilon N_2)^{-\frac12-s} \epsilon^{1+s} J(N_1,N_2,N_3),
\end{align*}
all of which are acceptable. We remark that the last two of these estimates are the first place that the factor of \(\<\epsilon\>^{\frac12}\) on \RHS{bilinear estimates} shows up.

For the remaining term, we again argue as in Case 3(b) to bound
\begin{align*}
&I[u_{N_1}^\ell,v_{N_2}^\ell,w_{N_3}^H]\lesssim\begin{cases}\tfrac{(\epsilon N_3)^{\frac56}\epsilon^s}{(1+\epsilon N_3)(\epsilon N_2)^{\frac16+s}}\|u_{N_1}\|_{\X_{\epsilon,\sigma_1}^{-\frac12,b}}\|v_{N_2}\|_{\X_{\epsilon,\sigma_2}^{s,b}}\|w_{N_3}\|_{Z_{\epsilon,\sigma_3}^1}&\text{if \(s\geq -\frac16\)},\\
\tfrac{(\epsilon N_3)^{1+s}\epsilon^{-\frac12-2s}}{1+\epsilon N_3}\|u_{N_1}\|_{\X_{\epsilon,\sigma_1}^{-\frac12,b}}\|v_{N_2}\|_{\X_{\epsilon,\sigma_2}^{s,b}}\|w_{N_3}\|_{Z_{\epsilon,\sigma_3}^1}&\text{if \(s\in[-\frac12,-\frac16)\)},\end{cases}
\end{align*}
which is the second place the factor of \(\<\epsilon\>^{\frac12}\) on \RHS{bilinear estimates} appears. For the remaining component, we proceed as in Case 3(b) once more and when \(\phi_\epsilon(1)\leq 1\) we estimate
\begin{align*}
&I[u_{N_1}^\ell,v_{N_2}^\ell,w_{N_3}^H]\\&\lesssim \begin{cases}\tfrac{(\epsilon N_3)^{\frac56}}{1+\epsilon N_3}(\epsilon N_2)^{-\frac16-s}\epsilon^s\|u_{N_1}\|_{\X_{\epsilon,\sigma_1}^{-\frac12,b}}\|v_{N_2}\|_{\X_{\epsilon,\sigma_2}^{s,b}}\Bigl[\tfrac1{\<\phi_\epsilon(1)\>^b}\|w_{N_3}\|_{X_{\epsilon,\sigma_3}^1}\Bigr]&\qt{if \(s\geq -\tfrac16\)},\\
\tfrac{(\epsilon N_3)^{1+s}}{1+\epsilon N_3}\epsilon^{-\frac12-2s}\|u_{N_1}\|_{\X_{\epsilon,\sigma_1}^{-\frac12,b}}\|v_{N_2}\|_{\X_{\epsilon,\sigma_2}^{s,b}}\Bigl[\tfrac1{\<\phi_\epsilon(1)\>^b}\|w_{N_3}\|_{X_{\epsilon,\sigma_3}^1}\Bigr]&\qt{if \(s\in[-\tfrac12,-\tfrac16)\)},\end{cases}
\end{align*}
whereas when \(\phi_\epsilon(1)>1\) we have
\begin{align*}
&I[u_{N_1}^\ell,v_{N_2}^\ell,w_{N_3}^H]\\&\lesssim\!\!\begin{cases}\!(\epsilon N_3)^{\frac7{12}-b}(\epsilon N_2)^{\frac1{12}-s-b}\epsilon^{s+b}\|u_{N_1}\|_{\X_{\epsilon,\sigma_1}^{-\frac12,b}}\|v_{N_2}\|_{\X_{\epsilon,\sigma_2}^{s,b}}\!\Bigl[\tfrac1{\<\phi_\epsilon(1)\>^b}\|w_{N_3}\|_{X_{\epsilon,\sigma_3}^1}\Bigr]\text{ if \(s\!\geq\!\tfrac1{12}\!-\!b\)},\\
\!(\epsilon N_3)^{\frac12+s}\epsilon^{\frac14-2s-2b}\|u_{N_1}\|_{\X_{\epsilon,\sigma_1}^{-\frac12,b}}\|v_{N_2}\|_{\X_{\epsilon,\sigma_2}^{s,b}}\!\Bigl[\tfrac1{\<\phi_\epsilon(1)\>^b}\|w_{N_3}\|_{X_{\epsilon,\sigma_3}^1}\Bigr]\text{ if \(s\in[-\tfrac12,\tfrac1{12}-b)\)},\end{cases}
\end{align*}
By decomposing into multiple cases depending on the size of \(\epsilon\) and \(s\) once more, all of these bounds can be seen to be acceptable with one exception: the final estimate loses logarithmically at the endpoint \(s = -\frac12\). We also note that the \(s\geq \frac1{12}-b\) case is the only place we use that the inequality \(b < \frac7{12}\) is strict.

To complete the proof in the case \(s=-\frac12\), we apply Bernstein's inequality and \eqref{X-L2} when \(\phi_\epsilon(1)>1\) to bound
\begin{align*}
I[u_{N_1}^\ell,v_{N_2}^\ell,w_{N_3}^H]&\lesssim N_3\|u_{N_1}^\ell\|_{L_t^\infty L_x^2} \|v_{N_2}^\ell\|_{L_t^\infty L_x^2} \|w_{N_3}^H\|_{L_t^2L_x^\infty}\\
&\lesssim N_3^{\frac32}N_2\epsilon^{-2b}\<M\>^{-1}\|u_{N_1}\|_{\X_{\epsilon,\sigma_1}^{-\frac12,b}}\|v_{N_2}\|_{\X_{\epsilon,\sigma_2}^{s,b}}\Bigl[\tfrac1{\<\phi_\epsilon(1)\>^b}\|w_{N_3}\|_{X_{\epsilon,\sigma_3}^1}\Bigr]\\
&\lesssim N_3^{\frac12}\epsilon^{1-2b}\|u_{N_1}\|_{\X_{\epsilon,\sigma_1}^{-\frac12,b}}\|v_{N_2}\|_{\X_{\epsilon,\sigma_2}^{s,b}}\Bigl[\tfrac1{\<\phi_\epsilon(1)\>^b}\|w_{N_3}\|_{X_{\epsilon,\sigma_3}^1}\Bigr].
\end{align*}
Combining this with our previous estimate, we obtain
\begin{align*}
I[u_{N_1}^\ell,v_{N_2}^\ell,w_{N_3}^H]\lesssim \epsilon^{1-2b}\min\bigl\{N_3^{\frac12},\epsilon^{\frac14}\bigr\}\|u_{N_1}\|_{\X_{\epsilon,\sigma_1}^{-\frac12,b}}\|v_{N_2}\|_{\X_{\epsilon,\sigma_2}^{s,b}}\Bigl[\tfrac1{\<\phi_\epsilon(1)\>^b}\|w_{N_3}\|_{X_{\epsilon,\sigma_3}^1}\Bigr],
\end{align*}
which is now acceptable. This completes the proof of \eqref{bilinear estimates}.
\epf

With suitable linear and nonlinear estimates in hand, we finally apply the contraction mapping theorem to complete the proof of Theorem~\ref{t:LWP-details}:

\bpf[Proof of Theorem~\ref{t:LWP-details}] We start by proving local well-posedness on the unit interval for sufficiently small initial data.

We fix a choice of \(\frac12<b<\frac 7{12}\) and subsequently ignore the dependence of constants on \(b\).  For \(s\geq -\frac12\), \(\epsilon>0\), \(v_0 = (v_0^+,v_0^-)\in H^s\times H^s\), and \(v = (v^+,v^-)\in \X_{\epsilon,+}^{s,b}\times  \X_{\epsilon,-}^{s,b}\), we denote
\[
\Psi^\epsilon(v) := \bigl(\Psi^{\epsilon,+}(v),\Psi^{\epsilon,-}(v)\bigr),
\]
where
\begin{align*}
\Psi^{\epsilon,\sigma}(v) &:= \chi(t)S^{\epsilon,\sigma}(t)v_0^\sigma\\
&\quad + \tfrac{\sigma}3\int_\R \chi(t)\chi\bigl(\tfrac{t'}4\bigr)\sgn(t')S^{\epsilon,\sigma}(t - t')\chi(t')\Lambda_\epsilon^{-1}\p_x\Bigl[\bigl(v^+(t')+v^-(t')\bigr)^2\Bigr]\,dt'\\
&\quad + \tfrac{\sigma}3\int_\R \chi\bigl(\tfrac{t-t'}4\bigr)\sgn(t-t')S^{\epsilon,\sigma}(t - t')\chi(t')\Lambda_\epsilon^{-1}\p_x\Bigl[\bigl(v^+(t')+v^-(t')\bigr)^2\Bigr]\,dt'.
\end{align*}
Applying \eqref{NRG1}, \eqref{NRG2a}, \eqref{NRG2b}, and then \eqref{bilinear estimates}, we may then estimate
\begin{align}
&\|\Psi^\epsilon(v)\|_{\X_{\epsilon,+}^{s,b}\times \X_{\epsilon,-}^{s,b}}\lesssim_s\|v_0\|_{H^s\times H^s} + \<\epsilon\>^{\frac12}\|v\|_{\X_{\epsilon,+}^{-\frac12,b}\times \X_{\epsilon,-}^{-\frac12,b}}\|v\|_{\X_{\epsilon,+}^{s,b}\times \X_{\epsilon,-}^{s,b}},\label{contraction 1}\\
&\|\Psi^\epsilon(v) - \Psi^\epsilon(w)\|_{\X_{\epsilon,+}^{s,b}\times \X_{\epsilon,-}^{s,b}}\label{contraction 2}\\
&\qquad\lesssim_s\<\epsilon\>^{\frac12}\Bigl[\|v\|_{\X_{\epsilon,+}^{-\frac12,b}\times \X_{\epsilon,-}^{-\frac12,b}} + \|w\|_{\X_{\epsilon,+}^{-\frac12,b}\times \X_{\epsilon,-}^{-\frac12,b}}\Bigr] \|v - w\|_{\X_{\epsilon,+}^{s,b}\times \X_{\epsilon,-}^{s,b}}\notag\\
&\qquad\quad + \<\epsilon\>^{\frac12}\Bigl[\|v\|_{\X_{\epsilon,+}^{s,b}\times \X_{\epsilon,-}^{s,b}} + \|w\|_{\X_{\epsilon,+}^{s,b}\times \X_{\epsilon,-}^{s,b}}\Bigr] \|v - w\|_{\X_{\epsilon,+}^{-\frac12,b}\times \X_{\epsilon,-}^{-\frac12,b}}
.\notag
\end{align}

As we wish the lifespan of our solution to be independent of \(s\), we start by considering solutions in \(C([-1,1];H^3\times H^3)\). We already know that local solutions exist and are unique in \(H^3\times H^3\) thanks to \cite{MR954673} and the fact that \eqref{GB-system rescaled} and \eqref{GB-system-2} are equivalent, cf.~\eqref{translate between systems}. Consequently, our main objective here is to quantify the lifespan of these solutions and prove that they have extensions lying in \(\X_{\epsilon,+}^{3,b}\times \X_{\epsilon,-}^{3,b}\).

For \(\epsilon,R,K>0\), we denote the ball
\begin{align*}
\bB_{R,K}^\epsilon &:= \Bigl\{v\in\X_{\epsilon,+}^{3,b}\times  \X_{\epsilon,-}^{3,b}:\|v\|_{\X_{\epsilon,+}^{-\frac12,b}\times  \X_{\epsilon,-}^{-\frac12,b}}\leq R\text{ and }\|v\|_{\X_{\epsilon,+}^{3,b}\times \X_{\epsilon,-}^{3,b}}\leq K\Bigr\}.
\end{align*}
Employing \eqref{contraction 1} and \eqref{contraction 2} for both \(s = -\frac12\) and \(s=3\), we see that if \(C\geq 1\) is chosen sufficiently large and \(0<\eta\ll1\) is chosen sufficiently small (both independently of \(\epsilon\)) then, if \(R_\epsilon = \eta\<\epsilon\>^{-\frac12}\) and \(v_0\in \fB_{R_\epsilon}\cap (H^3\times H^3)\) satisfies \(\|v_0\|_{H^3\times H^3}\leq K_\epsilon:=K\<\epsilon\>^{-\frac12}\) for some \(K>0\), the map \(\Psi^\epsilon\) has a unique fixed point \(v^\epsilon\in \bB_{CR_\epsilon,CK_\epsilon}^\epsilon\).

The estimates \eqref{XY sq sum} and \eqref{X-L2} ensure that \(\X_{\epsilon,\sigma}^{s,b}\subseteq C(\R;H^s)\) and hence the restriction of our fixed point \(v^\epsilon\) to the time interval \([-1,1]\) lies in \(C([-1,1];H^3\times H^3)\) and satisfies \eqref{Duhamel}.

Taking the solution map \(\Upsilon^\epsilon \colon [-1,1]\times \bigl(\fB_{R_\epsilon}\cap(H^3\times H^3)\bigr)\to H^3\times H^3\) to be given by
\[
\Upsilon_t^\epsilon(v_0) := v^\epsilon(t),
\]
we may again apply the estimates \eqref{NRG1}, \eqref{NRG2a}, \eqref{NRG2b}, \eqref{bilinear estimates}, and \eqref{X-L2}, to obtain the estimate
\begin{align}
\sup_{|t|\leq 1}\|\Upsilon_t^\epsilon(v_0) - \Upsilon_t^\epsilon(w_0)\|_{H^{-\frac12}\times H^{-\frac12}}\lesssim\|v_0 - w_0\|_{H^{-\frac12}\times H^{-\frac12}},\label{uni lip from contract}
\end{align}
for all \(v_0,w_0\in \fB_{R_\epsilon}\cap (H^3\times H^3)\). It follows that the solution map has a unique continuous extension to a map from \([-1,1]\times \fB_{R_\epsilon}\to H^{-\frac12}\times H^{-\frac12}\), still satisfying the estimate \eqref{uni lip from contract}.

Our next step is to lay the groundwork for persistence of regularity. Our application of the contraction mapping theorem shows that if \(v_0\in \fB_{R_\epsilon}\) then \((\Psi^\epsilon)^n(0)\to v^\epsilon\) in \(\X_{\epsilon,+}^{-\frac12,b}\times \X_{\epsilon,-}^{-\frac12,b}\).  If \(s>-\frac12\) then an identical contraction mapping argument using \eqref{contraction 1} and \eqref{contraction 2} shows that if we take \(R_{\epsilon,s} = \eta_s\<\epsilon\>^{-\frac12}\) for some possibly smaller \(0<\eta_s\leq \eta\), depending only on \(s\), then for all \(v_0\in \fB_{R_{\epsilon,s}}\cap (H^s\times H^s)\) we have \((\Psi^\epsilon)^n(0)\to v^\epsilon\) in \(\X_{\epsilon,+}^{s,b}\times \X_{\epsilon,-}^{s,b}\). It follows that (the restriction of) the solution map \(\Upsilon^\epsilon\colon [-1,1]\times \bigl(\fB_{R_{\epsilon,s}}\cap (H^s\times H^s)\bigr)\to H^s\times H^s\) satisfies the estimate
\begin{align}
&\sup_{|t|\leq 1}\|\Upsilon_t^\epsilon(v_0) - \Upsilon_t^\epsilon(w_0)\|_{H^s\times H^s}\label{uni lip from contract 2}\\
&\quad\lesssim_s \|v_0 - w_0\|_{H^s\times H^s} + \<\epsilon\>^{\frac12}\Bigl[\|v_0\|_{H^s\times H^s} + \|w_0\|_{H^s\times H^s}\Bigr]\|v_0 - w_0\|_{H^{-\frac12}\times H^{-\frac12}},\notag
\end{align}
for all \(v_0,w_0\in \fB_{R_{\epsilon,s}}\cap (H^s\times H^s)\).

It remains to extend our result to arbitrarily large data.  For \(\lambda>0\), define the spatial rescaling operator
\[
\cR_\lambda[v](x) = \lambda^{-2}v\bigl(\lambda^{-1}x\bigr),
\]
and observe that \(v^\epsilon\in C([-\lambda^{-3},\lambda^{-3}];H^s\times H^s)\) solves \eqref{GB-system-2} for some \(\epsilon>0\) if and only if \(\cR_\lambda[v^\epsilon(\lambda^{-3}t)] \in C([-1,1];H^s\times H^s)\) solves \eqref{GB-system-2} with \(\epsilon\) replaced by \(\lambda\epsilon\).  Crucially, if \(\lambda\geq 1\) we have
\begin{align}
\lambda^{-\frac32 - \max\{s,0\}}\|v\|_{H^s\times H^s}&\lesssim \|\cR_\lambda[v]\|_{H^s\times H^s}\lesssim \lambda^{-\frac32 - \min\{s,0\}}\|v\|_{H^s\times H^s}.\label{rescaling size}
\end{align}

Now, if \(v_0\in \fB_R\) for some \(R>0\) we may choose \(\lambda\simeq \<R\<R\epsilon\>\>\) to ensure that \(\cR_\lambda(v_0)\in \fB_{R_{\lambda\epsilon}}\).  Setting \(T = \lambda^{-3}\), we may then define the solution operator \(\Upsilon^\epsilon \colon [-T,T]\times \fB_R\to H^{-\frac12}\times H^{-\frac12}\) by
\[
\Upsilon_t^\epsilon = \cR_\lambda^{-1}\circ \Upsilon_{\lambda^3 t}^{\lambda\epsilon}\circ \cR_\lambda,
\]
and \eqref{uni lip from contract} yields the estimate \eqref{uniform lipschitz} in the case \(s=-\frac12\).

We note that our construction of \(\Upsilon^\epsilon\) shows that it restricts to a map \(\Upsilon^\epsilon \colon [-T,T]\times \bigl(\fB_R\cap (H^3\times H^3)\bigr)\to H^3\times H^3\). Unconditional uniqueness of solutions of \eqref{GB-system-2} in \(H^3\times H^3\) (see \cite{MR954673} or Lemma~\ref{l:uniqueness} below) and density then suffice to prove that
\eq{GP}{
\Upsilon_{t+s}^\epsilon = \Upsilon_t^\epsilon\circ \Upsilon_s^\epsilon,
}
whenever both sides of this identity are well-defined.

A similar rescaling argument proves that if \(s>-\frac12\) and \(v_0\in \fB_R\cap (H^s\times H^s)\) then there exists a time \(0<T_s = T_s(\lambda)\leq T\) so that \(\Upsilon^\epsilon \colon [-T_s,T_s]\times \bigl(\fB_R\cap (H^s\times H^s)\bigr)\to H^s\times H^s\).  Noting that this time depends only on \(\lambda\) and the size of the data in \(H^{-\frac12}\times H^{-\frac12}\), we may use the group property \eqref{GP} to iterate this argument a finite number of times to deduce that \(\Upsilon^\epsilon \colon [-T,T]\times \bigl(\fB_R\cap (H^s\times H^s)\bigr)\to H^s\times H^s\), i.e., we have persistence of regularity on the whole time interval \([-T,T]\).  The estimate \eqref{uniform lipschitz} for \(s>-\frac12\) now follows from applying \eqref{uni lip from contract 2} and \eqref{rescaling size} to each of these iterates, where we note that when \(0<\epsilon\leq \frac1R\) the time \(T_s\) depends only on \(s\) and \(R\).
\epf

We conclude this section with a simple argument that proves unconditional uniqueness directly for the system \eqref{GB-system-2} whenever \(s>\frac12\). This improves the corresponding statement in \cite{MR954673} that only holds for \(s>\frac52\):
\begin{lem}\label{l:uniqueness}
If \(s>\frac12\), \(v_0 = (v_0^+,v_0^-)\in H^s\times H^s\), and \(I\subseteq \R\) is an interval containing zero then for all \(\epsilon>0\) there exists at most one distributional solution \(v^\epsilon = (v^{\epsilon,+},v^{\epsilon,-})\in C(I;H^s\times H^s)\) of \eqref{GB-system-2} satisfying \(v^\epsilon(0) = v_0\).
\end{lem}
\bpf
Suppose for a contradiction that we have two solutions \(v^\epsilon,w^\epsilon\in C(I;H^s\times H^s)\) satisfying \eqref{GB-system-2} in the sense of distributions with the same initial data \(v^\epsilon(0) = v_0 = w^\epsilon(0)\).

As \(s>\frac12\), the space \(H^s\) is an algebra and hence \(v_t^\epsilon,w_t^\epsilon\in C(I;H^{s-2}\times H^{s-2})\). Taking
\[
V^\epsilon := v^\epsilon - w^\epsilon,
\]
we see that \(V^\epsilon\in C(I;H^s\times H^s)\cap C^1(I;H^{s-2}\times H^{s-2})\) satisfies
\eq{V eqn}{
V_t^{\epsilon,\sigma} = - \tfrac\sigma{\epsilon^2}\Lambda_\epsilon V_x^{\epsilon,\sigma}+\tfrac{2\sigma}3\Lambda_\epsilon^{-1}\p_x\bigl[(v^{\epsilon,+}+v^{\epsilon,-} + w^{\epsilon,+}+w^{\epsilon,-})(V^{\epsilon,+}+V^{\epsilon,-})\bigr].
}
We will prove that for all \(0<T\in I\) and all \(0\leq t\leq T\) we have \(V^\epsilon(t) = 0\). Applying a parallel argument for times \(0>T\in I\) then yields the desired result.

For \(t\in (0,T)\), we use \eqref{V eqn} and the fact that \(H^s\) is an algebra to compute that
\begin{align*}
\Bigl|\tfrac d{dt}\|V_{\leq N}^{\epsilon,\sigma}\|_{H^s}^2\Bigr|&\lesssim \|V_{\leq N}^{\epsilon,\sigma}\|_{H^s}\bigl\|\Lambda_\epsilon^{-1}\p_x\bigl[(v^{\epsilon,+}+v^{\epsilon,-} + w^{\epsilon,+}+w^{\epsilon,-})(V^{\epsilon,+}+V^{\epsilon,-})\bigr]\bigr\|_{H^s}\\
&\lesssim \tfrac1{\epsilon} C_T \|V^\epsilon\|_{H^s\times H^s}^2,
\end{align*}
where
\[
C_T := \sup_{0\leq t\leq T}\Bigl[\|v^\epsilon(t)\|_{H^s\times H^s} + \|w^\epsilon(t)\|_{H^s\times H^s}\Bigr].
\]

Integrating in time and then applying Fatou's lemma, we see that for all \(0\leq t\leq T\) we have
\begin{align*}
\|V^\epsilon(t)\|_{H^s\times H^s}^2 \leq \liminf_{N\to\infty}\Bigl[\|V_{\leq N}^{\epsilon,+}(t)\|_{H^s}^2 + \|V_{\leq N}^{\epsilon,-}(t)\|_{H^s}^2 \Bigr]\lesssim \tfrac1\epsilon C_T \int_0^t \|V^\epsilon(t')\|_{H^s\times H^s}^2\,dt'.
\end{align*}
The integral form of Gronwall's inequality now proves that \(V^\epsilon(t) = 0\) for all \(0\leq t\leq T\).
\epf

\section{A priori estimates}\label{s:AP}

In this section we first prove Theorem~\ref{t:AP} and then prove Theorems~\ref{t:GWP} and~\ref{t:equicontinuity}.
Recalling the definition of \(L_0\) from \eqref{free Lax}, for real-valued \(|\kappa|\geq 2\) we define the operator \(R_0(\kappa) := (iL_0 + i\kappa^3 - i\kappa)^{-1}\) so that
\begin{align*}
R_0(\kappa)f &= \tfrac1{\sqrt{2\pi}}\int_\R \tfrac1{\xi^3 + \xi + i\kappa^3 - i\kappa}\hat f(\xi)e^{ix\xi}\,d\xi,
\end{align*}
and then take
\begin{align*}
\sqrt{R_0(\kappa)}f &:= \tfrac1{\sqrt{2\pi}}\int_\R \tfrac1{\sqrt{\xi^3 + \xi + i\kappa^3 - i\kappa}}\hat f(\xi)e^{ix\xi}\,d\xi,
\end{align*}
where we use the branch of the square root for which \(\sqrt z>0\) whenever \(z>0\).

We start with the following lemma:
\begin{lem}\label{l:basic bounds}
For real-valued \(|\kappa|\geq 2\) and complex-valued \(f\in \Schwartz(\R)\) we have the estimates
\begin{align}
\Bigl\|\sqrt{R_0(\kappa)}f \sqrt{R_0(\kappa)}\Bigr\|_{\I_2}&\lesssim \tfrac1{|\kappa|}\|f\|_{H_\kappa^{-\frac32}},\label{H-3/2}\\
\Bigl\|\sqrt{R_0(\kappa)}(f\p_x + \p_xf)\sqrt{R_0(\kappa)}\Bigr\|_{\I_2}&\lesssim \tfrac1{|\kappa|}\|f\|_{H_\kappa^{-\frac12}}.\label{H-1/2}
\end{align}
\end{lem}
\bpf
Arguing as in \cite[Lemma 4.1]{MR3820439}, we compute that for \(f\in \Schwartz\) we have
\begin{align*}
\Bigl\|\sqrt{R_0(\kappa)}f \sqrt{R_0(\kappa)}\Bigr\|_{\I_2}^2=\tfrac1{2\pi}\iint_{\R^2} \tfrac{|\hat f(\xi-\eta)|^2}{\sqrt{(\xi^3 + \xi)^2 + (\kappa^3 - \kappa)^2}\sqrt{(\eta^3 + \eta)^2 + (\kappa^3 - \kappa)^2}}\,d\xi\,d\eta.
\end{align*}
For \(|\kappa|\geq 2\) we have
\[
\sqrt{(\xi^3 + \xi)^2 + (\kappa^3 - \kappa)^2}\simeq \bigl(|\kappa|^2 + |\xi|^2\bigr)^{\frac32},
\]
and hence
\begin{align*}
\Bigl\|\sqrt{R_0(\kappa)}f\sqrt{R_0(\kappa)}\Bigr\|_{\I_2}^2\simeq \int_\R \tfrac1{|\kappa|^5}m_0\bigl(\tfrac\xi{2|\kappa|}\bigr)|\hat f(\xi)|^2\,d\xi,
\end{align*}
where
\[
m_0(\xi) := \int_\R \tfrac1{\<\eta-\xi\>^3\<\xi+\eta\>^3}\,d\eta.
\]

Similarly, we may estimate
\begin{align*}
\Bigl\|\sqrt{R_0(\kappa)}(f\p_x + \p_xf) \sqrt{R_0(\kappa)}\Bigr\|_{\I_2}^2&=\tfrac1{2\pi}\iint_{\R^2} \tfrac{(\xi+\eta)^2|\hat f(\xi-\eta)|^2}{\sqrt{(\xi^3 + \xi)^2 + (\kappa^3 - \kappa)^2}\sqrt{(\eta^3 + \eta)^2 + (\kappa^3 - \kappa)^2}}\,d\xi\,d\eta\\
&\simeq \int_\R \tfrac1{|\kappa|^3}m_1\bigl(\tfrac \xi{2|\kappa|}\bigr)|\hat f(\xi)|^2\,d\xi,
\end{align*}
where
\[
m_1(\xi) := \int_\R \tfrac{\eta^2}{\<\eta-\xi\>^3\<\xi+\eta\>^3}\,d\eta.
\]

Both of the estimates \eqref{H-3/2} and \eqref{H-1/2} now follow from the fact that for \(j=0,1\) we have
\begin{align*}
m_j(\xi) &= 2\int_0^{2\<\xi\>}\frac{\eta^{2j}}{\<\eta-\xi\>^3\<\xi+\eta\>^3}\,d\eta + 2\int_{2\<\xi\>}^\infty\frac{\eta^{2j}}{\<\eta-\xi\>^3\<\xi+\eta\>^3}\,d\eta\\
&\lesssim \<\xi\>^{2j-3}\int_0^{2\<\xi\>}\frac1{\<\eta - |\xi|\>^3}\,d\eta + \int_{2\<\xi\>}^\infty \eta^{2j-6}\,d\eta\lesssim \<\xi\>^{2j-3}.\qedhere
\end{align*}
\epf

Taking \(L\) to be defined as in \eqref{good Lax pair}, for real-valued \(|\kappa|\geq 2\) we define the operator
\begin{align}
A(\kappa) &:= i\sqrt{R_0(\kappa)}(L - L_0)\sqrt{R_0(\kappa)}\label{d:A}\\
&= i\sqrt{R_0(\kappa)}(q\p_x+\p_xq + r)\sqrt{R_0(\kappa)}.\notag
\end{align}
We then have the following:

\begin{lem}\label{l:I bounds}
For complex-valued \((q,r)\in H^{-\frac12}\times H^{-\frac32}\), we have the estimate
\begin{align}
\|A(\kappa)\|_{\I_2} &\lesssim \tfrac 1{|\kappa|}\|(q,r)\|_{H_\kappa^{-\frac12}\times H_\kappa^{-\frac32}},\label{A I2}
\end{align}
uniformly for all real-valued \(|\kappa|\geq 2\).

Further, if \(\sigma\geq-\frac12\) and \((q,r)\in H^\sigma\times H^{\sigma-1}\) then we have the estimate
\begin{align}
\|A(\kappa)\|_{\I_4}^4\lesssim |\kappa|^{-7+\max\{2-2\sigma,0\}}\|(q,r)\|_{H^{-\frac12}\times H^{-\frac32}}^2\|(q,r)\|_{H^\sigma\times H^{\sigma-1}}^2,\label{A I4}
\end{align}
uniformly for all real-valued \(|\kappa|\geq 2\).
\end{lem}
\bpf
The estimate \eqref{A I2} follows directly from applying the estimates \eqref{H-3/2} and \eqref{H-1/2} to the definition \eqref{d:A}.

For \eqref{A I4}, we denote \(A = A(\kappa;q,r)\) and \(A_N = A(\kappa;q_N,r_N)\), where we recall our convention that \((q_1,r_1) = (P_{\leq 1}q,P_{\leq 1}r)\) and \((q_N,r_N) = (P_Nq,P_Nr)\) for dyadic \(N\geq 2\).

If \(\sigma=-\frac12\) then we apply \eqref{A I2} to bound
\[
\|A\|_{\I_4}^4\leq\|A\|_{\I_2}^4\lesssim \tfrac 1{|\kappa|^4}\|(q,r)\|_{H_\kappa^{-\frac12}\times H_\kappa^{-\frac32}}^4\lesssim \tfrac 1{|\kappa|^4}\|(q,r)\|_{H^{-\frac12}\times H^{-\frac32}}^4.
\]

If \(\sigma>-\frac12\) then we first apply \eqref{A I2} to obtain
\eq{freq loc I2}{
\|A_N\|_{\I_2} \lesssim \tfrac1{|\kappa|N^{\sigma}\sqrt{|\kappa|+N}}\|(q_N,r_N)\|_{H^\sigma\times H^{\sigma-1}}.
}
Next, we apply Bernstein's inequality to estimate
\begin{align*}
\|A_N\|_\op\lesssim \tfrac1{|\kappa|^2}\|q_N\|_{L^\infty} + \tfrac1{|\kappa|^3}\|r_N\|_{L^\infty}\lesssim \tfrac{(|\kappa| + N)N}{|\kappa|^3}\|(q_N,r_N)\|_{H^{-\frac12}\times H^{-\frac32}}.
\end{align*}
Using this estimate when \(N\leq |\kappa|\) and \eqref{A I2} when \(N> |\kappa|\) gives us
\eq{freq loc op}{
\|A_N\|_\op \lesssim \tfrac{N}{|\kappa|(|\kappa| + N)}\|(q_N,r_N)\|_{H^{-\frac12}\times  H^{-\frac32}}.
}
We now apply \eqref{freq loc I2}, \eqref{freq loc op}, and the Cauchy--Schwarz inequality to estimate
\begin{align*}
\|A\|_{\I_4}^4 &= \sum_{N_1,N_2,N_3,N_4\geq 1}\tr\bigl\{A_{N_1}^*A_{N_2}A_{N_3}^*A_{N_4}\bigr\}\\
&\lesssim \sum_{1\leq N_1\leq N_2\leq N_3\simeq N_4}\|A_{N_1}\|_\op\|A_{N_2}\|_\op\|A_{N_3}\|_{\I_2}\|A_{N_4}\|_{\I_2}\\
&\lesssim \|(q,r)\|_{H^{-\frac12}\times H^{-\frac32}}^2\|(q,r)\|_{H^\sigma\times H^{\sigma-1}}^2\\
&\quad\times \sup_{M\geq 1}\left[\sum_{1\leq N_1\leq N_2\leq M}\tfrac1{|\kappa|^8}\tfrac{N_1^2N_2^2}{(|\kappa| + N_1)^2(|\kappa| + N_2)^2(|\kappa|+M)^2M^{4\sigma}}\right]^{\frac12}\\
&\lesssim \|(q,r)\|_{H^{-\frac12}\times H^{-\frac32}}^2\|(q,r)\|_{H^\sigma\times H^{\sigma-1}}^2\times\sup_{M\geq 1}\tfrac1{|\kappa|^4}\tfrac{M^{2(1-\sigma)}}{(|\kappa| + M)^3}\log\bigl(2+\tfrac{M}{|\kappa|}\bigr),
\end{align*}
from which we obtain \eqref{A I4}.
\epf

We are now in a position to rigorously construct the perturbation determinant \(\alpha(\kappa)\) that we formally defined in \eqref{perturbation determinant formal}. Recalling the definition \eqref{balls} of the ball \(B_\delta\), we have the following:
\begin{prop}\label{p:alpha}
There exists \(\delta_0>0\) so that for all real-valued \(|\kappa|\geq 2\) the series
\eq{perturbation determinant}{
\alpha(\kappa) := \sum_{\ell=2}^\infty \tfrac{(-1)^\ell}\ell \tr\bigl\{A(\kappa)^\ell\bigr\}
}
defines an analytic function from \(B_{\delta_0}\to \C\).

Adopting the notation from \cite{MR4726498}, we denote
\begin{align*}
\alpha\sprack \ell(\kappa) := \tfrac{(-1)^\ell}\ell\tr\bigl\{A(\kappa)^\ell\bigr\}\qtq{and}\alpha\sprack{\geq L}(\kappa) := \sum_{\ell=L}^\infty \alpha\sprack \ell(\kappa).
\end{align*}
For \(\kappa\geq 2\) and real-valued \((q,r)\in H^{-\frac12}\times H^{-\frac32}\) we then have the identity
\begin{align}
\alpha\sprack2(\kappa) + \alpha\sprack2(-\kappa) &= \tfrac{2\kappa}{3\kappa^2-1}\int_\R \tfrac{(4+\xi^2)|\hat q(\xi)|^2 + 3|\hat r(\xi)|^2}{\xi^4 + (3\kappa^2+2)\xi^2 + (3\kappa^2-1)^2}\,d\xi.\label{tr A2}
\end{align}
Moreover, if \(\sigma\geq -\frac12\) and \((q,r)\in (H^\sigma\times H^{\sigma-1})\cap B_{\delta_0}\) are real-valued then we have the estimate
\begin{align}
\left|\alpha\sprack{\geq 3}(\kappa) + \alpha\sprack{\geq 3}(-\kappa)\right|\lesssim \kappa^{-5+\max\{1-2\sigma,0\}}\delta_0\|(q,r)\|_{H^\sigma\times H^{\sigma-1}}^2,\label{tr A3}
\end{align}
whereas if \(\sigma\geq \frac12\) then
\begin{align}
&\left|\alpha\sprack{\geq 3}(\kappa) + \alpha\sprack{\geq 3}(-\kappa) + \tfrac{16\kappa}{3(3\kappa^2-1)^3}\int_\R q^3\,dx\right| \label{tr A3 high reg}\\
&\qquad\lesssim \kappa^{-7+\max\{3-2\sigma,0\}}\delta_0\|(q,r)\|_{H^\sigma\times H^{\sigma-1}}^2,\notag
\end{align}
uniformly for all \(\kappa\geq 2\).

Finally, if \((q,r)\in (H^3\times H^2)\cap B_{\delta_0}\) are real-valued, the Hamiltonian \(\cH\) is defined as in \eqref{energy}, and the Poisson bracket \(\{\cdot,\cdot\}\) is defined as in \eqref{Poisson bracket} then for all real-valued \(|\kappa|\geq 2\) we have
\eq{Poisson commute}{
\bigl\{\alpha(\kappa),\cH\bigr\} = 0.
}
\end{prop}
\bpf
From \eqref{A I2}, for each \(\ell\geq 2\) and real-valued \(|\kappa|\geq 2\) the map \(H^{-\frac12}\times H^{-\frac32}\ni(q,r)\mapsto\alpha\sprack \ell\in \C\) is an \(\ell\)-homogeneous polynomial satisfying
\begin{align*}
\bigl|\alpha\sprack \ell\bigr|\leq \tfrac1\ell \|A\|_{\I_2}^\ell \leq C^\ell\|(q,r)\|_{H^{-\frac12}\times H^{-\frac32}}^\ell,
\end{align*}
for some \(C>0\). It follows that if \(\delta_0>0\) is chosen sufficiently small then for all real-valued \(|\kappa|\geq 2\) the series \eqref{perturbation determinant} converges to an analytic function from \(B_{\delta_0}\to \C\).

To prove \eqref{tr A2}, we compute that for real-valued \((q,r)\) we have
\begin{align*}
\tr\bigl\{A^2\bigr\} &= \tfrac1{2\pi}\iint_{\R^2} \tfrac{(\xi+\eta)^2|\hat q(\xi-\eta)|^2 - |\hat r(\xi-\eta)|^2 - 2i(\xi+\eta)\Re[\hat q(\xi-\eta)\overline{\hat r(\xi-\eta)}]}{(\xi^3 + \xi + i\kappa^3 - i\kappa)(\eta^3 + \eta + i\kappa^3 - i\kappa)}\,d\xi\,d\eta\\
&= \tfrac1{2\pi}\iint_{\R^2} \tfrac{(\xi+2\eta)^2|\hat q(\xi)|^2 - |\hat r(\xi)|^2 - 2i(\xi+2\eta)\Re[\hat q(\xi)\overline{\hat r(\xi)}]}{\left((\xi+\eta)^3 + (\xi+\eta) + i\kappa^3 - i\kappa\right)(\eta^3 + \eta + i\kappa^3 - i\kappa)}\,d\xi\,d\eta.
\end{align*}
We may readily factor
\[
\eta^3 + \eta + i\kappa^3 - i\kappa = (\eta - i\kappa)(\eta^2 + i\kappa \eta - \kappa^2 + 1).
\]
If \(|\kappa|\geq 2\) this polynomial has exactly one simple root at \(\eta = i\kappa\) in either the upper (when \(\kappa\geq 2\)) or lower  (when \(\kappa\leq -2\)) half plane. By contour integration, we obtain
\begin{align*}
\tr\bigl\{A^2\bigr\}= \tfrac{2\sgn(\kappa)}{3\kappa^2-1}\int_\R \tfrac{\kappa(4+\xi^2)|\hat q(\xi)|^2 + 3\kappa|\hat r(\xi)|^2 - 2(3\kappa^2 + 1 + \xi^2)\Re[\hat q(\xi)\overline{\hat r(\xi)}]}{\xi^4 + (3\kappa^2+2)\xi^2 + (3\kappa^2-1)^2}\,d\xi,
\end{align*}
from which \eqref{tr A2} readily follows.

For \eqref{tr A3} and \eqref{tr A3 high reg}, we start by decomposing
\[
\alpha\sprack{\geq 3}(\kappa) = \alpha\sprack 3(\kappa) + \alpha\sprack{\geq 4}(\kappa).
\]
To handle the second of these terms, we apply \eqref{A I2} and \eqref{A I4} to estimate
\begin{align}
\bigl|\alpha\sprack{\geq 4}\bigr| &\leq \sum_{\ell=4}^\infty \|A\|_{\I_2}^{\ell-4}\|A\|_{\I_4}^4\label{a4}\\
&\lesssim \sum_{\ell=4}^\infty C^{\ell-4}|\kappa|^{-7+\max\{2-2\sigma,0\}}\|(q,r)\|_{H^{-\frac12}\times H^{-\frac32}}^{\ell-2}\|(q,r)\|_{H^\sigma\times  H^{\sigma-1}}^2\notag\\
&\lesssim|\kappa|^{-7+\max\{2-2\sigma,0\}}\delta_0^2\|(q,r)\|_{H^\sigma\times  H^{\sigma-1}}^2.\notag
\end{align}

To control the contribution of \(\alpha\sprack 3(\kappa)\), we introduce
\[
B(\kappa;q) := A(\kappa;q,0)\qtq{and}C(\kappa;r) := A(\kappa;0,r),
\]
so that \(A = B + C\). We then compute that
\begin{align*}
\tr\{B^2C\} &= \tfrac{i}{(2\pi)^{3/2}}\iiint_{\R^3} \tfrac{(\xi+\eta)(\eta+\zeta)\hat q(\xi-\eta)\hat q(\eta-\zeta)\hat r(\zeta - \xi)}{\left(\xi^3 + \xi + i\kappa^3 - i\kappa\right)\left(\eta^3 + \eta + i\kappa^3 - i\kappa\right)\left(\zeta^3 + \zeta + i\kappa^3 - i\kappa\right)} \,d\xi\,d\eta\,d\zeta,\\
\tr\{C^3\} &= \tfrac1{i(2\pi)^{3/2}}\iiint_{\R^3} \tfrac{\hat r(\xi-\eta)\hat r(\eta-\zeta)\hat r(\zeta - \xi)}{\left(\xi^3 + \xi + i\kappa^3 - i\kappa\right)\left(\eta^3 + \eta + i\kappa^3 - i\kappa\right)\left(\zeta^3 + \zeta + i\kappa^3 - i\kappa\right)} \,d\xi\,d\eta\,d\zeta.
\end{align*}
By making the change of variables \((\xi,\eta,\zeta)\mapsto(-\zeta,-\eta,-\xi)\), we see that both of these traces are odd in \(\kappa\). Conversely, the traces
\begin{align*}
\tr\{B^3\} &= -\tfrac1{(2\pi)^{3/2}}\iiint_{\R^3} \tfrac{(\xi+\eta)(\eta+\zeta)(\zeta+\xi)\hat q(\xi-\eta)\hat q(\eta-\zeta) \hat q(\zeta - \xi)}{\left(\xi^3 + \xi + i\kappa^3 - i\kappa\right)\left(\eta^3 + \eta + i\kappa^3 - i\kappa\right)\left(\zeta^3 + \zeta + i\kappa^3 - i\kappa\right)} \,d\xi\,d\eta\,d\zeta,\\
\tr\{C^2B\} &= \tfrac1{(2\pi)^{3/2}}\iiint_{\R^3} \tfrac{(\zeta+\xi)\hat r(\xi-\eta) \hat r(\eta - \zeta)\hat q(\zeta-\xi) }{\left(\xi^3 + \xi + i\kappa^3 - i\kappa\right)\left(\eta^3 + \eta + i\kappa^3 - i\kappa\right)\left(\zeta^3 + \zeta + i\kappa^3 - i\kappa\right)} \,d\xi\,d\eta\,d\zeta,
\end{align*}
are even in \(\kappa\). It follows that
\begin{align}
\alpha\sprack 3(\kappa) + \alpha\sprack 3(-\kappa) &= -\tfrac 23\tr\{B(\kappa)^3\} - 2\tr\{C(\kappa)^2B(\kappa)\}.\label{a3 decomp}
\end{align}

For \(\kappa\geq 2\) we may use contour integration again to compute that
\begin{align*}
\tr\{B^3\}
&= - \tfrac i{\sqrt{2\pi}(3\kappa^2 - 1)}\iint_{\R^2} \tfrac{b(\xi,\eta,\kappa)\hat q(\xi-\eta)\hat q(\eta) \overline{\hat q(\xi)}}{a(\xi,\kappa)a(\eta,\kappa)a(\xi - \eta,\kappa)} \,d\xi\,d\eta,\\
\tr\{C^2B\}
& = \tfrac i{\sqrt{2\pi}(3\kappa^2 - 1)}\iint_{\R^2} \tfrac{c(\xi,\eta,\kappa)\hat r(\xi-\eta) \hat r(\eta)\overline{\hat q(\xi)} }{a(\xi,\kappa)a(\eta,\kappa)a(\xi - \eta,\kappa)} \,d\xi\,d\eta,
\end{align*}
where the polynomial
\[
a(\xi,\kappa) := \xi^4 + (3\kappa^2+2)\xi^2 + (3\kappa^2-1)^2,
\]
appears in the denominator of \RHS{tr A2} and the polynomials 
\begin{align*}
b(\xi,\eta,\kappa) &:=648i\kappa^9 + 54i\kappa^7(\xi^2 +\eta^2 + \zeta^2 - 24) - 9\kappa^6(\xi+\eta)(\xi+\zeta)(\eta - \zeta)\\
&\quad - 144i\kappa^5(\xi^2+\eta^2+\zeta^2-3)\\
&\quad + \tfrac32\kappa^4(\xi+\eta)(\xi+\zeta)(\eta - \zeta)(\xi^2+\eta^2+\zeta^2+26)\\
&\quad -6i\kappa^3\bigl[\xi^2\eta^2\zeta^2 + 2(\xi^4+\eta^4+\zeta^4) - 5(\xi^2 + \eta^2 + \zeta^2) - 8\bigr]\\
&\quad + \tfrac32\kappa^2(\xi+\eta)(\xi+\zeta)(\eta-\zeta)\bigl[\xi^4+\eta^4 + \zeta^4 + 9(\xi^2 + \eta^2 + \zeta^2) + 22\bigr]\\
&\quad - \tfrac12i\kappa\bigl[\xi^2\eta^2\zeta^2(\xi^2+\eta^2+\zeta^2) + 4(\xi^6+\eta^6+\zeta^6) - 18\xi^2\eta^2\zeta^2\\
&\quad\qquad\quad+ 20(\xi^4+\eta^4+\zeta^4) + 40(\xi^2+\eta^2+\zeta^2) + 48\bigr] \\
&\quad + (\xi+\eta)(\xi+\zeta)(\eta-\zeta)(1+\xi^2)(1+\eta^2)(1+\zeta^2),\\
c(\xi,\eta,\kappa) &:= 324i\kappa^7 - 135\kappa^6(\eta-\zeta)+108i\kappa^5(\xi^2+\eta^2+\zeta^2)\\
&\quad -\tfrac{27}2\kappa^4(\eta-\zeta)(\xi^2+\eta^2+\zeta^2-6)\\
&\quad +18i\kappa^3\bigl[2\xi^4-5\xi^2\eta\zeta - 2(\xi^2+\eta^2+\zeta^2)-6\bigr]\\
&\quad  - \tfrac32\kappa^2(\eta-\zeta)\bigl[\xi^4+\eta^4+\zeta^4 + \xi^2+\eta^2+\zeta^2 + 6\bigr]\\
&\quad +3i\kappa\bigl[2\xi^6 - 7\xi^4\eta\zeta+\xi^2\eta^2\zeta^2 + 4\xi^4-5\xi^2\eta\zeta\\
&\quad\qquad\quad + 2(\eta^4+\zeta^4)+6(\xi^2+\eta^2+\zeta^2)+8\bigr]\\
&\quad - (\eta-\zeta)(1+\xi^2)(1+\eta^2)(1+\zeta^2),
\end{align*}
where \(\zeta := \xi-\eta\).

If \(\<\xi\>\simeq N_1\), \(\<\eta\>\simeq N_2\), \(\<\zeta\>\simeq N_3\), and we denote \(N_{\max} = \max\{N_1,N_2,N_3\}\), \(N_{\min} = \min\{N_1,N_2,N_3\}\) then we have the estimates
\begin{align}
|b(\xi,\eta,\kappa)| &\lesssim (\kappa+N_{\min})^2(\kappa + N_{\max})^7,\label{b estimate}\\
|b(\xi,\eta,\kappa) - 648i\kappa^9| &\lesssim (\kappa+N_{\min})^2(\kappa + N_{\max})^5N_{\max}^2 ,\label{b high reg estimate}\\
|c(\xi,\eta,\kappa)| &\lesssim (\kappa + N_{\min})^2(\kappa+N_{\max})^5+\kappa N_1^4(N_1^2+N_2N_3).\label{c estimate}
\end{align}

Similar to the proof of Lemma~\ref{l:I bounds}, we take \(B_N = B(\kappa;q_N)\), \(C_N = C(\kappa;r_N)\), and then use H\"older's inequality, Young's convolution inequality, Plancherel's Theorem, and the estimate \eqref{b estimate} to bound
\begin{align}
\bigl|\tr\{B^3\}\bigr| &\lesssim \sum_{1\leq N_1\leq N_2\simeq N_3}\tfrac1{\kappa^2(\kappa+N_1)^2(\kappa+N_2)}\|\hat q_{N_1}\|_{L^1}\|\hat q_{N_2}\|_{L^2}\|\hat q_{N_3}\|_{L^2}\label{tr B3}\\
&\lesssim \sum_{1\leq N_1\leq N_2\simeq N_3} \tfrac{N_1N_2^{-2\sigma}}{\kappa^2(\kappa+N_1)^2(\kappa+N_2)}\|q_{N_1}\|_{H^{-\frac12}}\|q_{N_2}\|_{H^\sigma}\|q_{N_3}\|_{H^\sigma}\notag\\
&\lesssim \sum_{1\leq N_2\simeq N_3}\tfrac{N_2^{1-2\sigma}}{\kappa^3(\kappa+N_2)^2}\|q\|_{H^{-\frac12}}\|q_{N_2}\|_{H^\sigma}\|q_{N_3}\|_{H^{\sigma}}\notag\\
&\lesssim \kappa^{-5 + \max\{1-2\sigma,0\}}\|q\|_{H^{-\frac12}}\|q\|_{H^\sigma}^2,\notag
\end{align}
for all \(\sigma\geq -\frac12\). Employing an identical argument with \eqref{b estimate} replaced by \eqref{b high reg estimate} yields the estimate
\begin{align}
\biggl|\tr\{B^3\} &- \tfrac{648\kappa^9}{\sqrt{2\pi}(3\kappa^2-1)}\iint_{\R^2}\tfrac{\hat q(\xi-\eta)\hat q(\eta) \overline{\hat q(\xi)}}{a(\xi,\kappa)a(\eta,\kappa)a(\xi - \eta,\kappa)} \,d\xi\,d\eta\biggr|\label{tr B3 high reg}\\
&\lesssim \kappa^{-7 + \max\{3-2\sigma,0\}}\|q\|_{H^{-\frac12}}\|q\|_{H^\sigma}^2,\notag
\end{align}
for all \(\sigma\geq -\frac12\). Finally, a similar argument using \eqref{c estimate} yields
\begin{align}
\bigl|\tr\{C^2B\}\bigr| &\lesssim \sum_{1\leq N_2\leq N_3\simeq N_1}\tfrac1{\kappa^2(\kappa+N_2)^3(\kappa+N_1)^2}\|\hat q_{N_1}\|_{L^2}\|\hat r_{N_2}\|_{L^1}\|\hat r_{N_3}\|_{L^2}\label{tr B2C}\\
&\quad + \sum_{1\leq N_1\leq N_2\simeq N_3}\tfrac1{\kappa^2(\kappa+N_1)^2(\kappa+N_2)^3}\|\hat q_{N_1}\|_{L^1}\|\hat r_{N_2}\|_{L^2}\|\hat r_{N_3}\|_{L^2}\notag\\
&\lesssim \sum_{1\leq N_2\leq N_3\simeq N_1}\tfrac{N_2^2N_1^{1-2\sigma}}{\kappa^2(\kappa+N_2)^3(\kappa+N_1)^2}\|q_{N_1}\|_{H^\sigma}\|r_{N_2}\|_{H^{-\frac32}}\|r_{N_3}\|_{H^{\sigma-1}}\notag\\
&\quad + \sum_{1\leq N_1\leq N_2\simeq N_3}\tfrac{N_1N_2^{2-2\sigma}}{\kappa^2(\kappa+N_1)^2(\kappa+N_2)^3}\|q_{N_1}\|_{H^{-\frac12}}\|r_{N_2}\|_{H^{\sigma-1}}\|r_{N_3}\|_{H^{\sigma-1}}\notag\\
&\lesssim \kappa^{-7 + \max\{3-2\sigma,0\}}\|(q,r)\|_{H^{-\frac12}\times H^{-\frac32}}\|(q,r)\|_{H^\sigma\times H^{\sigma-1}}\|r\|_{H^{\sigma-1}},\notag
\end{align}
for all \(\sigma\geq -\frac12\). The estimate \eqref{tr A3} now follows from combining \eqref{a4}, \eqref{a3 decomp}, \eqref{tr B3}, and \eqref{tr B2C}.

To complete the proof of \eqref{tr A3 high reg}, we compute that
\[
\tfrac{648\kappa^9}{(3\kappa^2-1)a(\xi,\kappa)a(\eta,\kappa)a(\zeta,\kappa)} = \tfrac{8\kappa}{(3\kappa^2-1)^3} + \tfrac{d(\xi,\eta,\kappa)}{(3\kappa^2-1)^3a(\xi,\kappa)a(\eta,\kappa)a(\zeta,\kappa)},
\]
where
\begin{align*}
d(\xi,\eta,\kappa) & := -1944\kappa^{11}(\xi^2+\eta^2+\zeta^2-4)\\
&\quad -324\kappa^9\bigl[3(\xi^4+\eta^4+\zeta^4) - 4(\xi^2+\eta^2+\zeta^2) + 28\bigr]\\
&\quad -216\kappa^7\bigl[\xi^6+\eta^6+\zeta^6 - 5\xi^2\eta^2\zeta^2 - 3(\xi^4+\eta^4+\zeta^4) -2(\xi^2+\eta^2+\zeta^2)-20\bigr]\\
&\quad -36\kappa^5\bigl[\xi^8+\eta^8+\zeta^8-6\xi^2\eta^2\zeta^2(\xi^2+\eta^2+\zeta^2) + 12\xi^2\eta^2\zeta^2\\
&\quad\qquad\qquad + 9(\xi^4+\eta^4+\zeta^4) + 16(\xi^2+\eta^2+\zeta^2)+30\bigr]\\
&\quad - 12\kappa^3\bigl[\xi^2\eta^2\zeta^2(\xi^4+\eta^4+\zeta^4)-2(\xi^8+\eta^8+\zeta^8)\\
&\quad\qquad\qquad+24\xi^2\eta^2\zeta^2(\xi^2+\eta^2+\zeta^2)-6(\xi^6+\eta^6+\zeta^6)\\
&\quad\qquad\qquad+60\xi^2\eta^2\zeta^2-12(\xi^4+\eta^4+\zeta^4)-14(\xi^2+\eta^2+\zeta^2)-12\bigr]\\
&\quad - 8\kappa(1+\xi^2)^2(1+\eta^2)^2(1+\zeta^2)^2.
\end{align*}
If \(\<\xi\>\simeq N_1\), \(\<\eta\>\simeq N_2\), \(\<\zeta\>\simeq N_3\), \(N_{\max} = \max\{N_1,N_2,N_3\}\), and \(N_{\min} = \min\{N_1,N_2,N_3\}\), we may bound
\[
|d(\xi,\eta,\kappa)| \lesssim \kappa(\kappa+N_{\min})^4(\kappa+N_{\max})^6N_{\max}^2.
\]
Employing an identical argument to before, we then have
\begin{align}
&\left|\tfrac{648\kappa^9}{\sqrt{2\pi}(3\kappa^2-1)}\iint_{\R^2}\tfrac{\hat q(\xi-\eta)\hat q(\eta) \overline{\hat q(\xi)}}{a(\xi,\kappa)a(\eta,\kappa)a(\xi - \eta,\kappa)} \,d\xi\,d\eta - \tfrac{8\kappa}{(3\kappa^2-1)^3}\int_\R q^3\,dx\right|\label{a3 lot}\\
&\qquad=\left|\tfrac1{\sqrt{2\pi}(3\kappa^2-1)^3}\iint_{\R^2}\tfrac{d(\xi,\eta,\kappa)\hat q(\xi-\eta)\hat q(\eta) \overline{\hat q(\xi)}}{a(\xi,\kappa)a(\eta,\kappa)a(\xi - \eta,\kappa)} \,d\xi\,d\eta\right|\notag\\
&\qquad\lesssim\sum_{1\leq N_1\leq N_2\simeq N_3} \tfrac{N_2^2}{\kappa^5(\kappa+N_2)^2}\|\hat q_{N_1}\|_{L^1}\|\hat q_{N_2}\|_{L^2}\|\hat q_{N_3}\|_{L^2}\notag\\
&\qquad\lesssim \sum_{1\leq N_1\leq N_2\simeq N_3} \tfrac{N_2^{2(1-\sigma)}N_1}{\kappa^5(\kappa+N_2)^2}\|q_{N_1}\|_{H^{-\frac12}}\|q_{N_2}\|_{H^\sigma}\|q_{N_3}\|_{H^\sigma}\notag\\
&\qquad\lesssim \kappa^{-7+\max\{3-2\sigma,0\}}\|q\|_{H^{-\frac12}}\|q\|_{H^\sigma}^2,\notag
\end{align}
provided \(\sigma\geq \frac12\). The estimate \eqref{tr A3 high reg} is now obtained by combining \eqref{a4}, \eqref{a3 decomp}, \eqref{tr B3 high reg}, \eqref{tr B2C}, and \eqref{a3 lot}.

It remains to prove \eqref{Poisson commute}. We first note that if \(f\in L^2\) and \(j=0,1,2\) then we may apply \cite[Theorem 4.1]{MR2154153} to show that
\begin{align}
\bigl\|\tfrac{\p_x^j}{iL_0 +i\kappa^3 - i\kappa}f\bigr\|_{\I_2} & \lesssim \bigl\|\tfrac{\xi^j}{\xi^3 + \xi + i\kappa^3 - i\kappa}\bigr\|_{L_\xi^2} \|f\|_{L_x^2}\lesssim |\kappa|^{j-\frac52}\|f\|_{L^2}.\label{af}
\end{align}
For \((q,r)\in H^3\times H^2\) we compute that
\[
\frac{\delta\cH}{\delta r} = r\qtq{and}\frac{\delta\cH}{\delta q} = \tfrac43q - \tfrac13q_{xx} - \tfrac43 q^2,
\]
and hence the Poisson bracket
\begin{align}
\bigl\{\alpha,\cH\bigr\} &= \int_\R \frac{\delta\alpha}{\delta q} r_x + \frac{\delta\alpha}{\delta r} \bigl(\tfrac43q - \tfrac13q_{xx} - \tfrac43 q^2\bigr)_x\,dx\notag\\
&= \sum_{\ell=2}^\infty (-1)^\ell i\tr\bigl\{A^{\ell-1}\sqrt{R_0}(r_x\p_x+\p_xr_x + \tfrac43q_x - \tfrac13q_{xxx} - \tfrac83 qq_x)\sqrt{R_0}\bigr\}\notag\\
&= \sum_{\ell=2}^\infty (-i)^\ell\tr\left\{R_0(r_x\p_x + \p_xr_x + \tfrac43q_x - \tfrac13q_{xxx} - \tfrac83 qq_x)\left[R_0(L - L_0)\right]^{\ell-1}\right\}.\label{PB exp}
\end{align}
After possibly shrinking the size of \(\delta_0>0\), the series \eqref{PB exp} is readily seen to converge absolutely by applying \eqref{af} to bound the leading order term by
\begin{align*}
&\Bigl|\tr\left\{R_0(r_x\p_x + \p_xr_x + \tfrac43q_x - \tfrac13q_{xxx} - \tfrac83 qq_x)R_0(L - L_0)\right\}\Bigr|\\
&\qquad\lesssim\bigl\|R_0(r_x\p_x + \p_xr_x + \tfrac43q_x - \tfrac13q_{xxx} - \tfrac83 qq_x)\bigr\|_{\I_2}\bigl\|R_0(L - L_0)\bigr\|_{\I_2}\\
&\qquad\lesssim |\kappa|^{-3}\|(q,r)\|_{H^3\times H^2}^2\bigl[1 + \|(q,r)\|_{H^3\times H^2}\bigr] ,
\end{align*}
and applying \eqref{af}, \eqref{A I2}, and Sobolev embedding to bound the remaining terms by
\begin{align*}
&\Bigl|\tr\left\{R_0(r_x\p_x + \p_xr_x + \tfrac43q_x - \tfrac13q_{xxx} - \tfrac83 qq_x)\left[R_0(L - L_0)\right]^{\ell-1}\right\}\Bigr|\\
&\qquad\leq \bigl\|R_0(r_x\p_x + \p_xr_x + \tfrac43q_x - \tfrac13q_{xxx} - \tfrac83 qq_x)\bigr\|_{\I_2}\bigl\|\sqrt{R_0}\bigr\|_\op\|A\|_{\I_2}^{\ell-2}\bigl\|\sqrt{R_0}(L - L_0)\bigr\|_\op\\
&\qquad \lesssim C^\ell|\kappa|^{-\frac32-\ell}\delta_0^{\ell-2}\|(q,r)\|_{H^3\times H^2}^2\bigl[1 + \|(q,r)\|_{H^3\times H^2}\bigr],
\end{align*}
for some constant \(C>0\) and all \(\ell\geq 3\).

We now compute that
\begin{align*}
[\p_x^2,L - L_0] &= r_x\p_x + \p_xr_x + q_{xxx} + 4q_x\p_x^2 + 4q_{xx}\p_x,\\
[\tfrac43q,L_0 + i\kappa^3 - i\kappa] &= \tfrac43q_x - \tfrac43q_{xxx} - 4q_{xx}\p_x - 4q_x\p_x^2,\\
[\tfrac43q,L-L_0] &= -\tfrac83qq_x,
\end{align*}
and therefore we may decompose
\begin{align}
&(-i)^\ell\tr\left\{R_0(r_x\p_x + \p_xr_x + \tfrac43q_x - \tfrac13q_{xxx} - \tfrac83 qq_x)\left[R_0(L-L_0)\right]^{\ell-1}\right\}\label{tr ell}\\
&\qquad= (-i)^\ell\tr\left\{R_0[\p_x^2,L - L_0]\left[R_0(L - L_0)\right]^{\ell-1}\right\}\notag\\
&\qquad\quad + (-i)^\ell\tr\left\{R_0[\tfrac43q,L_0 + i\kappa^3 - i\kappa]\left[R_0(L - L_0)\right]^{\ell-1}\right\}\notag\\
&\qquad\quad + (-i)^\ell\tr\left\{R_0[\tfrac43q,L - L_0]\left[R_0(L - L_0)\right]^{\ell-1}\right\}.\notag
\end{align}

For the first summand on \RHS{tr ell}, we may write
\begin{align*}
&(-i)^\ell\tr\left\{R_0[\p_x^2,L - L_0]\left[R_0(L - L_0)\right]^{\ell-1}\right\}\\
&\qquad= (-i)^\ell\tr\left\{R_0\p_x^2(L - L_0)\left[R_0(L - L_0)\right]^{\ell-1}\right\}\\
&\qquad\quad - (-i)^\ell\tr\left\{R_0(L - L_0)R_0\p_x^2(L - L_0)\left[R_0(L - L_0)\right]^{\ell-2}\right\}\\
&\qquad= 0,
\end{align*}
where the second equality follows from cycling the trace.

For the second summand on \RHS{tr ell}, when \(\ell=2\) we explicitly compute that
\begin{align*}
&(-i)^\ell\tr\left\{R_0[\tfrac43q,L_0 + i\kappa^3 - i\kappa]R_0(L - L_0)\right\}\\
&\qquad= -\tfrac2{3\pi}(-i)^\ell\iint_{\R^2} \tfrac{(\xi-\eta)(\xi+\eta)\left[1 + (\xi-\eta)^2 + 3\xi\eta\right]|\hat q(\xi-\eta)|^2}{(\xi^3 + \xi + i\kappa^3 - i\kappa)(\eta^3 + \eta + i\kappa^3 - i\kappa)}\,d\xi\,d\eta\\
&\qquad\quad + \tfrac{2i}{3\pi}(-i)^\ell\iint_{\R^2} \tfrac{(\xi-\eta)\left[1 + (\xi-\eta)^2 + 3\xi\eta\right] \hat q(\xi-\eta)\overline{\hat r(\xi - \eta)}}{(\xi^3 + \xi + i\kappa^3 - i\kappa)(\eta^3 + \eta + i\kappa^3 - i\kappa)}\,d\xi\,d\eta.
\end{align*}
Symmetrizing under \(\xi\leftrightarrow\eta\), we immediately see that
\[
\iint_{\R^2} \tfrac{(\xi-\eta)(\xi+\eta)\left[1 + (\xi-\eta)^2 + 3\xi\eta\right]|\hat q(\xi-\eta)|^2}{(\xi^3 + \xi + i\kappa^3 - i\kappa)(\eta^3 + \eta + i\kappa^3 - i\kappa)}\,d\xi\,d\eta = 0,
\]
whereas by computing the contour integral as in the proof of \eqref{tr A2}, we obtain
\[
\iint_{\R^2} \tfrac{(\xi-\eta)\left[1 + (\xi-\eta)^2 + 3\xi\eta\right] \hat q(\xi-\eta)\overline{\hat r(\xi - \eta)}}{(\xi^3 + \xi + i\kappa^3 - i\kappa)(\eta^3 + \eta + i\kappa^3 - i\kappa)}\,d\xi\,d\eta = 0.
\]

Finally, when \(\ell\geq 3\) the second summand on \RHS{tr ell} satisfies
\begin{align*}
&(-i)^\ell\tr\left\{R_0[\tfrac43q,L_0 + i\kappa^3 - i\kappa]\left[R_0(L - L_0)\right]^{\ell-1}\right\}\\
&\qquad = - (-i)^{\ell-1}\tr\left\{R_0\tfrac43q(L-L_0)\left[R_0(L - L_0)\right]^{\ell-2}\right\}\\
&\qquad \quad + (-i)^{\ell-1} \tr\left\{\tfrac43q\left[R_0(L - L_0)\right]^{\ell-1}\right\}\\
&\qquad = - (-i)^{\ell-1} \tr\left\{R_0[\tfrac43q,(L-L_0)]\left[R_0(L - L_0)\right]^{\ell-2}\right\},
\end{align*}
which cancels the third summand on \RHS{tr ell} with \(\ell\) replaced by \(\ell-1\). Putting these together we see that the series \eqref{PB exp} telescopes, thereby completing the proof of \eqref{Poisson commute}.
\epf

We now turn to the

\bpf[Proof of Theorem~\ref{t:AP}]
From \eqref{tr A2}, for \(\vk\geq 2\) we have
\[
\alpha\sprack 2(\vk) + \alpha\sprack 2(-\vk)\simeq \tfrac1\vk\Bigl[\|\sqrt{1-\p_x^2}q\|_{H_\vk^{-2}}^2 + \|r\|_{H_\vk^{-2}}^2\Bigr].
\]
It follows that for all \(-\frac12\leq s<1\) and \(\kappa\geq 2\) we have
\eq{alpha 2 main}{
\int_{\kappa}^\infty \vk^{2s+2}\bigl[\alpha\sprack 2(\vk) + \alpha\sprack 2(-\vk)\bigr]\,d\vk\simeq_s \bigl\|\sqrt{1-\p_x^2}\,q\|_{H_\kappa^{s-1}}^2 + \|r\|_{H_\kappa^{s-1}}^2.
}

Further, taking \(\delta_0>0\) as in Proposition~\ref{p:alpha}, if \((q,r)\in B_{\delta_0}\) then we may apply \eqref{tr A3} to bound
\begin{align}
&\int_\kappa^\infty \vk^{2s+2}\bigl|\alpha\sprack{\geq 3}(\vk) + \alpha\sprack{\geq 3}(-\vk)\bigr|\,d\vk\label{alpha geq3 main}\\
&\qquad\lesssim_{s,\sigma} \kappa^{\max\{2s-2\sigma-1,2s-2\}}\delta_0\|(q,r)\|_{H^\sigma\times H^{\sigma-1}}^2,\notag
\end{align}
whenever \(-\frac12\leq \sigma\leq s<\min\{1,\sigma+\frac12\}\) and \(\kappa\geq 2\).

Taking \(s=-\frac12=\sigma\) and \(\kappa=2\), we may combine \eqref{alpha 2 main} and \eqref{alpha geq3 main} to show that, after possibly shrinking the size of \(\delta_0>0\), we have
\eq{alpha -1/2}{
\int_2^\infty \vk\bigl[\alpha(\vk) + \alpha(-\vk)\bigr]\,d\vk\simeq \|(q,r)\|_{H^{-\frac12}\times H^{-\frac32}}^2,
}
for all \((q,r)\in B_{\delta_0}\).

Now consider an interval \(I\subseteq \R\) containing zero and a solution \((q,r)\in C(I;H^3\times H^2)\) of \eqref{GB-system} satisfying \eqref{init small} with \(0<\delta\leq \delta_0\). By continuity, there exists a maximal open interval \(J\subseteq I\) containing zero so that \((q,r)(t)\in B_{\delta_0}\) for all \(t\in J\). From \eqref{Poisson commute}, we see that for all \(t\in J\) and \(\vk\geq 2\) we have
\[
\tfrac d{dt}\alpha(\vk) = \bigl\{\alpha(\vk),\cH\bigr\} = 0,
\]
and hence \LHS{alpha -1/2} is constant. It follows that for all \(t\in J\) we have
\[
\|(q,r)(t)\|_{H^{-\frac12}\times H^{-\frac32}}\simeq \|(q,r)(0)\|_{H^{-\frac12}\times H^{-\frac32}},
\]
so provided \(0<\delta\ll\delta_0\) is chosen sufficiently small we have \(J = I\) and obtain \eqref{AP bound} for \(s=-\frac12\).

To complete the proof of \eqref{AP bound}, we will show that for all \(-\frac12<s<2\) we have
\eq{AP bound half}{
\|(q,r)(t)\|_{H^s\times H^{s-1}}\lesssim_s \|(q,r)(t')\|_{H^s\times H^{s-1}},
}
uniformly for all \(t,t'\in I\). As \(t,t'\in I\) are arbitrary, this immediately implies \eqref{AP bound}.

For \(-\frac12<s<1\), we first use \eqref{alpha 2 main} with \(\kappa = 2\) to estimate
\[
\int_2^\infty \vk^{2s+2}\bigl[\alpha\sprack 2(\vk) + \alpha\sprack 2(-\vk)\bigr]\,d\vk\simeq_s\|(q,r)\|_{H^s\times H^{s-1}}^2.
\]
We will control the contribution of \(\alpha\sprack{\geq 3}\) using a lower regularity norm than \(H^s\times H^{s-1}\). This will ensure that our choice of \(\delta>0\) is independent of the value of \(s\). To this end, if \(-\frac12<s<0\) we may apply \eqref{alpha geq3 main} with \(\sigma=-\frac12\) to bound
\[
\int_2^\infty \vk^{2s+2}\bigl|\alpha\sprack{\geq 3}(\vk) + \alpha\sprack{\geq 3}(-\vk)\bigr|\,d\vk\lesssim_s\delta_0\|(q,r)\|_{H^{-\frac12}\times H^{-\frac32}}^2.
\]
As we have already proved that \eqref{AP bound} holds in the case \(s=-\frac12\), this yields \eqref{AP bound half} for \(-\frac12<s<0\).

For \(0\leq s<1\) we instead take \(\sigma = \frac14\lfloor 4s-1\rfloor\) and estimate
\[
\int_2^\infty \vk^{2s+2}\bigl|\alpha\sprack{\geq 3}(\vk) + \alpha\sprack{\geq 3}(-\vk)\bigr|\,d\vk\lesssim_s\delta_0\|(q,r)\|_{H^\sigma\times H^{\sigma-1}}^2.
\]
By inductively considering \(s\) in each of the intervals \([0,\frac14)\), \([\frac14,\frac12)\), \([\frac12,\frac34)\), and \([\frac34,1)\), we obtain \eqref{AP bound half} for the range \(0\leq s<1\).

For the case \(s=1\), we use that \(H^{\frac34}\) is an algebra to bound
\begin{align*}
\left|\int_\R q^3\,dx\right| \lesssim \|q\|_{H^{-\frac12}}\|q^2\|_{H^{\frac12}}\lesssim \delta_0\|q\|_{H^{\frac34}}^2.
\end{align*}
We may then combine conservation of the Hamiltonian \(\cH\) with the \(s=\frac34\) case of \eqref{AP bound half} to complete the proof of \eqref{AP bound half} in this case.

Finally, we turn to the case \(1<s<2\). We start by decomposing the Hamiltonian \(\cH = \cH\sprack 2 + \cH\sprack 3\), where
\[
\cH\sprack 2 = \int_\R \tfrac23 q^2 + \tfrac16 q_x^2 +  \tfrac12r^2\,dx \qtq{and}\cH\sprack 3 = - \int_\R\tfrac49q^3\,dx.
\]
Using \eqref{tr A2}, we compute that for \(\vk\geq 2\) we have
\begin{align}\label{corrected alpha ID}
&\tfrac{12\vk}{(3\vk^2-1)^3}\cH\sprack 2 - \alpha\sprack 2(\vk) - \alpha\sprack 2(-\vk)\\
&\qquad= \tfrac{2\vk}{(3\vk^2-1)^3}\int_\R \tfrac{(\xi^2 + 3\vk^2+2)\xi^2\left[(4+\xi^2)|\hat q(\xi)|^2 + 3|\hat r(\xi)|^2\right]}{\xi^4 + (3\vk^2+2)\xi^2 + (3\vk^2-1)^2}\,d\xi .\notag
\end{align}
It follows that for \(1<s<2\) and \(\kappa\geq 2\) we have the estimate
\begin{align}
&\int_\kappa^\infty \vk^{2s+2}\Bigl[\tfrac{12\vk}{(3\vk^2-1)^3}\cH\sprack 2 - \alpha\sprack 2(\vk) - \alpha\sprack 2(-\vk)\Bigr]\,d\vk\label{alpha 2 high}\\
&\qquad\simeq_s \|(\sqrt{1-\p_x^2}\,q_x,r_x)\|_{H_\kappa^{s-2}\times H_\kappa^{s-2}}^2.\notag
\end{align}
From \eqref{tr A3 high reg}, we also have
\begin{align}
&\int_\kappa^\infty \vk^{2s+2}\Bigl|\tfrac{12\vk}{(3\vk^2-1)^3}\cH\sprack 3 - \alpha\sprack{\geq 3}(\vk) - \alpha\sprack{\geq 3}(-\vk)\Bigr|\,d\vk\label{alpha geq3 high}\\
&\qquad\lesssim_{s,\sigma} \kappa^{\max\{2s-2\sigma-1,2s-4\}}\delta_0\|(q,r)\|_{H^\sigma\times H^{\sigma-1}}^2,\notag
\end{align}
whenever \(\frac12\leq \sigma\leq s<\min\{2,\sigma+\frac12\}\) and \(\kappa\geq 2\).

By taking \(\kappa = 2\), \(\sigma = \frac14\lfloor 4s-1\rfloor\), and using that
\[
\tfrac d{dt}\int_\kappa^\infty \vk^{2s+2}\Bigl[\tfrac{12\vk}{(3\vk^2-1)^3}\cH - \alpha(\vk) - \alpha(-\vk)\Bigr]\,d\vk = 0,
\]
we obtain the estimate
\[
\|(q_x,r_x)(t)\|_{H^{s-1}\times H^{s-2}}^2\lesssim_s \|(q,r)(t')\|_{H^s\times H^{s-1}}^2 + \delta_0\|(q,r)(t)\|_{H^\sigma\times H^{\sigma-1}}^2,
\]
uniformly for all \(t,t'\in I\). To control the low frequencies, we combine this with \eqref{AP bound half} in the case \(s=1\) to get
\eq{first correction}{
\|(q,r)(t)\|_{H^s\times H^{s-1}}^2\lesssim_s \|(q,r)(t')\|_{H^s\times H^{s-1}}^2 + \delta_0\|(q,r)(t)\|_{H^\sigma\times H^{\sigma-1}}^2,
}
for all \(t,t'\in I\).

To complete the proof of \eqref{AP bound half} in the range \(1\leq s<2\), we inductively apply \eqref{first correction} for \(s\) in each of the intervals \((1,\frac54)\), \([\frac54,\frac32)\), \([\frac32,\frac74)\), and \([\frac74,2)\).
\epf

With the a priori estimates of Theorem~\ref{t:AP} in hand, let us now give the
\bpf[Proof of Theorem~\ref{t:GWP}]
By combining Theorem~\ref{t:AP} with \cite[Theorem 4]{MR954673}, we see that the solution map \(\Phi\colon\R\times (H^3\times H^2)\cap B_\delta\to H^3\times H^2\) is well-defined and for all \((q,r)\in (H^3\times H^2)\cap B_\delta\) we have the estimate
\eq{soln unif}{
\sup_{t\in \R}\|\Phi_t(q,r)\|_{H^{-\frac12}\times H^{-\frac32}}\lesssim \delta.
}

As the time of existence in Theorem~\ref{t:LWP} depends only on the size of the initial data in \(H^{-\frac12}\times H^{-\frac32}\), we may then iteratively combine the local well-posedness of Theorem~\ref{t:LWP} with our uniform bound \eqref{soln unif} to show that the solution map \(\Phi\) can be either restricted (if \(s>3\)) or uniquely extended (if \(-\frac12\leq s<3\)) to a continuous map from \(\R\times \bigl((H^s\times H^{s-1})\cap B_\delta\bigr)\to H^s\times H^{s-1}\). Finally, we note that local Lipschitz continuity of the solution map, uniformly on bounded time intervals, follows directly from local well-posedness.
\epf

By combining the continuity of the solution map established in Theorem~\ref{t:GWP} with the a priori estimates from Theorem~\ref{t:AP}, we see that if \(-\frac12 \leq s<2\) and \((q,r)(0)\in (H^s\times H^{s-1})\cap B_\delta\) then the corresponding global solution \((q,r)\in C(\R;H^s\times H^{s-1})\) satisfies the a priori estimate \eqref{AP bound} for all \(t\in \R\). Using this fact, we conclude this section with the

\bpf[Proof of Theorem~\ref{t:equicontinuity}]
We recall that a bounded set \(\cQ\subseteq H^s\times H^{s-1}\) is equicontinuous if and only if
\[
\lim_{N\to \infty}\sup_{(q,r)\in \cQ}\|P_{>N}(q,r)\|_{H^s\times H^{s-1}} = 0.
\]
Proceeding as in the proof of \cite[Lemma 2.4]{MR4743514}, we first observe that for all \(s\in \R\) and \(\kappa\geq 2\) we have
\begin{align*}
\|P_{>\kappa}(q,r)\|_{H^s\times H^{s-1}}&\lesssim \|(\sqrt{1-\p_x^2}\,q_x,r_x)\|_{H_\kappa^{s-2}\times H_\kappa^{s-2}}\lesssim \|(\sqrt{1-\p_x^2}q,r)\|_{H_\kappa^{s-1}\times H_\kappa^{s-1}}.
\end{align*}
Conversely, if \(s<1\) and \(\kappa\geq 2\) then
\[
\|(\sqrt{1-\p_x^2}\,q,r)\|_{H_\kappa^{s-1}\times H_\kappa^{s-1}} \lesssim \bigl(\tfrac N\kappa\bigr)^{1-s}\|(q,r)\|_{H^s\times H^{s-1}} + \|P_{>N}(q,r)\|_{H^s\times H^{s-1}},
\]
whereas if \(s<2\) then
\[
\|(\sqrt{1-\p_x^2}\,q_x,r_x)\|_{H_\kappa^{s-2}\times H_\kappa^{s-2}} \lesssim \bigl(\tfrac N\kappa\bigr)^{2-s}\|(q,r)\|_{H^s\times H^{s-1}} + \|P_{>N}(q,r)\|_{H^s\times H^{s-1}}.
\]
It follows that if \(s<1\) then a bounded set \(\cQ\subseteq H^s\times H^{s-1}\) is equicontinuous if and only if
\begin{align}\label{equi char}
\lim_{\kappa\to\infty}\sup_{(q,r)\in \cQ}\|(\sqrt{1-\p_x^2}\,q,r)\|_{H_\kappa^{s-1}\times H_\kappa^{s-1}} &= 0,
\end{align}
and if \(s<2\) then a bounded set 
\(\cQ\subseteq H^s\times H^{s-1}\) is equicontinuous if and only if
\begin{align}\label{equi char hi}
\lim_{\kappa\to\infty}\sup_{(q,r)\in \cQ}\|(\sqrt{1-\p_x^2}\,q_x,r_x)\|_{H_\kappa^{s-2}\times H_\kappa^{s-2}} &= 0.
\end{align}

We now proceed using several estimates established in the proof of Theorem~\ref{t:AP}. First, we consider the case where \(-\frac12\leq s<1\). By combining \eqref{alpha 2 main} and \eqref{alpha geq3 main} (for \(\sigma=s\)) with the a priori estimate \eqref{AP bound} and conservation of \(\alpha\), we obtain the estimate
\begin{align*}
\sup_{t\in \R}\bigl\|\bigl(\sqrt{1-\p_x^2}\,q,r\bigr)(t)\bigr\|_{H_\kappa^{s-1}\times H_\kappa^{s-1}}^2&\lesssim \bigl\|\bigl(\sqrt{1-\p_x^2}\,q,r\bigr)(0)\bigr\|_{H_\kappa^{s-1}\times H_\kappa^{s-1}}^2\\
&\quad + \delta_0\kappa^{\max\{-1,2s-2\}}\|(q,r)(0)\|_{H^s\times H^{s-1}}^2,
\end{align*}
whenever \(\kappa\geq 2\). Consequently, if \(\cQ\subseteq (H^s\times H^{s-1})\cap B_\delta\) is equicontinuous then we may apply \eqref{equi char} to obtain
\[
\lim_{\kappa\to\infty}\sup_{(q,r)(0)\in \cQ}\sup_{t\in \R}\bigl\|\bigl(\sqrt{1-\p_x^2}\,q,r\bigr)(t)\bigr\|_{H_\kappa^{s-1}\times H_\kappa^{s-1}} = 0,
\]
and hence \(\cQ^*\) is equicontinuous.

For \(1<s<2\), we instead use \eqref{alpha 2 high} and \eqref{alpha geq3 high} (again with \(\sigma=s\)) followed by \eqref{AP bound} to obtain
\begin{align*}
\sup_{t\in \R}\bigl\|\bigl(\sqrt{1-\p_x^2}\,q_x,r_x\bigr)(t)\bigr\|_{H_\kappa^{s-2}\times H_\kappa^{s-2}}^2&\lesssim \bigl\|\bigl(\sqrt{1-\p_x^2}\,q_x,r_x\bigr)(0)\bigr\|_{H_\kappa^{s-2}\times H_\kappa^{s-2}}^2\\
&\quad + \delta_0\kappa^{\max\{-1,2s-4\}}\|(q,r)(0)\|_{H^s\times H^{s-1}}^2,
\end{align*}
and an application of \eqref{equi char hi} now proves that \(\cQ^*\) is equicontinuous.

Finally, for the case \(s=1\) we observe that \eqref{corrected alpha ID} gives us the estimate
\[
\tfrac{12\kappa}{(3\kappa^2-1)^3}\cH\sprack 2 - \alpha\sprack 2(\kappa) - \alpha\sprack 2(-\kappa)\simeq \kappa^{-5}\bigl\|\bigl(\sqrt{1-\p_x^2}\,q_x,r_x\bigr)(t)\bigr\|_{H_\kappa^{-1}\times H_\kappa^{-1}}^2,
\]
whereas \eqref{tr A3 high reg} gives us
\begin{align*}
\Bigl|\tfrac{12\kappa}{(3\kappa^2-1)^3}\cH\sprack{3} - \alpha\sprack{\geq 3}(\kappa) - \alpha\sprack{\geq 3}(-\kappa)\Bigr|\lesssim \kappa^{-6}\delta_0\|(q,r)\|_{H^1\times L^2}^2.
\end{align*}
Combining these bounds with the conservation of \(\alpha\), \(\cH\) and \eqref{AP bound} in the case \(s=1\) we obtain
\begin{align*}
\sup_{t\in \R}\bigl\|\bigl(\sqrt{1-\p_x^2}\,q_x,r_x\bigr)(t)\bigr\|_{H_\kappa^{-1}\times H_\kappa^{-1}}^2&\lesssim \bigl\|\bigl(\sqrt{1-\p_x^2}\,q_x,r_x\bigr)(0)\bigr\|_{H_\kappa^{-1}\times H_\kappa^{-1}}^2\\
&\quad + \delta_0\kappa^{-1}\|(q,r)(0)\|_{H^1\times L^2}^2,
\end{align*}
and applying \eqref{equi char hi} once again completes the proof.
\epf

\section{The long wave limit}\label{s:long wave}

In this section we prove Theorem~\ref{t:long wave}.  We start by compiling a couple of technical estimates into the following lemma:
\begin{lem}
Let \(s\in \R\). Then we have the estimates
\begin{align}
\|fg\|_{H^s}&\lesssim_s \|f\|_{W^{\lceil|s|\rceil,\infty}}\|g\|_{H^s},\label{product-1}\\
\|[\<\p_x\>^s,f]g\|_{L^2} &\lesssim_s \|f\|_{H^{|s|+2}}\|g\|_{H^{s-1}}.\label{commutator}
\end{align}
\end{lem}
\bpf
The estimate \eqref{product-1} should be viewed as a version of the fractional Leibniz rule (see, e.g., \cite{MR3200091}) and the estimate \eqref{commutator} as a version of the Kato--Ponce commutator estimate that originally appeared in~\cite{MR951744}.  For completeness, we give a self-contained proof of the precise (not especially sharp) versions of these estimates that we require.

For \eqref{product-1} we may argue by duality to restrict attention to the case that \(s\geq 0\).

When \(s\geq 0\) is an integer, \eqref{product-1} follows from the fact that
\[
\|fg\|_{H^s}^2\simeq_s \|fg\|_{L^2}^2 + \|\p_x^s(fg)\|_{L^2}^2,
\]
and the Leibniz rule.

For non-integer \(s>0\), we introduce the paraproducts
\[
T_fg := \sum_{N\geq 8}f_{\leq \frac N8}g_N\qtq{and}
\Pi(f,g) := 
fg - T_fg - T_gf
\]
For the low-high interactions we estimate
\begin{align*}
\|T_fg\|_{H^s}^2 \lesssim_s \sum_{N\geq 8} N^{2s}\|f_{\leq \frac N8}\|_{L^\infty}^2\|g_N\|_{L^2}^2\lesssim_s \|f\|_{L^\infty}^2\|g\|_{H^s}^2,
\end{align*}
whereas for the high-low interactions, we use the Cauchy--Schwarz inequality and that \(\lceil s\rceil>s\) to bound
\begin{align*}
\|T_gf\|_{H^s}^2 &\lesssim_s \sum_{N\geq 8} N^{2s}\|f_N\|_{L^\infty}^2\|g_{\leq \frac N8}\|_{L^2}^2\\
&\lesssim_s \|f\|_{W^{\lceil s\rceil,\infty}}^2\|g\|_{L^2}^2\sum_{N\geq 8}N^{2(s - \lceil s\rceil)}\lesssim_s \|f\|_{W^{\lceil s\rceil,\infty}}^2\|g\|_{L^2}^2.
\end{align*}
For the high-high interactions we again use that \(s<\lceil s\rceil\) to bound
\begin{align*}
\|\Pi(f,g)\|_{H^s}&\lesssim_s \sum_{N_1\simeq N_2}\sum_{M\lesssim N_1} M^s\|f_{N_1}\|_{L^\infty}\|g_{N_2}\|_{L^2}\\
&\lesssim_s \|f\|_{W^{\lceil s\rceil,\infty}}\sum_{N\geq 1} N^{s-\lceil s\rceil}\|g_N\|_{L^2}\lesssim_s \|f\|_{W^{\lceil s\rceil,\infty}}\|g\|_{L^2}.
\end{align*}
Combining these estimates completes the proof of \eqref{product-1}.

For the commutator estimate \eqref{commutator} we decompose
\begin{align*}
[\<\p_x\>^s,f]g = [\<\p_x\>^s,T_f]g + \<\p_x\>^sT_g f - T_{\<\p_x\>^sg}f + \<\p_x\>^s\Pi[f,g] - \Pi[f,\<\p_x\>^s g].
\end{align*}

For the first term, we apply Plancherel's Theorem followed by Young's convolution inequality and the Cauchy--Schwarz inequality to estimate
\begin{align*}
\|[\<\p_x\>^s,f_{\leq\frac N8}]g_N\|_{L^2} &\simeq\left\|\int_\R \bigl[\<\xi\>^s - \<\eta\>^s\bigr]\hat f_{\leq\frac N8}(\xi-\eta)\hat g_N(\eta)\,d\eta\right\|_{L_\xi^2}\\
&\lesssim_s N^{s-1}\|\xi\hat f\|_{L_\xi^1}\|\hat g_N\|_{L_\xi^2}\lesssim_s N^{s-1} \|f\|_{H^2}\|g_N\|_{L^2},
\end{align*}
where we have used the Mean Value Theorem to deduce that
\[
\bigl|\<\xi\>^s - \<\eta\>^s\bigr|\lesssim_s N^{s-1}|\xi-\eta|
\]
on the support of the \(\hat f_{\leq\frac N8}(\xi-\eta)\hat g_N(\eta)\). It follows that
\begin{align*}
\|[\<\p_x\>^s,T_f]g\|_{L^2}^2 &\lesssim \sum_{N\geq 8}\|[\<\p_x\>^s,f_{\leq\frac N8}]g_N\|_{L^2}^2\lesssim_s \|f\|_{H^{|s|+2}}^2\|g\|_{H^{s-1}}^2.
\end{align*}

For the remaining terms, we first apply Bernstein's inequality to estimate
\begin{align*}
\|\<\p_x\>^s T_gf\|_{L^2}^2 &\lesssim_s \sum_{N\geq 8}N^{2s}\|f_N\|_{L^2}^2 \|g_{\leq\frac N8}\|_{L^\infty}^2\\
&\lesssim_s \sum_{N\geq 8}N^{\max\{3,2s+1\}}\|f_N\|_{L^2}^2\|g\|_{H^{s-1}}^2\lesssim_s \|f\|_{H^{|s|+2}}^2\|g\|_{H^{s-1}}^2, \\
\|T_{\<\p_x\>^sg}f\|_{L^2}^2&\lesssim_s \sum_{N\geq 8} N^3\|f_N\|_{L^2}^2\|g_{\leq \frac N8}\|_{H^{s-1}}^2\lesssim_s \|f\|_{H^{|s|+2}}^2\|g\|_{H^{s-1}}^2,\\
\|\Pi[f,\<\p_x\>^sg]\|_{L^2} &\lesssim_s \sum_{N_1\simeq N_2}\sum_{M\lesssim N_1}M^{\frac12}N_1\|f_{N_1}\|_{L^2}\|g_{N_2}\|_{H^{s-1}}\lesssim_s \|f\|_{H^{|s|+2}}\|g\|_{H^{s-1}}.
\end{align*}

For the final term, when \(s>-\frac12\) we again apply Bernstein's inequality and then the Cauchy--Schwarz inequality to bound
\begin{align*}
\|\<\p_x\>^s\Pi[f,g]\|_{L^2}&\lesssim \sum_{N_1\simeq N_2}\sum_{M\lesssim N_1}M^{s+\frac12}\|f_{N_1}\|_{L^2}\|g_{N_2}\|_{L^2}\\
&\lesssim \sum_{N_1\simeq N_2}N_1^{\frac32}\|f_{N_1}\|_{L^2}\|g_{N_2}\|_{H^{s-1}}\lesssim \|f\|_{H^{|s|+2}}\|g\|_{H^{s-1}}.
\end{align*}
The case \(s<-\frac12\) is treated similarly to obtain
\begin{align*}
\|\<\p_x\>^s\Pi[f,g]\|_{L^2}
&\lesssim \sum_{N_1\simeq N_2}N_1^{|s|+1}\|f_{N_1}\|_{L^2}\|g_{N_2}\|_{H^{s-1}}\lesssim \|f\|_{H^{|s|+2}}\|g\|_{H^{s-1}}.
\end{align*}
Finally, in the case \(s=-\frac12\) we have
\begin{align*}
\|\<\p_x\>^{-\frac12}\Pi[f,g]\|_{L^2}
&\lesssim \sum_{N_1\simeq N_2}\log(1 + N_1)N_1^{\frac32}\|f_{N_1}\|_{L^2}\|g_{N_2}\|_{H^{-\frac32}}\lesssim \|f\|_{H^{\frac52}}\|g\|_{H^{s-1}},
\end{align*}
which completes the proof.
\epf

The first point of order in our proof of Theorem~\ref{t:long wave} will be proving that sufficiently smooth solutions \(w^{\epsilon,\pm}\) of \eqref{KdV} are approximate solutions  of \eqref{GB-system-2}. To prove this, we will require the following:

\begin{lem}\label{l:ap stat}
Let \(\sigma\geq -1\) and \(0<\epsilon_0\leq 1\).  Suppose that for each \(0<\epsilon\leq \epsilon_0\) we have solutions \(w^{\epsilon,\pm}\in C(\R;H^\sigma)\) of \eqref{KdV} so that for some \(w_0^\pm \in H^\sigma\) we have
\eq{KdV init cvgce}{
\lim_{\epsilon \to 0}\|(w^{\epsilon,+},w^{\epsilon,-})(0) - (w_0^+,w_0^-)\|_{H^\sigma\times H^\sigma} = 0.
}
Then we have the estimate
\eq{uniform consequence KdV}{
\limsup_{\epsilon\to0}\sup_{t\in \R}\|(w^{\epsilon,+},w^{\epsilon,-})(t)\|_{H^\sigma\times H^\sigma}\lesssim_{\sigma,\|w_0^\pm\|_{H^\sigma}}1,
}
and if we define
\eq{chi geq R eps}{
\chi_{>R}^{\epsilon,\pm}(t,x) := \chi_{>R}\bigl(x\mp\tfrac2{\sqrt 3\epsilon^2}t\bigr)\qt{for \(R\geq 1\)},
}
then for all \(T>0\) we have the identity
\begin{align}\label{tight}
\lim_{R\to \infty}\limsup_{\epsilon\to0}\sup_{|t|\leq T}\|(\chi_{>R}^{\epsilon,+}\,w^{\epsilon,+},\chi_{>R}^{\epsilon,-}\,w^{\epsilon,-})(t)\|_{H^\sigma\times H^\sigma} = 0.
\end{align}
\end{lem}
\bpf
From \eqref{KdV init cvgce}, we may find some \(0<\epsilon_1\leq \epsilon_0\) so that
\[
\sup_{0<\epsilon\leq \epsilon_1}\|(w^{\epsilon,+},w^{\epsilon,-})(0)\|_{H^\sigma\times H^\sigma}\leq 1 + \|(w_0^+,w_0^-)\|_{H^\sigma\times H^\sigma}.
\]
We may then combine the a priori estimates of \cite{MR3874652,MR3820439} with the well-posedness of \eqref{KdV stat} from \cite{MR3990604,MR1086966,MR2018661,MR1969209,MR2531556,MR2501679} and the change of variables \eqref{static rescaling} to obtain
\[
\sup_{0<\epsilon\leq \epsilon_1}\sup_{t\in \R}\|(w^{\epsilon,+},w^{\epsilon,-})(t)\|_{H^\sigma\times H^\sigma}\lesssim_{\sigma,\|w_0^\pm\|_{H^\sigma}} 1,
\]
which yields \eqref{uniform consequence KdV}.

To prove \eqref{tight}, we take a sequence \(\epsilon_0\geq \epsilon_n\to 0\) so that, by \eqref{KdV init cvgce},
\[
\mc W := \bigl\{(w^{\epsilon_n,+},w^{\epsilon_n,-})(0):n\geq 1\bigr\}
\]
is a precompact subset of \(H^\sigma\times H^\sigma\). From the well-posedness of \eqref{KdV stat}, we see that if \(w_n^\pm\in C(\R;H^\sigma)\) solve \eqref{KdV stat} with initial data \(w_n^\pm (0) = w^{\epsilon_n,\pm}(0)\) then for all \(T>0\) the set
\[
\bigl\{(w_n^+,w_n^-)(t):|t|\leq T\text{ and }n\geq 1\bigr\}
\]
is the continuous image of \([-T,T]\times \mc W\) and therefore also precompact in \(H^\sigma\times H^\sigma\).  Using that precompact subsets of \(H^\sigma\times H^\sigma\) are uniformly tight (see, e.g., the proof of \cite[Lemma 2.2]{harropgriffiths2026homogenizationnonlinearschrodingerequation}) we have
\[
\lim_{R\to\infty}\sup_{n\geq 1}\sup_{|t|\leq T}\|(\chi_{>R}\,w_n^+,\chi_{>R}\,w_n^-)(t)\|_{H^\sigma\times H^\sigma} = 0.
\]
The identity \eqref{tight} now follows from the change of variables \eqref{static rescaling}.
\epf

As a corollary, we have the following:
\begin{lem}
Let \(s\geq -\frac12\), \(m \geq \lceil s+1\rceil\) be an integer, and \(0<\epsilon_0\leq 1\). Suppose that for each \(0<\epsilon\leq \epsilon_0\) we have solutions \(w^{\epsilon,\pm}\in C(\R;H^m)\) of \eqref{KdV} satisfying \eqref{KdV init cvgce} with \(\sigma = m\) for some \(w_0^\pm\in H^m\). Then, for all \(T>0\) we have the estimate
\eq{cancellation}{
\lim_{\epsilon\to 0}\int_{-T}^T\|\Lambda_\epsilon^{-1}\p_x(w^{\epsilon,+}w^{\epsilon,-})\|_{H^s}\,dt = 0.
}
\end{lem}
\bpf
From \eqref{uniform consequence KdV}, we may find some \(0<\epsilon_1\leq \epsilon_0\) so that
\begin{align*}
A:= \sup_{0<\epsilon\leq \epsilon_1}\sup_{t\in \R}\|(w^{\epsilon,+},w^{\epsilon,-})(t)\|_{H^m\times H^m}<\infty.
\end{align*}

As \(m\geq 1\), the space \(H^m\) is an algebra and hence for any \(0<\nu<T\) we may estimate
\begin{align}
\int_{-\nu}^\nu\|\Lambda_\epsilon^{-1}\p_x(w^{\epsilon,+}w^{\epsilon,-})\|_{H^s}\,dt\lesssim \int_{-\nu}^\nu\|w^{\epsilon,+}\,w^{\epsilon,-}\|_{H^m}\,dt\lesssim \nu A^2.\label{short time}
\end{align}

Taking \(R\geq 1\), \(\chi_{>R}^{\epsilon,\pm}\) to be defined as in \eqref{chi geq R eps}, and setting \(\chi_{\leq R}^{\epsilon,\pm} := 1 - \chi_{>R}^{\epsilon,\pm}\), we may use that \(H^m\) is an algebra with the estimate \eqref{product-1}, to bound
\begin{align*}
\|\Lambda_\epsilon^{-1}\p_x(w^{\epsilon,+}w^{\epsilon,-})\|_{H^s}&\lesssim \|\chi_{>R}^{\epsilon,+}w^{\epsilon,+}\|_{H^m}\|\chi_{>R}^{\epsilon,-}\,w^{\epsilon,-}\|_{H^m}+ A^2\|\chi_{\leq R}^{\epsilon,+}\chi_{\leq R}^{\epsilon,-}\|_{W^{m,\infty}}\\
&\quad + A\Bigl[\|\chi_{>R}^{\epsilon,+}\,w^{\epsilon,+}\|_{H^m} + \|\chi_{>R}^{\epsilon,-}\,w^{\epsilon,-}\|_{H^m}\Bigr] .
\end{align*}
Integrating, we then have
\begin{align*}
&\int_{\nu\leq |t|\leq T}\|\Lambda_\epsilon^{-1}\p_x(w^{\epsilon,+}w^{\epsilon,-})\|_{H^s}\,dt\\
&\qquad\lesssim T\|\chi_{>R}^{\epsilon,+}\,w^{\epsilon,+}\|_{C_TH^m}\|\chi_{>R}^{\epsilon,-}\,w^{\epsilon,-}\|_{C_TH^m} + A^2\int_{\nu\leq |t|\leq T}\|\chi_{\leq R}^{\epsilon,+}\,\chi_{\leq R}^{\epsilon,-}\|_{W^{m,\infty}}\,dt\\
&\qquad\quad + AT\Bigl[\|\chi_{>R}^{\epsilon,+}\,w^{\epsilon,+}\|_{C_TH^m} + \|\chi_{>R}^{\epsilon,-}\,w^{\epsilon,-}\|_{C_TH^m}\Bigr].
\end{align*}

For fixed \(R\geq 1\), we have \(\chi_{\leq R}^{\epsilon,+}\,\chi_{\leq R}^{\epsilon,-} = 0\) whenever \(|t|\geq \nu\) and \(\epsilon \ll \frac\nu R\). Consequently, we may first pass to the limit as \(\epsilon \to 0\) to obtain
\begin{align*}
&\limsup_{\epsilon\to0}\int_{\nu\leq |t|\leq T}\|\Lambda_\epsilon^{-1}\p_x(w^{\epsilon,+}w^{\epsilon,-})\|_{H^s}\,dt\\
&\qquad\lesssim AT\limsup_{\epsilon\to0}\Bigl[\|\chi_{>R}^{\epsilon,+}\,w^{\epsilon,+}\|_{C_TH^m} + \|\chi_{>R}^{\epsilon,-}\,w^{\epsilon,-}\|_{C_TH^m}\Bigr].
\end{align*}
We may then pass to the limit as \(R\to\infty\) using \eqref{tight} to obtain
\[
\lim_{\epsilon \to 0}\int_{\nu\leq |t|\leq T}\|\Lambda_\epsilon^{-1}\p_x(w^{\epsilon,+}\,w^{\epsilon,-})\|_{H^s}\,dt = 0.
\]
Combining this with \eqref{short time}, we arrive at
\[
\limsup_{\epsilon\to0}\int_{-T}^T\|\Lambda_\epsilon^{-1}\p_x(w^{\epsilon,+}\,w^{\epsilon,-})\|_{H^s}\,dt \lesssim \nu A^2,
\]
and \eqref{cancellation} now follows from passing to the limit as \(\nu\to 0\).
\epf

Our final preparatory lemma approximates \(\Lambda_\epsilon\) by differential operators obtained from its Taylor expansion:
\begin{lem}\label{l:Taylor}
If \(\sigma\in \R\) then we have the estimates
\begin{align}
\|\tfrac1{\epsilon^2}\Lambda_\epsilon f_x - \tfrac 2{\sqrt 3 \epsilon^2}f_x + \tfrac1{4\sqrt 3}f_{xxx}\|_{H^\sigma} & \lesssim \|f\|_{H^{\sigma+3}},\label{linear convergence 2}\\
\|\tfrac1{\epsilon^2}\Lambda_\epsilon f_x - \tfrac 2{\sqrt 3 \epsilon^2}f_x + \tfrac1{4\sqrt 3}f_{xxx}\|_{H^\sigma} & \lesssim \epsilon^2 \|f\|_{H^{\sigma+5}},\label{linear convergence}\\
\|\Lambda_\epsilon f - \tfrac2{\sqrt 3}f\|_{H^\sigma} &\lesssim \epsilon^2\|f\|_{H^{\sigma+2}},\label{nonlinear convergence-1}\\
\|\Lambda_\epsilon^{-1}f - \tfrac{\sqrt 3}2f\|_{H^\sigma} &\lesssim \epsilon^2\|f\|_{H^{\sigma+2}}.\label{nonlinear convergence-2}
\end{align}
\end{lem}
\bpf
We recall that \(\Lambda_\epsilon\) has symbol \(\frac1{\sqrt 3}\sqrt{4 + \epsilon^2\xi^2}\). All of the desired estimates now follow from Plancherel's Theorem after computing that
\begin{align*}
\tfrac1{\sqrt 3\epsilon^2}\xi \sqrt{4+\epsilon^2\xi^2} - \tfrac2{\sqrt 3\epsilon^2}\xi - \tfrac1{4\sqrt 3}\xi^3 &= -\tfrac1{4\sqrt 3}\tfrac{\epsilon^2\xi^5}{8 + 4\sqrt{4+ \epsilon^2\xi^2} + \epsilon^2\xi^2},\\
\tfrac1{\sqrt 3}\sqrt{4+\epsilon^2\xi^2} - \tfrac2{\sqrt 3} &= \tfrac1{\sqrt 3}\tfrac{\epsilon^2\xi^2}{2 + \sqrt{4+\epsilon^2\xi^2}},\\
\tfrac{\sqrt 3}{\sqrt{4+\epsilon^2\xi^2}} - \tfrac{\sqrt 3}2 &= -\tfrac{\sqrt 3}2 \tfrac{\epsilon^2\xi^2}{4 + 2\sqrt{4+\epsilon^2\xi^2} + \epsilon^2\xi^2}.\qedhere
\end{align*}
\epf

We are now in a position to show that the solutions of \eqref{KdV} are approximate solutions of \eqref{GB-system-2}:

\begin{prop}
Let \(s\geq -\frac12\), \(m\geq \lceil s+5\rceil\) be an integer, and \(0<\epsilon_0\leq 1\). Suppose that for each \(0<\epsilon\leq \epsilon_0\) we have solutions \(w^{\epsilon,\pm}\in C(\R;H^m)\) of \eqref{KdV} that satisfy \eqref{KdV init cvgce} with \(\sigma = m\) for some \(w_0^\pm\in H^m\). Then
\eq{bv def}{
\bv^{\epsilon,\pm} := w^{\epsilon,\pm} \mp \tfrac{\epsilon^2}4\bigl[(w^{\epsilon,+})^2 - (w^{\epsilon,-})^2\bigr]\in C(\R;H^m)
}
satisfy
\begin{align}
\lim_{\epsilon \to0}\sup_{t\in \R}\|(w^{\epsilon,+},w^{\epsilon,-})(t) - (\bv^{\epsilon,+},\bv^{\epsilon,-})(t)\|_{H^m\times H^m} =0,\label{from u to bv}
\end{align}
and are approximate \(H^s\)-solutions of \eqref{GB-system-2}, in the sense that
\eq{bv eqn}{
\bv_t^{\epsilon,\pm} =  \mp \tfrac1{\epsilon^2}\Lambda_\epsilon \bv_x^{\epsilon,\pm} \pm \tfrac23\Lambda_\epsilon^{-1}\bigl[(\bv^{\epsilon,+} + \bv^{\epsilon,-})^2\bigr]_x \pm \err^{\epsilon,\pm},
}
where
\begin{align}
\limsup_{\epsilon\to0}\sup_{t\in \R}\|\err^{\epsilon,\pm}(t)\|_{H^1}<\infty,\label{crude err bound}
\end{align}
and for all \(T>0\) we have
\begin{align}
\lim_{\epsilon\to 0}\int_{-T}^T\|\err^{\epsilon,\pm}(t)\|_{H^s}\,dt = 0.\label{err bound}
\end{align}
\end{prop}
\bpf
The fact that \(\bv^{\epsilon,\pm}\in C(\R;H^m)\) follows immediately from the definition \eqref{bv def} and that \(H^m\) is an algebra, as does the estimate
\begin{align*}
\sup_{t\in \R}\|(w^{\epsilon,+},w^{\epsilon,-})(t) - (\bv^{\epsilon,+},\bv^{\epsilon,-})(t)\|_{H^m\times H^m} &\lesssim \epsilon^2 \sup_{t\in \R}\|(w^{\epsilon,+},w^{\epsilon,-})(t)\|_{H^m\times H^m}^2,
\end{align*}
which may be combined with \eqref{uniform consequence KdV} to yield \eqref{from u to bv}.

To prove \eqref{err bound}, we start by computing that
\begin{align*}
\err^{\epsilon,\pm} &= \Bigl[\tfrac1{\epsilon^2}\Lambda_\epsilon w_x^{\epsilon,\pm} - \tfrac 2{\sqrt 3\epsilon^2}w_x^{\epsilon,\pm} + \tfrac1{4\sqrt 3}w_{xxx}^{\epsilon,\pm}\Bigr] \mp \tfrac14 (\Lambda_\epsilon - \tfrac 2{\sqrt 3})\bigl[(w^{\epsilon,+})^2 - (w^{\epsilon,-})^2\bigr]_x\\
&\quad - \tfrac23(\Lambda_\epsilon^{-1} - \tfrac{\sqrt 3}2)\bigl[(w^{\epsilon,+})^2 + (w^{\epsilon,-})^2\bigr]_x - \tfrac43\Lambda_\epsilon^{-1}\bigl[w^{\epsilon,+}w^{\epsilon,-}\bigr]_x\\
&\quad - \tfrac{\epsilon^2}2\Bigl[w^{\epsilon,+}\bigl(w_t^{\epsilon,+}+\tfrac2{\sqrt 3\epsilon^2}w_x^{\epsilon,+}\bigr) - w^{\epsilon,-}\bigl(w_t^{\epsilon,-} - \tfrac2{\sqrt 3\epsilon^2}w_x^{\epsilon,-}\bigr)\Bigr].
\end{align*}

Applying \eqref{uniform consequence KdV}, we may find some \(0<\epsilon_1\leq \epsilon_0\) so that
\begin{align*}
A:= \sup_{0<\epsilon\leq \epsilon_1}\|(w^{\epsilon,+},w^{\epsilon,-})\|_{C_tH^m}<\infty.
\end{align*}
Taking \(0<\epsilon\leq \epsilon_1\), we may then apply Lemma~\ref{l:Taylor} with the fact that \(H^m\) is an algebra to estimate
\begin{align*}
\bigl\|\tfrac1{\epsilon^2}\Lambda_\epsilon w_x^{\epsilon,\pm} - \tfrac 2{\sqrt 3\epsilon^2}w_x^{\epsilon,\pm} + \tfrac1{4\sqrt 3}w_{xxx}^{\epsilon,\pm}\bigr\|_{C_tH^1} &\lesssim A,\\
\bigl\|\tfrac1{\epsilon^2}\Lambda_\epsilon w_x^{\epsilon,\pm} - \tfrac 2{\sqrt 3\epsilon^2}w_x^{\epsilon,\pm} + \tfrac1{4\sqrt 3}w_{xxx}^{\epsilon,\pm}\Bigr\|_{C_tH^s} &\lesssim \epsilon^2 A,\\
\bigl\|(\Lambda_\epsilon - \tfrac 2{\sqrt 3})\bigl[(w^{\epsilon,+})^2 - (w^{\epsilon,-})^2\bigr]_x\bigr\|_{C_tH^{\max\{s,1\}}}&\lesssim \epsilon^2 A^2,\\
\bigl\|(\Lambda_\epsilon^{-1} - \tfrac{\sqrt 3}2)\bigl[(w^{\epsilon,+})^2 + (w^{\epsilon,-})^2\bigr]_x\bigr\|_{C_tH^{\max\{s,1\}}}&\lesssim \epsilon^2 A^2.
\end{align*}
Again using that \(H^m\) is an algebra we have
\[
\bigl\|\Lambda_\epsilon^{-1}\bigl[w^{\epsilon,+}w^{\epsilon,-}\bigr]_x\bigr\|_{C_tH^{\max\{s,1\}}}\lesssim A^2,
\]
and also using that \(w^{\epsilon,\pm}\) solve \eqref{KdV} we have
\[
\|w_t^{\epsilon,\pm} \pm \tfrac2{\sqrt 3\epsilon^2}w_x^{\epsilon,\pm}\|_{C_tH^{\max\{s,1\}}}\lesssim A(1+A).
\]
Combining this last estimate with \eqref{product-1} and Sobolev embedding then yields
\begin{align*}
&\bigl\|w^{\epsilon,\pm}\bigl(w_t^{\epsilon,\pm}\pm\tfrac2{\sqrt 3\epsilon^2}w_x^{\epsilon,\pm}\bigr)\bigr\|_{C_tH^{\max\{s,1\}}}\\
&\qquad\lesssim \|w^{\epsilon,\pm}\|_{C_tW^{\max\{\lceil s\rceil,1\},\infty}} \times\|w_t^{\epsilon,\pm} \pm \tfrac2{\sqrt 3\epsilon^2}w_x^{\epsilon,\pm}\|_{C_tH^{\max\{s,1\}}}\lesssim A^2(1+A).
\end{align*}

Putting these estimates together, we see that for all \(0<\epsilon\leq \epsilon_1\) we may bound
\[
\|\err^{\epsilon,\pm}\|_{C_tH^1}\lesssim A(1+A)^2,
\]
from which we obtain \eqref{crude err bound}. Further, after integrating in time we have
\[
\int_{-T}^T\|\err^{\epsilon,\pm}(t)\|_{H^s}\,dt\lesssim \int_{-T}^T\bigl\|\Lambda_\epsilon^{-1}\bigl[w^{\epsilon,+}w^{\epsilon,-}\bigr]_x\bigr\|_{H^s}\,dt + \epsilon^2TA(1+A)^2.
\]
The identity \eqref{err bound} now follows from \eqref{cancellation}.
\epf

Now we have shown that sufficiently smooth solutions of \eqref{KdV} can be used to construct approximate solutions of \eqref{GB-system-2}, we compare our approximate solution to a genuine solution:

\begin{prop}\label{p:modified energy}
Let \(s\geq -\frac12\), \(m \geq \lceil s+5\rceil\), \(0<\epsilon_0\leq 1\), and \(w_0^\pm\in H^m\). Suppose that for all \(0<\epsilon\leq \epsilon_0\) we have solutions \(w^{\epsilon,\pm}\in C(\R;H^m)\) of \eqref{KdV} satisfying \eqref{KdV init cvgce} with \(\sigma = m\), as well as solutions \(v^{\epsilon,\pm}\in C(\R;H^m)\) of \eqref{GB-system-2} so that
\eq{init is close}{
\lim_{\epsilon \to 0}\|(v^{\epsilon,+},v^{\epsilon,-})(0) - (w_0^+,w_0^-)\|_{H^m\times H^m} = 0.
}
Then for all times \(T>0\) we have
\eq{cvgce}{
\lim_{\epsilon \to 0}\sup_{|t|\leq T}\|(v^{\epsilon,+},v^{\epsilon,-})(t) - (w^{\epsilon,+},w^{\epsilon,-})(t)\|_{H^s\times H^s} = 0.
}
\end{prop}
\bpf
From \eqref{uniform consequence KdV}, we may find some \(R>0\) and \(0<\epsilon_1\leq \epsilon_0\) so that
\eq{state}{
\sup_{0<\epsilon\leq \epsilon_1}\sup_{t\in \R}\|(w^{\epsilon,+},w^{\epsilon,-})(t)\|_{H^{-\frac12}\times H^{-\frac12}}\leq R.
}

We now suppose that for some \(T_0\geq 0\) we have
\eq{boot}{
\lim_{\epsilon\to 0}\sup_{|t|\leq T_0}\|(v^{\epsilon,+},v^{\epsilon,-})(t) - (w^{\epsilon,+},w^{\epsilon,-})(t)\|_{H^s\times H^s} = 0.
}
This evidently holds for \(T_0 = 0\) due to \eqref{KdV init cvgce} and \eqref{init is close}.  We will prove that there exists some \(T_1>T_0\) so that \(T_1 - T_0\simeq \<R\>^{-3}\) and
\[
\lim_{\epsilon\to 0}\sup_{|t|\leq T_1}\|(v^{\epsilon,+},v^{\epsilon,-})(t) - (w^{\epsilon,+},w^{\epsilon,-})(t)\|_{H^s\times H^s} = 0.
\]
The identity \eqref{cvgce} then follows from iterating this argument.

From \eqref{state} and \eqref{boot}, after possibly shrinking \(\epsilon_1\) we may assume that \(\epsilon_1\leq \tfrac1{2R}\) and
\begin{align*}
\sup_{0<\epsilon\leq \epsilon_1}\sup_{|t|\leq T_0}\|(v^{\epsilon,+},v^{\epsilon,-})(t)\|_{H^{-\frac12}\times H^{-\frac12}}&\leq 2R,\\
\sup_{0<\epsilon\leq \epsilon_1}\|(v^{\epsilon,+},v^{\epsilon,-})(0)\|_{H^m\times H^m}&<\infty.
\end{align*}
Iteratively applying Theorem~\ref{t:LWP-details} on time intervals of length \(\simeq \<R\>^{-3}\), we may then find some \(T_1>T_0\) so that \(T_1 - T_0\simeq \<R\>^{-3}\) and
\[
\sup_{0<\epsilon\leq \epsilon_1}\sup_{|t|\leq T_1}\|(v^{\epsilon,+},v^{\epsilon,-})(t)\|_{H^m\times H^m}<\infty.
\]
After possibly shrinking the size of \(\epsilon_1\) even further, we may combine this with \eqref{uniform consequence KdV}, \eqref{from u to bv}, and \eqref{crude err bound} to ensure that, writing \(v^\epsilon = (v^{\epsilon,+},v^{\epsilon,-})\), etc.,
\begin{align*}
A := \sup_{0<\epsilon\leq \epsilon_1}\sup_{|t|\leq T_1}\Bigl[&\|v^\epsilon(t)\|_{H^m\times H^m} + \|w^\epsilon(t)\|_{H^m\times H^m}\\
&+ \|\bv^\epsilon(t)\|_{H^m\times H^m} + \|\err^\epsilon(t)\|_{H^1\times H^1}\Bigr]<\infty.
\end{align*}

For \(0<\epsilon\leq \epsilon_1\), let us define
\eq{V defn}{
V^{\epsilon,\pm} := v^{\epsilon,\pm} - \bv^{\epsilon,\pm},
}
which satisfy
\eq{V diff eq}{
V_t^{\epsilon,\pm} = \mp \tfrac1{\epsilon^2}\Lambda_\epsilon V_x^{\epsilon,\pm} \pm\tfrac1{\sqrt3}\bigl[W^\epsilon(V^{\epsilon,+} + V^{\epsilon,-})\bigr]_x \pm \Err^{\epsilon,\pm},
}
where
\[
W^\epsilon = v^{\epsilon,+} + v^{\epsilon,-} + \bv^{\epsilon,+} + \bv^{\epsilon,-},
\]
with error terms
\begin{align*}
\Err^{\epsilon,\pm} &:= \tfrac23\bigl(\Lambda_\epsilon^{-1} - \tfrac{\sqrt 3}2\bigr)\bigl[W^\epsilon(V^{\epsilon,+} + V^{\epsilon,-})\bigr]_x - \err^{\epsilon,\pm}.
\end{align*}

Applying \eqref{nonlinear convergence-2} with the fact that \(H^m\) is an algebra, we have
\[
\int_{-T_1}^{T_1} \|\Err^{\epsilon,\pm}(t)\|_{H^s}\,dt\lesssim \int_{-T_1}^{T_1} \|\err^{\epsilon,\pm}(t)\|_{H^s}\,dt + \epsilon^2T_1A^2.
\]
The identity \eqref{err bound} then gives us
\eq{Err bound}{
\lim_{\epsilon\to 0}\int_{-T_1}^{T_1}\|\Err^{\epsilon,\pm}(t)\|_{H^s}\,dt = 0.
}

Using that both \(H^m\) and \(H^{m-3}\) are algebras, that \(v^{\epsilon,\pm},w^{\epsilon,\pm}\in C(\R;H^m)\), the equations \eqref{KdV} and \eqref{GB-system-2}, and the definition \eqref{bv def}, we conclude that \(V^{\epsilon,\pm}\in C(\R;H^m)\cap C^1(\R;H^{m-3})\). Writing \(V^\epsilon = (V^{\epsilon,+},V^{\epsilon,-})\), we then integrate by parts to obtain
\begin{align*}
\tfrac d{dt}\tfrac12\|V^\epsilon\|_{H^s\times H^s}^2 &= \tfrac1{\sqrt 3}\int_\R \<\p_x\>^sV^{\epsilon,+}\times \<\p_x\>^s\bigl(W_x^\epsilon (V^{\epsilon,+}+V^{\epsilon,-})\bigr)\,dx\\
&\quad + \tfrac1{\sqrt 3}\int_\R \<\p_x\>^sV^{\epsilon,+}\times [\<\p_x\>^s,W^\epsilon](V^{\epsilon,+}+V^{\epsilon,-})_x\,dx\\
&\quad - \tfrac1{2\sqrt 3}\int_\R W_x^\epsilon |\<\p_x\>^sV^{\epsilon,+}|^2\,dx + \tfrac1{\sqrt 3}\int_\R W^\epsilon\times\<\p_x\>^s V^{\epsilon,+}\times \<\p_x\>^sV_x^{\epsilon,-}\,dx\\
&\quad + \int_\R\<\p_x\>^s V^{\epsilon,+}\times \<\p_x\>^s\Err^{\epsilon,+}\,dx\\
&\quad -  \tfrac1{\sqrt 3}\int_\R \<\p_x\>^sV^{\epsilon,-}\times \<\p_x\>^s\bigl(W_x^\epsilon (V^{\epsilon,+}+V^{\epsilon,-})\bigr)\,dx\\
&\quad - \tfrac1{\sqrt 3}\int_\R \<\p_x\>^sV^{\epsilon,-}\times [\<\p_x\>^s,W^\epsilon](V^{\epsilon,+}+V^{\epsilon,-})_x\,dx\\
&\quad + \tfrac1{2\sqrt 3}\int_\R W_x^\epsilon |\<\p_x\>^sV^{\epsilon,-}|^2\,dx - \tfrac1{\sqrt 3}\int_\R W^\epsilon\times\<\p_x\>^s V^{\epsilon,-}\times \<\p_x\>^sV_x^{\epsilon,+}\,dx\\
&\quad - \int_\R\<\p_x\>^s V^{\epsilon,-}\times \<\p_x\>^s\Err^{\epsilon,-}\,dx.
\end{align*}
Employing \eqref{product-1}, Sobolev embedding, and \eqref{commutator}, we may then bound
\begin{align*}
&\left|\tfrac d{dt}\tfrac12\|V^\epsilon\|_{H^s\times H^s}^2 - \tfrac1{\sqrt 3}\int_\R W^\epsilon\times\Bigl[\<\p_x\>^s V^{\epsilon,+}\times \<\p_x\>^sV_x^{\epsilon,-} - \<\p_x\>^s V^{\epsilon,-}\times \<\p_x\>^sV_x^{\epsilon,+}\Bigr]\,dx\right|\\
&\qquad \lesssim A\|V^\epsilon\|_{H^s\times H^s}^2 + \|V^\epsilon\|_{H^s\times H^s}\|\Err^\epsilon\|_{H^s\times H^s}.
\end{align*}

To handle the remaining term, we compute that
\begin{align*}
&\tfrac d{dt}\tfrac{\epsilon^2}4\int_\R W^\epsilon |\<\p_x\>^s(V^{\epsilon,+} + V^{\epsilon,-})|^2\,dx\\
&\qquad = \tfrac 1{\sqrt 3}\int_\R W^\epsilon \times \Bigl[\<\p_x\>^sV^{\epsilon,+}\times\<\p_x\>^sV_x^{\epsilon,-} - \<\p_x\>^sV^{\epsilon,-}\times \<\p_x\>^sV_x^{\epsilon,+}\Bigr]\,dx\\
&\qquad\quad + \tfrac 1{2\sqrt 3}\int_\R W_x^\epsilon \times |\<\p_x\>^sV^{\epsilon,+}|^2\,dx - \tfrac 1{2\sqrt 3}\int_\R W_x^\epsilon \times |\<\p_x\>^sV^{\epsilon,-}|^2\,dx\\
&\qquad\quad + \tfrac{\epsilon^2}4\int_\R W_t^\epsilon |\<\p_x\>^s(V^{\epsilon,+} + V^{\epsilon,-})|^2\,dx\\
&\qquad\quad + \tfrac12\int_\R W^\epsilon \times \<\p_x\>^s(V^{\epsilon,+} + V^{\epsilon,-})\times \<\p_x\>^s\bigl(\epsilon^2 V_t^{\epsilon,+} + \tfrac2{\sqrt 3}V_x^{\epsilon,+}\bigr)\,dx\\
&\qquad\quad + \tfrac12\int_\R W^\epsilon \times \<\p_x\>^s(V^{\epsilon,+} + V^{\epsilon,-})\times \<\p_x\>^s\bigl(\epsilon^2 V_t^{\epsilon,-} - \tfrac2{\sqrt 3}V_x^{\epsilon,-}\bigr)\,dx.
\end{align*}

From \eqref{GB-system-2}, \eqref{bv eqn}, the fact that \(H^m\) is an algebra, and Sobolev embedding we have
\[
\|\epsilon^2 W_t^\epsilon\|_{L^\infty}\lesssim\|\epsilon^2 W_t^\epsilon\|_{H^1}\lesssim A(1+A).
\]
Similarly, by using \eqref{V diff eq} and \eqref{nonlinear convergence-1} we have
\begin{align*}
\|\epsilon^2V_t^{\epsilon,\pm}\pm \tfrac2{\sqrt 3}V_x^{\epsilon,\pm}\|_{H^s}&\lesssim \epsilon^2 A(1+A) + \epsilon^2 \|\Err^\epsilon\|_{H^s\times H^s}.
\end{align*}
It follows that
\begin{align*}
&\left|\tfrac d{dt}\biggl[\tfrac12\|V^\epsilon\|_{H^s\times H^s}^2 - \tfrac{\epsilon^2}4\int_\R W^\epsilon |\<\p_x\>^s(V^{\epsilon,+} + V^{\epsilon,-})|^2\,dx\biggr]\right|\\
&\qquad\lesssim A(1+A)\|V^\epsilon\|_{H^s\times H^s}^2 + \|V^\epsilon\|_{H^s\times H^s}\Bigl[\epsilon^2A^2(1+A) + (1+\epsilon^2A)\|\Err^\epsilon\|_{H^s\times H^s}\Bigr].
\end{align*}
Noting that
\[
\left|\tfrac{\epsilon^2}4\int_\R W^\epsilon |\<\p_x\>^s(V^{\epsilon,+} + V^{\epsilon,-})|^2\,dx\right|\lesssim \epsilon^2A\|V^\epsilon\|_{H^s\times H^s}^2,
\]
we may then apply Gronwall's inequality to bound
\begin{align*}
\sup_{|t|\leq T_1}\|V^\epsilon(t)\|_{H^s\times H^s}^2 &\lesssim_{A,T_1} \|V^\epsilon(0)\|_{H^s\times H^s}^2 + \epsilon^2\sup_{|t|\leq T_1}\|V^\epsilon(t)\|_{H^s\times H^s}^2\\
&\qquad+ \sup_{|t|\leq T_1}\|V^\epsilon(t)\|_{H^s\times H^s}\left[\epsilon^2 +\int_{-T_1}^{T_1}\|\Err^\epsilon(t)\|_{H^s\times H^s}\,dt\right].
\end{align*}
For sufficiently small \(0<\epsilon \ll\epsilon_1\) we may then apply Cauchy's inequality to obtain
\[
\sup_{|t|\leq T_1}\|V^\epsilon(t)\|_{H^s\times H^s}\lesssim_{A,T_1}\|V^\epsilon(0)\|_{H^s\times H^s} + \epsilon^2 + \int_{-T_1}^{T_1}\|\Err^\epsilon(t)\|_{H^s\times H^s}\,dt,
\]
and \eqref{cvgce} now follows from the definition \eqref{V defn} and the identities \eqref{from u to bv}, \eqref{Err bound}.
\epf

At long last, we are in a position to complete the

\bpf[Proof of Theorem~\ref{t:long wave}] We start by considering part (i). Let \((q_0,r_0) := (u_0,\p_x^{-1}u_1) \in H^s\times H^{s-1}\).  Let \(m = \lceil s+5\rceil\) and choose a sequence \((q_{0,n},r_{0,n})\in H^m\times H^{m-1}\) converging to \((q_0,r_0)\) in \(H^s\times H^{s-1}\).  After passing to a tail of the sequence we may assume that
\eq{uniform init bound}{
\sup_{n\geq1}\|(q_{0,n},r_{0,n})\|_{H^s\times H^{s-1}}\leq 1 + \|(q_0,r_0)\|_{H^s\times H^{s-1}}.
}

Using the change of variables \eqref{d:veps}, we take
\begin{align*}
v_0^{\epsilon,\pm} := \tfrac12\bigl(q_0 \mp \epsilon^2\Lambda_\epsilon^{-1}r_0\bigr)\qtq{and}v_{0,n}^{\epsilon,\pm} := \tfrac12\bigl(q_{0,n} \mp \epsilon^2\Lambda_\epsilon^{-1}r_{0,n}\bigr).
\end{align*}
It follows that
\eq{uniform init lim}{
\lim_{\epsilon \to 0}\|v_0^\epsilon - \tfrac12(q_0,q_0)\|_{H^s\times H^s} = 0 = \lim_{\epsilon \to 0}\sup_{n\geq 1}\|v_{0,n}^\epsilon - \tfrac12(q_{0,n},q_{0,n})\|_{H^s\times H^s},
}
where we continue to use the convention that \(v_0^\epsilon = (v_0^{\epsilon,+},v_0^{\epsilon,-})\) and \(v_{0,n}^\epsilon = (v_{0,n}^{\epsilon,+},v_{0,n}^{\epsilon,-})\). Further, for each fixed \(n\geq 1\) we have
\eq{prop former}{
\lim_{\epsilon \to 0}\|v_{0,n}^\epsilon - \tfrac12(q_{0,n},q_{0,n})\|_{H^m\times H^m} = 0.
}

From \eqref{uniform init bound}, \eqref{uniform init lim}, and Corollary~\ref{cor:rescaled well-posedness}, we may find some \(0<\epsilon_0\leq 1\) so that for all \(0<\epsilon\leq \epsilon_0\) and \(n\geq 1\) we have solutions \(v^{\epsilon,\pm}\in C(\R;H^s)\) and \(v_n^{\epsilon,\pm}\in C(\R;H^m)\) of \eqref{GB-system-2} with initial data \(v^{\epsilon,\pm}(0) = v_0^{\epsilon,\pm}\) and \(v_n^{\epsilon,\pm}(0) = v_{0,n}^{\epsilon,\pm}\). Similarly, we may use \cite{MR1086966,MR1969209} to find solutions \(w^{\epsilon,\pm}\in C(\R;H^s)\) and \(w_n^{\epsilon,\pm}\in C(\R;H^m)\) of \eqref{KdV} with initial data \(w^{\epsilon,\pm}(0) = v_0^{\epsilon,\pm}\) and \(w_n^{\epsilon,\pm}(0) = v_{0,n}^{\epsilon,\pm}\).

The identity \eqref{prop former} allows us to apply Proposition~\ref{p:modified energy} to show that for all \(n\geq 1\) and \(T>0\) we have
\eq{middle}{
\lim_{\epsilon\to 0}\sup_{|t|\leq T}\|v_n^\epsilon(t) - w_n^\epsilon(t)\|_{H^s\times H^s} = 0.
}

For any sequence \(n_j\to\infty\) and \(\epsilon_0\geq \epsilon_j\to 0\) we may use \eqref{uniform init lim} to obtain
\begin{align*}
\lim_{j\to\infty}\|v_{0,n_j}^{\epsilon_j} - v_0^{\epsilon_j}\|_{H^s\times H^s} &\lesssim \lim_{j\to\infty}\sup_{n\geq 1}\|v_{0,n}^{\epsilon_j} - \tfrac12(q_{0,n},q_{0,n})\|_{H^s\times H^s} + \lim_{j\to\infty}\|q_{0,n_j} - q_0\|_{H^s}\\
&\quad + \lim_{j\to\infty}\|\tfrac12(q_0,q_0) - v_0^{\epsilon_j}\|_{H^s\times H^s} \\
&= 0.
\end{align*}
We may then use the change of variables \eqref{static rescaling} and continuity of the solution map for \eqref{KdV stat} to prove that
\[
\lim_{j\to\infty}\sup_{|t|\leq T}\|w_{n_j}^{\epsilon_j}(t) - w^{\epsilon_j}(t)\|_{H^s\times H^s} = 0.
\]
This yields the identity
\eq{right end}{
\lim_{n\to\infty}\limsup_{\epsilon\to0}\sup_{|t|\leq T}\|w_n^\epsilon(t) - w^\epsilon(t)\|_{H^s\times H^s} = 0.
}

To complete the proof, we will employ a similar bootstrap argument to Proposition~\ref{p:modified energy}.  From \eqref{uniform consequence KdV} and \eqref{uniform init lim}, we may find some \(0<\epsilon_1\leq \epsilon_0\) and \(R>0\) so that
\eq{ap boots}{
\sup_{0<\epsilon\leq \epsilon_1}\sup_{t\in \R}\|w^\epsilon(t)\|_{H^{-\frac12}\times H^{-\frac12}}\leq R.
}
We then suppose that for some \(T_0\geq 0\) we have
\eq{boots}{
\lim_{n\to\infty}\limsup_{\epsilon \to 0}\sup_{|t|\leq T_0}\|v^\epsilon(t) - v_n^\epsilon(t)\|_{H^s\times H^s} = 0.
}
It follows from \eqref{uniform init lim} that \eqref{boots} is valid for \(T_0 = 0\). We will once again prove that there exists some \(T_1>T_0\) so that \(T_1 - T_0\simeq \<R\>^{-3}\) and
\eq{straps}{
\lim_{n\to\infty}\limsup_{\epsilon \to 0}\sup_{|t|\leq T_1}\|v^\epsilon(t) - v_n^\epsilon(t)\|_{H^s\times H^s} = 0.
}
Iterating this argument, for all \(T>0\) we have
\eq{left end}{
\lim_{n\to\infty}\limsup_{\epsilon \to 0}\sup_{|t|\leq T}\|v^\epsilon(t) - v_n^\epsilon(t)\|_{H^s\times H^s} = 0.
}
We may then apply the triangle inequality and \eqref{middle} to show that for all \(n\geq 1\) we have
\begin{align}
\limsup_{\epsilon\to0}\sup_{|t|\leq T}\|v^\epsilon(t) - w^\epsilon(t)\|_{H^s\times H^s}&\leq \limsup_{\epsilon\to0}\sup_{|t|\leq T}\|v^\epsilon(t) - v_n^\epsilon(t)\|_{H^s\times H^s}\label{triangle}\\
&\quad + \limsup_{\epsilon\to0}\sup_{|t|\leq T}\|w_n^\epsilon(t) - w^\epsilon(t)\|_{H^s\times H^s},\notag
\end{align}
and the identity \eqref{long wave approx} then follows from taking \(n\to\infty\) using \eqref{left end} and \eqref{right end}.

To prove \eqref{straps}, we first note that using \eqref{boots} in place of \eqref{left end} in \eqref{triangle} proves that
\[
\lim_{\epsilon\to0}\sup_{|t|\leq T_0}\|v^\epsilon(t) - w^\epsilon(t)\|_{H^s\times H^s} = 0.
\]
It then follows from \eqref{ap boots} and \eqref{boots} that for each \(n\geq 2\) we may find some \(0<\epsilon_n\leq \min\{\epsilon_1,\frac1{2R}\}\) so that, after possibly passing to a tail of our sequence, we have
\[
\sup_{n\geq 2}\sup_{0<\epsilon\leq \epsilon_n}\sup_{|t|\leq T_0}\|v_n^\epsilon(t)\|_{H^{-\frac12}\times H^{-\frac12}}\leq 2R.
\]
We may then iteratively apply Theorem~\ref{t:LWP-details} on timescales of size \(\simeq \<R\>^{-3}\) to show that there exists some time \(T_1>T_0\) such that \(T_1 - T_0\simeq\<R\>^{-3}\) and for all \(n\geq 2\) we have
\eq{final iteration}{
\sup_{|t|\leq T_1}\|v_n^\epsilon(t) - v^\epsilon(t)\|_{H^s\times H^s}\leq C\bigl(T_1,R,\|(q_0,r_0)\|_{H^s\times H^{s-1}}\bigr)\|v_{0,n}^\epsilon - v_0^\epsilon\|_{H^s\times H^s},
}
uniformly for all \(0<\epsilon\leq \epsilon_n\).  The estimate \eqref{straps} now follows from taking \(\epsilon\to0\) and then \(n\to\infty\).

Part (ii) is proved in an almost identical manner, first taking sequences \(w_{0,n}^\pm\in H^m\) converging to \(w_0^\pm \in H^s\) so that
\[
\sup_{n\geq 1}\|(w_{0,n}^+,w_{0,n}^-)\|_{H^s\times H^s}\leq 1 + \|(w_0^+,w_0^-)\|_{H^s\times H^s},
\]
and then setting \(v_0^{\epsilon,\pm} := w_0^\pm\), \(v_{0,n}^{\epsilon,\pm} := w_{0,n}^\pm\).  Inverting \eqref{d:veps}, we obtain data
\begin{align*}
(q_0^\epsilon,r_0^\epsilon) &:= \bigl(w_0^++w_0^-,\epsilon^{-2}\Lambda_\epsilon(w_0^--w_0^+)\bigr),\\
(q_{0,n}^\epsilon,r_{0,n}^\epsilon) &:=\bigl(w_{0,n}^++w_{0,n}^-, \epsilon^{-2}\Lambda_\epsilon(w_{0,n}^--w_{0,n}^+)\bigr),
\end{align*}
for \eqref{GB-system rescaled}. Rescaling, we obtain
\[
\|(\epsilon^2 q_0^\epsilon,\epsilon^4 r_0^\epsilon)(\epsilon x)\|_{H^{-\frac12}\times H^{-\frac32}}\lesssim \epsilon \|(w_0^+,w_0^-)\|_{H^{-\frac12}\times H^{-\frac12}},
\]
with a similar estimate for the rescaling of \((q_{0,n}^\epsilon,r_{0,n}^\epsilon)\). From Theorem~\ref{t:GWP}, this suffices to show that there exists some \(0<\epsilon_0\leq 1\) so that for all \(0<\epsilon \leq \epsilon_0\) the global solutions \(v^{\epsilon,\pm}\), \(v_n^{\epsilon,\pm}\) of \eqref{GB-system-2} with corresponding data \(v_0^{\epsilon,\pm}\), \(v_{0,n}^{\epsilon,\pm}\) are well-defined.

The proof is now completed by applying an identical argument to part (i) after replacing \(\frac12(q_0,q_0)\) by \((w_0^+,w_0^-)\) and \(\frac12(q_{0,n},q_{0,n})\) by \((w_{0,n}^+,w_{0,n}^-)\) in \eqref{uniform init lim} and \eqref{prop former}, and taking the constant in \eqref{final iteration} to depend on \(\|(w_0^+,w_0^-)\|_{H^s\times H^s}\) instead of \(\|(q_0,r_0)\|_{H^s\times H^{s-1}}\).
\epf

\bibliographystyle{habbrv}
\bibliography{refs}

\begin{thebibliography}{10}

\bibitem{MR4747862}
H.~Bahouri, T.~M. Leslie, and G.~Perelman.
\newblock {$H^s$} bounds for the derivative nonlinear {S}chr\"odinger equation.
\newblock {\em Math. Res. Lett.}, 30(5):1299--1333, 2023.

\bibitem{MR3953031}
R.~F. Barostichi, R.~O. Figueira, and A.~A. Himonas.
\newblock Well-posedness of the ``good'' {B}oussinesq equation in analytic
  {G}evrey spaces and time regularity.
\newblock {\em J. Differential Equations}, 267(5):3181--3198, 2019.

\bibitem{MR3922812}
R.~Bauer, W.-P. D\"ull, and G.~Schneider.
\newblock The {K}orteweg--de {V}ries, {B}urgers and {W}hitham limits for a
  spatially periodic {B}oussinesq model.
\newblock {\em Proc. Roy. Soc. Edinburgh Sect. A}, 149(1):191--217, 2019.

\bibitem{MR2204680}
I.~Bejenaru and T.~Tao.
\newblock Sharp well-posedness and ill-posedness results for a quadratic
  non-linear {S}chr\"odinger equation.
\newblock {\em J. Funct. Anal.}, 233(1):228--259, 2006.

\bibitem{MR954673}
J.~L. Bona and R.~L. Sachs.
\newblock Global existence of smooth solutions and stability of solitary waves
  for a generalized {B}oussinesq equation.
\newblock {\em Comm. Math. Phys.}, 118(1):15--29, 1988.

\bibitem{MR1209299}
J.~Bourgain.
\newblock Fourier transform restriction phenomena for certain lattice subsets
  and applications to nonlinear evolution equations. {I}. {S}chr\"odinger
  equations.
\newblock {\em Geom. Funct. Anal.}, 3(2):107--156, 1993.

\bibitem{MR1215780}
J.~Bourgain.
\newblock Fourier transform restriction phenomena for certain lattice subsets
  and applications to nonlinear evolution equations. {II}. {T}he {K}d{V}
  equation.
\newblock {\em Geom. Funct. Anal.}, 3(3):209--262, 1993.

\bibitem{MR3363411}
J.~Boussinesq.
\newblock Th\'eorie des ondes et des remous qui se propagent le long d'un canal
  rectangulaire horizontal, en communiquant au liquide contenu dans ce canal
  des vitesses sensiblement pareilles de la surface au fond.
\newblock {\em J. Math. Pures Appl. (2)}, 17:55--108, 1872.

\bibitem{MR4304314}
B.~Bringmann, R.~Killip, and M.~Vi\c{s}an.
\newblock Global well-posedness for the fifth-order {K}d{V} equation in
  {$H^{-1}(\Bbb R)$}.
\newblock {\em Ann. PDE}, 7(2):Paper No. 21, 46, 2021.

\bibitem{MR5130833}
C.~Charlier.
\newblock Blowup {S}olutions of the ``{B}ad'' {B}oussinesq {E}quation.
\newblock {\em SIAM J. Math. Anal.}, 58(5):5012--5048, 2026.

\bibitem{MR4481092}
C.~Charlier and J.~Lenells.
\newblock The ``good'' {B}oussinesq equation: a {R}iemann-{H}ilbert approach.
\newblock {\em Indiana Univ. Math. J.}, 71(4):1505--1562, 2022.

\bibitem{charlier2023boussinesqsequationwaterwaves}
C.~{C}harlier and J.~{L}enells.
\newblock On {B}oussinesq's equation for water waves.
\newblock \textit{Preprint}, 2023.
\newblock \href{https://arxiv.org/abs/2204.02365}{\texttt{arXiv:2204.02365}}.

\bibitem{MR4692724}
C.~Charlier and J.~Lenells.
\newblock Miura transformation for the ``good'' {B}oussinesq equation.
\newblock {\em Stud. Appl. Math.}, 152(1):73--110, 2024.

\bibitem{MR4632376}
C.~Charlier, J.~Lenells, and D.-S. Wang.
\newblock The ``good'' {B}oussinesq equation: long-time asymptotics.
\newblock {\em Anal. PDE}, 16(6):1351--1388, 2023.

\bibitem{MR2018661}
M.~Christ, J.~Colliander, and T.~Tao.
\newblock Asymptotics, frequency modulation, and low regularity ill-posedness
  for canonical defocusing equations.
\newblock {\em Amer. J. Math.}, 125(6):1235--1293, 2003.

\bibitem{MR1969209}
J.~Colliander, M.~Keel, G.~Staffilani, H.~Takaoka, and T.~Tao.
\newblock Sharp global well-posedness for {K}d{V} and modified {K}d{V} on
  {$\Bbb R$} and {$\Bbb T$}.
\newblock {\em J. Amer. Math. Soc.}, 16(3):705--749, 2003.

\bibitem{MR2411922}
P.~Deift.
\newblock Some open problems in random matrix theory and the theory of
  integrable systems.
\newblock In {\em Integrable systems and random matrices}, volume 458 of {\em
  Contemp. Math.}, pages 419--430. Amer. Math. Soc., Providence, RI, 2008.

\bibitem{PhysRevB.36.3031}
F.~Falk, E.~W. Laedke, and K.~H. Spatschek.
\newblock Stability of solitary-wave pulses in shape-memory alloys.
\newblock {\em Phys. Rev. B}, 36:3031--3041, 1987.

\bibitem{MR730182}
G.~E. Fal'kovich, M.~D. Spector, and S.~K. Turitsyn.
\newblock Destruction of stationary solutions and collapse in the nonlinear
  string equation.
\newblock {\em Phys. Lett. A}, 99(6-7):271--274, 1983.

\bibitem{MR1399198}
Y.-F. Fang and M.~G. Grillakis.
\newblock Existence and uniqueness for {B}oussinesq type equations on a circle.
\newblock {\em Comm. Partial Differential Equations}, 21(7-8):1253--1277, 1996.

\bibitem{MR2512853}
L.~G. Farah.
\newblock Local solutions in {S}obolev spaces with negative indices for the
  ``good'' {B}oussinesq equation.
\newblock {\em Comm. Partial Differential Equations}, 34(1-3):52--73, 2009.

\bibitem{MR2566562}
L.~G. Farah and M.~Scialom.
\newblock On the periodic ``good'' {B}oussinesq equation.
\newblock {\em Proc. Amer. Math. Soc.}, 138(3):953--964, 2010.

\bibitem{MR2756068}
L.~C.~F. Ferreira.
\newblock Existence and scattering theory for {B}oussinesq type equations with
  singular data.
\newblock {\em J. Differential Equations}, 250(5):2372--2388, 2011.

\bibitem{MR3200091}
L.~Grafakos and S.~Oh.
\newblock The {K}ato-{P}once inequality.
\newblock {\em Comm. Partial Differential Equations}, 39(6):1128--1157, 2014.

\bibitem{MR2531556}
Z.~Guo.
\newblock Global well-posedness of {K}orteweg-de {V}ries equation in
  {$H^{-3/4}(\Bbb R)$}.
\newblock {\em J. Math. Pures Appl. (9)}, 91(6):583--597, 2009.

\bibitem{MR4923538}
S.~Haque, R.~Killip, M.~Vi\c{s}an, and Y.~Zhang.
\newblock Global well-posedness and equicontinuity for modified {K}orteweg--de
  {V}ries equations in modulation spaces.
\newblock {\em Pure Appl. Anal.}, 7(3):615--637, 2025.

\bibitem{MR5031430}
B.~Harrop-Griffiths, R.~Killip, M.~Ntekoume, and M.~Vi\c{s}an.
\newblock Global well-posedness for the derivative nonlinear {S}chr\"odinger
  equation in {$L^2(\Bbb R)$}.
\newblock {\em J. Eur. Math. Soc. (JEMS)}, 28(2):843--924, 2026.

\bibitem{MR4726498}
B.~Harrop-Griffiths, R.~Killip, and M.~Vi\c{s}an.
\newblock Sharp well-posedness for the cubic {NLS} and m{K}d{V} in {$H^s(\Bbb
  R)$}.
\newblock {\em Forum Math. Pi}, 12:Paper No. e6, 86, 2024.

\bibitem{MR5083420}
B.~Harrop-Griffiths, R.~Killip, and M.~Vi\c{s}an.
\newblock A priori bounds and equicontinuity of orbits for the intermediate
  long wave equation.
\newblock {\em J. Evol. Equ.}, 26(3):Paper No. 81, 2026.

\bibitem{harropgriffiths2026homogenizationnonlinearschrodingerequation}
B.~Harrop-Griffiths and M.~Ntekoume.
\newblock Homogenization for the nonlinear {S}chr\"odinger equation with
  sprinkled nonlinearity.
\newblock \textit{Preprint}, 2026.
\newblock \href{https://arxiv.org/abs/2505.00589}{\texttt{arXiv:2505.00589}}.

\bibitem{hong2026longtimevaliditykortewegdevries}
Y.~Hong and J.~Jang.
\newblock Long-time validity of the {K}orteweg-de {V}ries approximation for the
  {B}oussinesq equation.
\newblock \textit{Preprint}, 2026.
\newblock \href{https://arxiv.org/abs/2609.30786}{\texttt{arXiv:2609.30786}}.

\bibitem{MR4767456}
Y.~Hong and C.~Yang.
\newblock On the {K}orteweg--de {V}ries limit for the {B}oussinesq equation.
\newblock {\em J. Differential Equations}, 408:94--116, 2024.

\bibitem{Kako}
M.~Kako and T.~Watanabe.
\newblock Nonlinear behaviors of ion-wave.
\newblock {\em Journal of the Physical Society of Japan}, 32(2):535--542, 1972.

\bibitem{MR604036}
V.~K. Kalantarov and O.~A. Lady\v{z}enskaja.
\newblock Formation of collapses in quasilinear equations of parabolic and
  hyperbolic types.
\newblock {\em Zap. Nau\v cn. Sem. Leningrad. Otdel. Mat. Inst. Steklov.
  (LOMI)}, 69:77--102, 274, 1977.

\bibitem{MR951744}
T.~Kato and G.~Ponce.
\newblock Commutator estimates and the {E}uler and {N}avier-{S}tokes equations.
\newblock {\em Comm. Pure Appl. Math.}, 41(7):891--907, 1988.

\bibitem{MR1101221}
C.~E. Kenig, G.~Ponce, and L.~Vega.
\newblock Oscillatory integrals and regularity of dispersive equations.
\newblock {\em Indiana Univ. Math. J.}, 40(1):33--69, 1991.

\bibitem{MR1086966}
C.~E. Kenig, G.~Ponce, and L.~Vega.
\newblock Well-posedness of the initial value problem for the {K}orteweg-de
  {V}ries equation.
\newblock {\em J. Amer. Math. Soc.}, 4(2):323--347, 1991.

\bibitem{MR1329387}
C.~E. Kenig, G.~Ponce, and L.~Vega.
\newblock A bilinear estimate with applications to the {K}d{V} equation.
\newblock {\em J. Amer. Math. Soc.}, 9(2):573--603, 1996.

\bibitem{MR4743514}
R.~Killip, T.~Laurens, and M.~Vi\c{s}an.
\newblock Sharp well-posedness for the {B}enjamin-{O}no equation.
\newblock {\em Invent. Math.}, 236(3):999--1054, 2024.

\bibitem{MR4922705}
R.~Killip, T.~Laurens, and M.~Vi\c{s}an.
\newblock Scaling-critical well-posedness for continuum {C}alogero-{M}oser
  models on the line.
\newblock {\em Commun. Am. Math. Soc.}, 5:284--320, 2025.

\bibitem{killip2026hamiltonianformulationcontinuumcalogeromoser}
R.~Killip, K.~Marsden, and M.~Vi\c{s}an.
\newblock The {H}amiltonian formulation of continuum {C}alogero--{M}oser
  models.
\newblock \textit{Preprint}, 2026.
\newblock \href{https://arxiv.org/abs/2604.09479}{\texttt{arXiv:2604.09479}}.

\bibitem{MR4628747}
R.~Killip, M.~Ntekoume, and M.~Vi\c{s}an.
\newblock On the well-posedness problem for the derivative nonlinear
  {S}chr\"odinger equation.
\newblock {\em Anal. PDE}, 16(5):1245--1270, 2023.

\bibitem{MR3990604}
R.~Killip and M.~Vi\c{s}an.
\newblock Kd{V} is well-posed in {$H^{-1}$}.
\newblock {\em Ann. of Math. (2)}, 190(1):249--305, 2019.

\bibitem{MR3820439}
R.~Killip, M.~Vi\c{s}an, and X.~Zhang.
\newblock Low regularity conservation laws for integrable {PDE}.
\newblock {\em Geom. Funct. Anal.}, 28(4):1062--1090, 2018.

\bibitem{MR2501679}
N.~Kishimoto.
\newblock Well-posedness of the {C}auchy problem for the {K}orteweg-de {V}ries
  equation at the critical regularity.
\newblock {\em Differential Integral Equations}, 22(5-6):447--464, 2009.

\bibitem{MR3016208}
N.~Kishimoto.
\newblock Sharp local well-posedness for the ``good'' {B}oussinesq equation.
\newblock {\em J. Differential Equations}, 254(6):2393--2433, 2013.

\bibitem{MR2654246}
N.~Kishimoto and K.~Tsugawa.
\newblock Local well-posedness for quadratic nonlinear {S}chr\"odinger
  equations and the ``good'' {B}oussinesq equation.
\newblock {\em Differential Integral Equations}, 23(5-6):463--493, 2010.

\bibitem{klaus}
F.~Klaus.
\newblock A priori estimates for a quadratic {dNLS}.
\newblock \textit{Preprint}, 2023.
\newblock \href{https://arxiv.org/abs/2305.11524}{\texttt{arXiv:2305.11524}}.

\bibitem{klaus2026wellposednesskdvhierarchy}
F.~Klaus, H.~Koch, and B.~Liu.
\newblock Wellposedness for the {K}d{V} hierarchy.
\newblock \textit{Preprint}, 2026.
\newblock \href{https://arxiv.org/abs/2309.12773}{\texttt{arXiv:2309.12773}}.

\bibitem{MR4515528}
F.~Klaus and R.~Schippa.
\newblock A priori estimates for the derivative nonlinear {S}chr\"odinger
  equation.
\newblock {\em Funkcial. Ekvac.}, 65(3):329--346, 2022.

\bibitem{MR3874652}
H.~Koch and D.~Tataru.
\newblock Conserved energies for the cubic nonlinear {S}chr\"odinger equation
  in one dimension.
\newblock {\em Duke Math. J.}, 167(17):3207--3313, 2018.

\bibitem{MR1251854}
F.~Linares.
\newblock Global existence of small solutions for a generalized {B}oussinesq
  equation.
\newblock {\em J. Differential Equations}, 106(2):257--293, 1993.

\bibitem{MR3032295}
S.~Oh and A.~Stefanov.
\newblock Improved local well-posedness for the periodic ``good'' {B}oussinesq
  equation.
\newblock {\em J. Differential Equations}, 254(10):4047--4065, 2013.

\bibitem{MR1048956}
R.~L. Sachs.
\newblock On the blow-up of certain solutions of the ``good'' {B}oussinesq
  equation.
\newblock {\em Appl. Anal.}, 36(3-4):145--152, 1990.

\bibitem{MR1628680}
G.~Schneider.
\newblock The long wave limit for a {B}oussinesq equation.
\newblock {\em SIAM J. Appl. Math.}, 58(4):1237--1245, 1998.

\bibitem{MR3702025}
G.~Schneider and H.~Uecker.
\newblock {\em Nonlinear {PDE}s}, volume 182 of {\em Graduate Studies in
  Mathematics}.
\newblock American Mathematical Society, Providence, RI, 2017.
\newblock A dynamical systems approach.

\bibitem{MSPC}
M.~Shan.
\newblock Endpoint well-posedness for the {B}oussinesq equation.
\newblock \textit{Preprint}.

\bibitem{MStalk}
M.~Shan.
\newblock Low regularity conservation laws for the {B}oussinesq equation.
\newblock One hour lecture delivered on August 17\textsuperscript{th}, 2022, at
  the ``{W}orkshop on {N}onlinear dispersive wave equations'' held at {P}eking
  {U}niversity, {C}hina.

\bibitem{MR2154153}
B.~Simon.
\newblock {\em Trace ideals and their applications}, volume 120 of {\em
  Mathematical Surveys and Monographs}.
\newblock American Mathematical Society, Providence, RI, second edition, 2005.

\bibitem{MR4270277}
B.~Talbut.
\newblock Low regularity conservation laws for the {B}enjamin-{O}no equation.
\newblock {\em Math. Res. Lett.}, 28(3):889--905, 2021.

\bibitem{MR4380591}
B.~A.~F. Talbut.
\newblock {\em Benjamin-{O}no at {L}ow {R}egularity: {A}n {I}ntegrability
  {A}pproach}.
\newblock ProQuest LLC, Ann Arbor, MI, 2021.
\newblock Thesis (Ph.D.)--University of California, Los Angeles.

\bibitem{MR2233925}
T.~Tao.
\newblock {\em Nonlinear dispersive equations}, volume 106 of {\em CBMS
  Regional Conference Series in Mathematics}.
\newblock Conference Board of the Mathematical Sciences, Washington, DC; by the
  American Mathematical Society, Providence, RI, 2006.
\newblock Local and global analysis.

\bibitem{MR1095753}
M.~Tsutsumi and T.~Matahashi.
\newblock On the {C}auchy problem for the {B}oussinesq type equation.
\newblock {\em Math. Japon.}, 36(2):371--379, 1991.

\bibitem{wang2025longtimeasymptoticsgoodboussinesq}
D.-S. Wang and X.~Zhu.
\newblock Long-time asymptotics of the good {B}oussinesq equation and its
  modified version: {P}ainlev\'{e} region.
\newblock \textit{Preprint}, 2025.
\newblock \href{https://arxiv.org/abs/2511.13382}{\texttt{arXiv:2511.13382}}.

\bibitem{MR1968371}
C.~E. Wayne and J.~D. Wright.
\newblock Higher order modulation equations for a {B}oussinesq equation.
\newblock {\em SIAM J. Appl. Dyn. Syst.}, 1(2):271--302, 2002.

\bibitem{zakharov}
V.~E. {Zakharov}.
\newblock {On Stochastization of one-dimensional chains of nonlinear
  oscillators}.
\newblock {\em Soviet Phys. JETP}, 38:108--110, 1974.

\end{thebibliography}
\end{document}